\pdfoutput=1
\documentclass[10pt]{amsart}
\usepackage{amsmath,amsfonts,amsthm,amssymb,graphicx,mathtools,bm,tikz,tikz-cd,epigraph, enumitem, xcolor,url,mathrsfs,stackrel, wasysym}
\usepackage[
colorlinks=true,
linkcolor = black,
citecolor= black
]{hyperref}
\usepackage[
noabbrev,nameinlink,capitalize
]{cleveref}
\usepackage[backend=bibtex,style=alphabetic,sorting=nyt]{biblatex}

\theoremstyle{definition}
\newtheorem{theorem}{Theorem}[section]

\newtheorem{lemma}[theorem]{Lemma}
\newtheorem{definition}[theorem]{Definition}

\newtheorem{proposition}[theorem]{Proposition}
\newtheorem{remarklabeled}[theorem]{Remark}
\newtheorem{corollary}[theorem]{Corollary}
\newtheorem{example}[theorem]{Example}

\numberwithin{equation}{section}
\DeclareMathOperator{\GL}{GL}
\DeclareMathOperator{\sheafO}{\mathcal{O}}

\DeclareMathOperator{\Sp}{Sp}
\DeclareMathOperator{\tr}{tr}
\DeclareMathOperator{\id}{id}
\DeclareMathOperator{\dR}{dR}
\DeclareMathOperator{\gr}{gr}
\DeclareMathOperator{\ad}{ad}
\DeclareMathOperator{\Ad}{Ad}
\DeclareMathOperator{\CP}{CP}
\DeclareMathOperator{\eu}{eu}
\DeclareMathOperator{\nil}{nil}
\DeclareMathOperator{\diag}{diag}

\DeclareMathOperator{\codim}{codim}
\DeclareMathOperator{\Spec}{Spec}

\DeclareMathOperator{\Hom}{Hom}

\DeclareMathOperator{\Span}{Span}
\DeclareMathOperator{\rk}{rk}
\DeclareMathOperator{\End}{End}
\DeclareMathOperator{\Loc}{Loc}

\DeclareMathOperator{\red}{red}

\DeclareMathOperator{\Eu}{Eu}
\DeclareMathOperator{\HC}{HC}
\DeclareMathOperator{\sym}{sym}

\DeclareMathOperator{\Stab}{Stab}
\DeclareMathOperator{\Der}{Der}

\newcommand{\inverse}{^{-1}}
\newcommand{\gitquo}{/\!/}
\newcommand{\idealm}{\mathfrak{m}}
\newcommand{\sliceA}{\underline{\A}}
\newcommand{\sliceR}{\underline{R}}
\newcommand{\sliceQ}{\underline{Q}}
\newcommand{\hamiltonianreduction}{/\!/\!/}
\newcommand{\RR}{\mathbb{R}}
\newcommand{\g}{\mathfrak{g}}
\newcommand{\C}{\mathbb{C}}
\newcommand{\idealI}{\mathcal{I}}
\newcommand{\D}{\mathcal{D}}
 
\newcommand{\ZZ}{\mathbb{Z}}

\newcommand{\VV}{\mathbf{V}}
\newcommand{\WW}{\mathbf{W}}
\newcommand{\PP}{\mathbb{P}}
\newcommand{\A}{\mathcal{A}}

\newcommand{\B}{\mathcal{B}}
\newcommand{\F}{\mathcal{F}}

\newcommand{\cartanh}{\mathfrak{h}}

\newcommand{\M}{\mathcal{M}}
\renewcommand{\L}{\mathcal{L}}
\newcommand{\mf}{\mathfrak}

\newcommand{\wt}{\widetilde}
\newcommand{\completeotimes}{\widehat{\otimes}}
\newcommand{\gl}{\mathfrak{gl}}

\renewcommand{\sl}{\mathfrak{sl}}

\newcommand{\WeylalgebraA}{\mathbf{A}}

\newcommand{\<}{\langle}
\renewcommand{\>}{\rangle}

\renewcommand{\d}{\partial}

\title[HC bimodules over flower quiver varieties]{Harish-Chandra bimodules over quantized flower quiver varieties with minimal support}
\author{Yaochen Wu}

\begin{document}
	\begin{abstract}
		We classify Harish-Chandra bimodules over the quantized flower quiver varieties with minimal support. We show that if the dimension vector is $n$, then there are $n!$ minimally supported simple Harish-Chandra bimodules for integral quantization parameters that are large enough or small enough, and there are none for other quantization parameters. The main tool used for this classification is an enhanced version of restriction functor. 
	\end{abstract}
	\maketitle
	
	\setcounter{tocdepth}{2}
	\tableofcontents
	\section{Introduction}
	\subsection{Quiver varieties}\label{quiver varieties introduction}
	We recall the notion of quiver varieties, following \cite{nakajima1998quiver}, with minor modifications. 
	Let $Q$ be a quiver which may contain multiple edges and edge loops, $Q_0$ be the set of vertices and $Q_1$ the set of arrows. 
	For $a\in Q_1$ let $t(a),h(a)$ denote the tail and head of $a$ respectively. 
	For $i,j\in Q_0$ let $n_{ij}$ denote the number of arrows (regardless of orientation) between $i,j$. 
	To each $i\in Q_0$ we associate a simple root $\alpha_i$. We identify $\ZZ_{\ge 0 }^{Q_0}$ with the root lattice of $Q$ and define the Tits form on $\ZZ_{\ge 0}^{Q_0 }$ by 
	\begin{equation} \ZZ_{\ge 0}^{Q_0 } \times \ZZ_{\ge 0}^{Q_0} \to \ZZ, (v,v') =\sum_{i\in Q_0} v_i v'_i (2- 2n_{ii} )- \sum_{i\neq j\in Q_0} v_iv'_j n_{ij}\end{equation}
	where for each $i\in Q_0$, $v_i,v'_i$ are the $i$-components of $v,v'$ respectively. 
	
	Let $v\in \ZZ_{\ge 0}^{Q_0 }$. Following \cite{crawley2001geometry}, we define
	\begin{equation}\label{definition of form p}
		p(v): = 1-\frac{1}{2}(v,v).
	\end{equation}
	
	Let $v,w\in\ZZ_{\ge 0}^{Q_0}$; let $V_i,W_i$ be vector spaces, $\dim V_i = v_i, \dim W_i = w_i$. Define the coframed representation space 
	\begin{equation}{\label{definition of R}}
		R(Q,v,w) := \bigoplus_{\substack{a\in Q_1,\\ t(a)\neq h(a)}} \Hom (V_{t(a)}, V_{h(a)}) \oplus
		\bigoplus_{\substack{a\in Q_1,\\ t(a)= h(a)}} \sl (V_{t(a)})
		\oplus \bigoplus_{i\in Q_0} \Hom (V_i,W_i).
	\end{equation}
	We emphasize that we attach \textit{traceless} endomorphisms to edge loops. 
	We write $R$ for $R(Q,v,w)$ when $Q,v,w$ is clear from the context. 
	
	The group $G = \GL(v) := \prod_{i\in Q_0} \GL(V_i)$ acts on $R$ naturally and induces a Hamiltonian $G$-action on $T^*R$. 
	Let $(x,y,p,q) = (x_a,y_a,p_i,q_i)_{i\in Q_0, a\in Q_1} \in T^*R$ where 
	\begin{equation}\label{linear maps listed in an element of T*R}
		\begin{split}
			x_a \in \Hom (V_{t(a)}, V_{h(a)}), y_a \in \Hom (V_{h(a)}, V_{t(a)}), \text{ if } t(a) \neq h(a) \\
			x_a \in \sl(V_{t(a)}), y_a \in  \sl(V_{t(a)}), \text{ if } t(a) = h(a)\\ 
			q_i \in \Hom (V_i,W_i), p_i\in \Hom(W_i,V_i).    
		\end{split}
	\end{equation} 
	Upon identifying $\g$ and $\g^*$ via the trace pairing, the moment map $\mu: T^*R \to \g^*$ for the Hamiltonian $G$-action on $T^*R$ is given by
	\begin{equation}
		\mu((x_a,y_a,p_i,q_i)) = [x,y]-pq.
	\end{equation} 
	More precisely,
	the $i$-component of the right hand side is
	\begin{equation}
		\sum_{h(a) = i} x_a y_a - \sum_{t(a) = i} y_a x_a -p_iq_i. 
	\end{equation}
	
	Let $\theta$ be a character of $G$; it has the form $\theta (g) = \prod_{i\in Q_0} \det(g_i)^{\theta_i}$, $g_i \in \GL(v_i)$, $\theta_i \in\ZZ$. We can thus identify $\ZZ^{Q_0}$ with the character lattice of $G$ and write $\theta \in \ZZ^{Q_0}$. Similarly, we identify $(\g/[\g,\g])^*$ with $\C^{Q_0}$. 
	
	Let $ (T^*R)^{\theta-ss}$ denote the set of semistable points of $T^*R$. 
	Let $\lambda \in (\g/[\g,\g])^* $. 
	The GIT quotient
	\begin{equation}\M_\lambda^\theta (v,w) : = (\mu\inverse(\lambda)^{\theta-ss})\gitquo G\end{equation}
	is called a Nakajima quiver variety. 
	The variety $\M_\lambda^0 (v,w)$ is affine, and there is a projective morphism $\rho: \M_\lambda^\theta (v,w) \to \M_\lambda^0 (v,w)$. 
	For $x\in \M_\lambda^\theta(v,w)$, we can write take a representative in $T^*R$ and write
	\begin{equation}\label{representative of an element of quiver variety in T*R}
		x = [x_a,y_a,p_i,q_i]_{a\in Q_1,i\in Q_0}
	\end{equation}
	where $x_a,y_a,p_i,q_i$ are as in \eqref{linear maps listed in an element of T*R}. 
	
	For $\theta\in\ZZ^{Q_0}, \lambda\in \C^{Q_0}$, the pair $(\theta,\lambda)$ is said to be \textit{generic} if there are no positive roots $v'<v$ of $Q$ such that $\lambda\cdot v'=\theta\cdot v'= 0$. 
	We say $\lambda$ is generic if $(0,\lambda)$ is generic and we say $\theta$ is generic if $(\theta,0)$ is generic.
	When $(\theta,\lambda)$ is generic,  $\M_\lambda^\theta (v,w)$ is smooth. 
	
	The algebra $\C[\M^0_0(v,w)]$ has a natural Poisson bracket inherited from the Poisson structure on $\C[T^*R]$. 
	It is equipped with a contracting $\C^*$-action, also inherited from that on $T^*R$, and the Poisson bracket has degree -2 with respect to this $\C^*$-action. Thus,  $\C[\M^0_0(v,w)]$ is a graded Poisson algebra. 
	
	It is shown in \cite[Theorem 1.2]{bellamy2021symplectic} that affine quiver varieties have symplectic singularities in the sense of \cite{beauville2000symplectic}. Therefore, by \cite[Theorem 2.3]{kaledin2006symplectic}, $\M_\lambda^0(v,w)$ has finitely many locally closed symplectic leaves. 
	We will describe these leaves in \Cref{subsection Symplectic leaves and local structures}.
	
	For simplicity of notations we write $\M_\lambda^\theta$ for $\M_\lambda^\theta(v,w)$ when $Q, v,w$ are clear from the context. We also write $\M$ for $\M_0^0$. 
	
	\subsection{Quantizations}\label{introduction section quantization}
	We recall the notion of quantizations of smooth quiver varieties and affine quiver varieties via quantum Hamiltonian reduction. More details on quantizations will be given in \Cref{more on quantizations subsection}. 
	
	Let $A$ be a graded, Noetherian, Poisson algebra such that the Poisson bracket has degree $-d$, $d\ge 0$. For example, we can take $A = \C[\M]$ and $d=2$. 
	\begin{definition}
		A filtered quantization of $A$ is a pair $(\A,\iota)$, where 
		\begin{itemize}
			\item $\A$ is an associative unital algebra equipped with a filtration $\A = \bigcup_{i\ge 0} \A_{\le i}$ such that $[\A_{\le i},\A_{\le j}] \subset \A_{\le i+j-d}$; 
			\item $\iota: \gr \A \to A$ is an isomorphism of graded Poisson algebras, where the Poisson bracket on $\gr \A$ is defined by 
			$\{a+\A_{\le i-1}, b+\A_{\le j-1}\} = [a,b] + \A_{\le i+j-d-1}$. 
		\end{itemize}
	\end{definition}
	\begin{definition}
		A graded formal quantization of $A$ is a pair $(\A_\hbar,\iota)$, where
		\begin{itemize}
			\item $\A_\hbar$ is an associative unital flat $\C[\hbar]$ algebra, such that $[\A_\hbar,\A_\hbar] \subset \hbar^d \A_\hbar$. 
			\item $\A_\hbar$ is equipped with a grading where $\hbar$ has degree 1. 
			\item $\iota: \A_\hbar / \hbar^d \A_\hbar \to A$ is an isomorphism of graded Poisson algebras, where the Poisson bracket on $\A_\hbar/\hbar^d \A_\hbar$ is defined by 
			\begin{equation}\label{definition of poisson bracket for formal graded quantization}
				\{a+\hbar^d \A_\hbar ,b+\hbar^d \A_\hbar\} = \frac{1}{\hbar^d}[a,b]+ \hbar^d \A_\hbar.
			\end{equation}
		\end{itemize}
	\end{definition}
	Given a filtered quantization $\A$ of $A$, we can take its Rees algebra $\A_\hbar = \bigoplus_{i\ge 0} \hbar^i \A_{\le i}$. Then $\A_\hbar$ is a graded formal quantization of $A$. 
	
	Let $A$ be as above, and let $I\subset A$ be a graded ideal. Let $A^{\wedge I} = \varprojlim_n A/I^n$ be the completion. 
	\begin{definition}
		A formal quantization of $A^{\wedge I}$ is a pair $(\A_\hbar^\wedge, \iota)$, where
		\begin{itemize}
			\item $\A_\hbar^\wedge$ is a flat $\C[[\hbar]]$-algebra 
			such that $[\A_\hbar^\wedge,\A_\hbar^\wedge] \subset \hbar^d \A_\hbar^\wedge$. 
			\item $\A_\hbar^\wedge$ is
			equipped with a $\C^*$-action, and $\hbar$ has degree 1, and the action on each $\A_\hbar^\wedge / (\hbar^k)$ is rational. 
			\item $\iota: \A_\hbar^\wedge / \hbar^d \A_\hbar^\wedge \to A^{\wedge I}$ is a $\C^*$-equivariant isomorphism of Poisson algebras, where the Poisson bracket on $A_\hbar^\wedge / \hbar^d \A_\hbar^\wedge $ is defined similarly to \eqref{definition of poisson bracket for formal graded quantization}, and the Poisson bracket on $A^{\wedge I}$ is extended from that on $A$ by continuity. 
			\item Let $\idealI_\hbar$ be the preimage of $I$ under the map $\A_\hbar^\wedge \twoheadrightarrow \A_\hbar^\wedge / \hbar^d \A_\hbar^\wedge \xrightarrow{\iota} A^{\wedge I}$. We require $\A_\hbar^\wedge$ is complete and separated with respect to the $\idealI_\hbar$-adic topology. 
		\end{itemize}
	\end{definition}
	For example, if $(\A_\hbar, \iota)$ is a graded formal quantization of $A$, then let
	\begin{equation}
		\A_{\hbar}^{\wedge \idealI_\hbar} := \varprojlim_n \A_\hbar / \idealI_\hbar^n.
	\end{equation}
	It is straightforward to verify that $\iota$ induces an isomorphism $\iota^\wedge: \A_{\hbar}^{\wedge \idealI_\hbar} / \hbar^d \to A^{\wedge I}$, which is $\C^*$-equivariant and Poisson. Thus, $(\A_\hbar^{\wedge \idealI_\hbar},\iota^\wedge)$ is a formal quantization of $A^{\wedge I}$.

	Note that the pair $(\Spec A/I, A^{\wedge{I}})$ defines an affine formal scheme. We also say $\A_\hbar^\wedge$ is a formal quantization of the affine formal scheme $(\Spec A/I, A^{\wedge I})$. 
	
	Let $X$ be a Poisson variety, i.e. $\sheafO_X$ has the structure of sheaf of Poisson algebras. Suppose $\sheafO_X$ is graded and the Poisson bracket has degree $-d$. 
	\begin{definition}
		A graded formal quantization of $\sheafO_X$ is a pair $(\D_\hbar,\iota)$, where
		\begin{itemize}
			\item $\D_\hbar$ is a sheaf on $X$ of flat $\C[[\hbar]]$-algebras, complete and seperated in the $\hbar$-adic topology; 
			\item $\D_\hbar$ is equipped with a $\C^*$-action, $\hbar$ has degree 1, and the action on each $\D_\hbar / \hbar^kD_\hbar$ is rational; 
			\item $[\D_\hbar,\D_\hbar] \subset \hbar^d \D_\hbar$; 
			\item $\iota : \D_\hbar / \hbar^d\D_\hbar \xrightarrow{\sim} \sheafO_X$ is an isomorphism of sheaves of Poisson algebras; the Poisson bracket on $ \D_\hbar /\hbar^d\D_\hbar$ is defined by 
			$\{a+\hbar^d\D_\hbar,b+\hbar^d\D_\hbar\} = \frac{1}{\hbar^d} [a, b]$, for $a,b\in \D_\hbar$.
		\end{itemize}
	\end{definition}

	Recall that $R = R(Q,v,w)$ is the representation space of a quiver, and let $D(R)$ be the algebra of differential operators on $R$. We can view it as a sheaf $D_R$ on $T^*R$ via microlocalization; that is, we only assign sections to $\C^*$-stable open subsets of $T^*R$. Both $D(R)$ and $D_R$ are equipped with the Bernstein filtration ($\deg R = \deg R^* = 1$). Then $D_R$ is a filtered quantization of the structure sheaf $\sheafO_{T^*R}$ with its natural Poisson structure, where $d = 2$.  
	
	There is a quantum comoment map $\Phi$ quantizing the usual comoment map $\mu^*$:
	\begin{equation}\label{not symmetrized qcm}
		\Phi: \g \to D(R), x\mapsto x_R
	\end{equation}
	where $x_R$ is the vector field on $R$ corresponding to $x$ induced by the $G$-action. We remark that $x_R \in D(R)_{\le 2}$. 
	
	Fix a generic $\theta$. Following \cite[Section 1.4]{bezrukavnikov2021etingof}, we define 
	\begin{equation}\A_\lambda^\theta := [D_R/D_R\{\Phi(x) - \<\lambda,x\>|x\in \g\}|_{(T^*R)^{\theta-ss}}]^G\end{equation}
	which is a sheaf of algebras on $\M_0^\theta$. Similarly, we define 
	\begin{equation}\label{definition of quantum algebra Alambda}
		\A_\lambda := [D(R)/D(R)\{\Phi(x) - \<\lambda,x\>|x\in \g\}]^G\end{equation}
	which is an algebra. 
	We define
	\begin{equation}
		\A_\lambda^0 := [D_R/D_R\{\Phi(x) - \<\lambda,x\>|x\in \g\}]^G
	\end{equation}
	which is a sheaf of algebras on $\M$ whose global sections coincides with $\A_\lambda $. Note that the notion here is slightly different from that of \Cite{bezrukavnikov2021etingof}.
	
	Each of $\A_\lambda^\theta, \A_\lambda$ and $\A_\lambda^0$ inherit a filtration from the Bernstein filtration of $D_R$, and we have 
	\begin{equation}
		\gr \A_\lambda^\theta = \sheafO_{\M^\theta_0} .
	\end{equation}
	Moreover, if the moment map $\mu$ is flat, then
	\begin{equation}
		\gr \A_\lambda = \C[\M^0_0], \gr \A_\lambda^0 = \sheafO_{\M_0^0}. 
	\end{equation}
	It is easy to check that the algebra $\A_\lambda$ is a filtered quantization of $\C[\M_0^0(v)]$.

	\subsection{Harish-Chandra and Poisson bimodules}
	Suppose $\A$ is a filtered algebra such that $\gr \A$ is a Noetherian Poisson algebra $A$, where the Poisson bracket has degree $-2$; for example, we can take $A = \C[\M_0^0]$ and $\A = \A_\lambda$. 
	\begin{definition}
		An $\A$-bimodule $\B$ is called a Harish-Chandra (HC) $\A$-bimodule if it admits a bimodule filtration $\B = \bigcup_{i\ge 0} \B_{\le i}$ such that $\gr \B$ is a finitely generated $A$-module. Such a filtration is called a \textit{good filtration}. 
	\end{definition}
	
	\begin{definition}
		An $A$-module $M$ is called a Poisson $A$-module if there is a bracket $\{\bullet,\bullet\}: A\otimes M\to M$ satisfying the Jacobi and Leibniz identities: 
		\begin{enumerate}
			\item $\{\{a_1,a_2\},m\} = \{a_1,\{a_2,m\}\} - \{a_2,\{a_1,m\}\}$;
			\item $\{a_1,a_2m \} = \{a_1,a_2\}m+a_2\{a_1,m\}$. 
		\end{enumerate}
		
	\end{definition}
	Assume $\Spec A$ is a symplectic singularity, so that it has finitely many symplectic leaves. 
	Let $\B$ be a HC $\A$-bimodule. Equip $\B$ with a good filtration and set $B = \gr \B$. Then the commutator map $[\bullet,\bullet]: \A \otimes \B \to \B$ induces a graded Poisson $A$-module structure on $B$. 
	\begin{definition}
		The \textit{support} of $\B$ is defined to be the support of $B$ in $\Spec A$.
	\end{definition}
	The support of $\B$ is a closed Poisson subvariety of $\Spec A$, and therefore is a union of symplectic leaves. Let $\L$ be a symplectic leaf of $\Spec A$. 
	\begin{definition}\label{definition of cats of HC bimodules with support condition}
		Let $\HC_{\overline{\L}}(\A)$ denote the full subcategory of $\HC (\A)$ whose objects are supported on $\overline{\L}$. Let $\HC_{\d {\L}}(\A)$ denote the Serre subcategory of $\HC_{\overline{\L}} (\A)$ whose objects are supported on $\partial\L$. Finally, let $\HC_\L(\A)$ denote the quotient category of $\HC_{\overline{\L}}(\A)$ by $\HC_{\d {\L}}(\A)$. 
	\end{definition} 
	
	\subsection{Goal and structure of the paper}
	We consider the flower quiver, i.e. quiver with one vertex and $\ell$ edge-loops. We assume $\ell \ge 2$. 
	The goal of this paper is to classify Harish-Chandra bimodules over the quantizations $\A_\lambda$ which have minimal support. The definition of ``minimal support" will be made clear in \Cref{Relevant leaves}. Our main result is the following. 
	\begin{theorem}\label{summary result in introduction}
		Suppose $\lambda \in \ZZ$. If $\lambda \ge (n-1)({\ell}-1)+w$, 
		or $\lambda\le -(n-1)({\ell}-1)$,  
		there are precisely $n!$ non-isomorphic minimally supported HC bimodules over $\A_\lambda$. 
		
		Otherwise, if $\lambda\not\in\ZZ$ or $ -(n-1)({\ell}-1) < \lambda < (n-1)({\ell}-1)+w$, there is none. 
	\end{theorem}
	
	Let us explain the strategy of the proof of \Cref{summary result in introduction}. The key tool is an enhanced restriction functor. 
	
	\subsubsection{Classical restriction functor} 
	Let us briefly recall the construction of classical restriction functor for quantized quiver varieties, following \Cite[Section 3.3]{bezrukavnikov2021etingof}.
	
	Let $x$ be a point in $\M$, and assume $x$ lies in a symplectic leaf $\L$ of $\M$.
	Then 
	\begin{equation}\label{formal nbhd of point in introduction}
		\C[\M]^{\wedge x} \cong \C[[R_0]] \completeotimes \C[\underline{\M}_0^0(\underline{v},\underline{w})]^{\wedge 0}
	\end{equation}
	where $R_0$ is the tangent space at $x$, and $\underline{\M}_0^0(\underline{v},\underline{w})$ is another quiver variety, to be recalled in \Cref{subsection Symplectic leaves and local structures}. We call it the slice quiver variety and denote it by $\underline{\M}$ for simplicity of notations. 
	
	This isomorphism has a quantum analogue. There is an affine linear map $r: \C^{Q_0}\to \C^{\underline{Q}_0}$ to be defined in \eqref{definition of symmetrized restriction}. For any $x\in \M$, let $\idealI_x$ be the preimage of the maximal ideal of $x$ under the projection $	(\A^0_\lambda)_\hbar \to \C[\M_0^0(v,w)]$, and let $^{\wedge x}$ denote the completion with respect to $\idealI_x$. 
	By {\cite[Lemma 3.7]{bezrukavnikov2021etingof}} we have the following isomorphism of algebras 
	\begin{equation}\label{decomposition of quantum algebra at a point}
		(\A_\lambda)_\hbar^{\wedge_x} \cong (\underline{\A}_{r(\lambda)})^{\wedge_0}_\hbar \completeotimes_{\C[[\hbar]]} \WeylalgebraA_\hbar^{\wedge_0}.
	\end{equation}
	Here, $\underline{\A} _{r(\lambda)}$ is the quantization of $\underline{\M} $ with parameter $r(\lambda)$, and  $\WeylalgebraA_\hbar$ is the homogeneous Weyl algebra of the symplectic vector space $R_0$. 
	
	Let $\B$ be a Harish-Chandra bimodule over $\A_\lambda$; pick a good filtration of $\B$ and form the Rees bimodule $\B_\hbar$. 
	By \Cite[Section 3.3.3]{bezrukavnikov2021etingof},
	we have a decomposition 
	\begin{equation}\label{decomposition classical restriction functor}
		(\B_\hbar)^{\wedge x} \cong \underline{\B}' \completeotimes_{\C[[\hbar]]} \WeylalgebraA_\hbar^{\wedge_0}
	\end{equation}
	where $\underline{\B}'$ is a finitely generated $(\underline{\A}_{r(\lambda)})_\hbar^{\wedge_0}$-bimodule equipped with an Euler derivation $\hat{\Eu}$ that is compatible with the Euler derivation $\hbar\partial_\hbar$ on $(\underline{\A} _{r(\lambda)})_\hbar^{\wedge_0}$. Let $\underline{\B}_\hbar$ be the $\hat{\Eu}$-finite part of $\underline{\B}'_\hbar$. Define $\B_{\dagger,x}:= \underline{\B}_\hbar/(\hbar-1)$; this is a Harish-Chandra bimodule over $\underline{\A} _{r(\lambda)}$. 
	
	\begin{theorem}[{\cite[Lemma 3.9]{bezrukavnikov2021etingof}}]{\label{restriction functor and support property}}
		The functor $\bullet_{\dagger,x}: \HC(\A^0_\lambda) \to \HC(\underline{\A} _{r(\lambda)})$ satisfies the following properties. 
		\begin{enumerate}
			\item $\bullet_{\dagger,x}$ is exact. 
			\item If $\B\in \HC(\A^0_\lambda)$, then the associated variety $V(\B_{\dagger,x})$ in $\underline{\M}_0^0$ of $\B_{\dagger,x}$ is uniquely characterized by 
			\begin{equation}(V(\B_{\dagger,x})\times \L )^{\wedge_x} = V(\B)^{\wedge_x}\end{equation}
			where $\L$ is the symplectic leaf of $\M_0^0$ that contains $x$. 
		\end{enumerate}
		
	\end{theorem}
	
	We will generalize this construction by taking the completion not at a point $x$, but along (open subsets of) a symplectic leaf. 
	
	\subsubsection{Relevant leaves and slices}
	First, in \Cref{Relevant leaves}, we give criteria for a symplectic leaf to be the closure of the support of a Harish-Chandra bimodule. We call such leaves \textit{relevant}. In particular, we will describe the minimal relevant leave, and will focus only on this leaf. 
	
	In the bookkeeping \Cref{section slice to minimal leaf}, we collect facts about the slice quiver variety (and its quantizations) associated to the minimal relevant leaf. We also introduce certain group actions on the slice quiver variety that will be important later. 
	\subsubsection{Formal neighbourhood of minimal leaf}
	In the geometric technical heart \Cref{Formal neighbourhood of minimal leaf}, we give a description of the formal neighbourhood of affine open subsets of the minimal relevant leaf, i.e. we obtain an analogue of \eqref{formal nbhd of point in introduction}.
	Instead of completing along $\L$, we need to take affine open subsets $\L_K \subset \L$ (to be defined in \Cref{subsection good affine open cover}) and affine open subsets $\M_K$ of $\M_0^0(v,w)$, so that $\L_K\subset \M_K$ is a closed subset.
	We will describe $\M_K^{\L_K}$. 
	
	First, we consider the formal neighbourhood before reduction, i.e. we pick open subsets $(T^*R)_K$ of $T^*R$ so that $\mu\inverse(0)\cap (T^*R)_K \gitquo G = \M_K$. 
	Let $H\subset G$ be the stabilizer of a point in $\mu\inverse(0)$ with closed $G$-orbit that maps to a point in $\L$; it turns out that $H \cong (\C^*)^n$ for the minimal leaf. In \Cref{subsection model varieties}, following the ideas of \Cite{losev2006symplectic}, we construct a ``model variety"
	\begin{equation}
		X\cong G\times^H (U\oplus V) \times \tilde{\L}.
	\end{equation} 
	Replacing $\tilde{\L}$ with $\tilde{\L}_K$, we get affine opens $X_K\subset X$. 
	Here $U\cong (\g/\cartanh)^*$, $V$ is identified with the representation space of the slice quiver, $\tilde{\L}$ is the universal cover of $\L$ with Galois group $S_n$, and $\tilde{\L}_K$ is the preimage of $\L_K$. There is a natural closed embedding $\tilde{\L}_K \hookrightarrow (T^*R)_K$. The model variety $X_K$ carries a free action of $S_n$. 
	The main part of \Cref{subsection model varieties} is devoted to an explicitly construction of isomorphisms of formal schemes 
	\begin{equation}
		\phi_K^\wedge: X_K^{\wedge \tilde{\L}_K}/S_n\to (T^*R)_K^{\wedge \tilde{\L}_K}. 
	\end{equation}
	
	We will pick two affine opens $\L_I$ and $\L_J$ such that $\codim_{\L} \L\setminus(\L_I\cup \L_J) = 2$, and construct $\phi_J^\wedge$ and $\phi_J^\wedge$. 
	In \Cref{Unipotency of transition map}, We show that the composition $(\phi_J^\wedge)\inverse\circ \phi_I^\wedge $, satisfies some (pro-)unipotency condition, so that we can take its $\log$ and study the corresponding vector field. 
	
	The isomorphisms $\phi_K^\wedge$ do not intertwine the symplectic forms on the source and target. 
	Therefore, in \Cref{subsection Moser's theorem}, we generalize the ideas of \Cite{losev2006symplectic} and use a Moser's theorem argument to produce an automorphism $\Phi_I$ of $X_I^\wedge$, so that $\phi_I^\wedge \circ \Phi_I^\wedge: X_K^{\wedge \tilde{\L}_I}/S_n\to (T^*R)_I^{\wedge \tilde{\L}_I}$ intertwines the symplectic forms; similarly we construct an automorphism $\Phi_J$ of $X_J^\wedge$. 
	
	By taking Hamiltonian reductions of the symplectomorphism $\phi_I^\wedge \circ \Phi_I^\wedge$, we obtain an isomorphism of Poisson formal schemes (see (1) of \Cref{structure of nbhd after reduction})
	\begin{equation}
		\theta_I^{\red}: (\tilde{\L}_I \times \underline{\M}^{\wedge 0})/S_n \xrightarrow{\sim} \M_I^{\wedge \L_I} 
	\end{equation}
	This is our analogue of \eqref{formal nbhd of point in introduction}. 
	
	We will show that the composition $(\theta^{\red}_J) \inverse\circ\theta^{\red}_I$ is ``inner" in the following sense. 
	There is a function $f\in (\C[\tilde{\L}_{I\cap J}]\completeotimes \C[\underline{\M}]^{\wedge_0})^{S_n}$ such that, as algebra automorphisms of $(\C[\tilde{\L}_{I\cap J}]\completeotimes \C[\underline{\M}]^{\wedge_0})^{S_n}$,
	\begin{equation}\label{poisson inner derivation in introduction}
		(\theta^{\red}_J|_{I\cap J}\inverse\circ\theta^{\red}_I|_{I\cap J})^* = \exp(\{f,\bullet\}).
	\end{equation}
	Moreover, $\{f,-\}$ is a pro-nilpotent operator. 
	Sections \ref{Unipotency of transition map} through \ref{inner derivation} are devoted to establishing these properties. The function $f$ and its quantum analogue will be important in \Cref{section Minimally-supported HC bimodules}. 
	
	\subsubsection{Bimodules over quantized algebras}\label{sbsbsection introduction Bimodules over quantized algebras}
	
	In \Cref{subsection Quantizations of affine opens} and \ref{subsection Lifting automorphisms}, we obtain an analogue of the decomposition \eqref{decomposition of quantum algebra at a point}. 
	More precisely, let $\M_K, \L_K$ be as in \Cref{subsection good affine open cover}. Let $\idealI_K$ be the preimage of the ideal of $\L_K\subset \M_K$ under the projection $\A_{\lambda,\hbar}|_{\M_K} \xrightarrow{/ \hbar} \C[\M_K]$, and let $^{\wedge\L_K}$ denote the $\idealI_K$-adic completion. 
	We will construct isomorphisms of Poisson algebras
	\begin{equation}
		\vartheta_{K,\hbar}^{\red}:
		\A_{\lambda,\hbar}^{\wedge \L_K} \xrightarrow{\sim} ( \WeylalgebraA_\hbar(\tilde{\L}_K) \completeotimes \underline{\A}_{r(\lambda), \hbar}^{\wedge_0})^{S_n}.
	\end{equation}
	Here, $\WeylalgebraA_\hbar(\tilde{\L}_K)$ is the restriction of a homogeneous Weyl algebra to $\tilde{\L}_K$. 
	Moreover, for the special choice of $\L_I,\L_J$, there is some $F$ lifting $f$ in \eqref{poisson inner derivation in introduction} such that
	\begin{equation}\label{quantum inner automorphism in introduction}
		\vartheta_{I,\hbar}^{\red} \circ  (\vartheta_{J,\hbar}^{\red} )\inverse = \exp \left( \frac{1}{\hbar^2}\ad F\right).
	\end{equation}
	The construction of $\vartheta_{K,\hbar}^{\red}$ and $F$ involves quantum versions of the results in \Cref{Formal neighbourhood of minimal leaf}. 
	
	We will define various categories and functors analogous to the construction of $(-)_{\dagger,x}$. 
	We very briefly state the steps of this construction; \Cref{subsection on coherent poisson bimodule over algebra} through \ref{subsection various functors} are devoted to the rigorous formulations. 
	
	Let $\B \in \HC(\A_\lambda^0)$; take a good filtration of $\B$ and form its Rees module $\B_\hbar$; 
	restrict it to an affine open $\M_K$ defined in \Cref{Formal neighbourhood of minimal leaf}, 
	and take its completion along the closed subset $\L_K$ of $\M_K$; we get an $\A_{\lambda,\hbar}^{\wedge \L_K}$-bimodule $\B_{\hbar,K}^\wedge$. 
	Then, use the isomorphisms $\vartheta_{K,\hbar}^{\red}$, we equip $\B_{\hbar,K}^\wedge$ with an $(\WeylalgebraA_\hbar(\tilde{\L}_K) \completeotimes \underline{\A}_{r(\lambda), \hbar}^{\wedge_0})^{S_n}$-bimodule structure. 
	Take $K= I,J$; we use the function $F$ in \eqref{quantum inner automorphism in introduction} to ``glue" $\B_{\hbar,I}^\wedge$ and $\B_{\hbar,J}^\wedge$ into a bimodule over $( \WeylalgebraA_\hbar(\tilde{\L}_I\cup \tilde{\L}_J) \completeotimes \underline{\A}_{r(\lambda), \hbar}^{\wedge_0})^{S_n}$; let us denote the resulting bimodule by $\mathcal{N}$. By equivariant descent we get an $S_n$-bimodule $\tilde{\mathcal{N}}$ over $\WeylalgebraA_\hbar(\tilde{\L}_I\cup \tilde{\L}_J) \completeotimes \underline{\A}_{r(\lambda), \hbar}^{\wedge_0}$. The condition $\codim_\L (\L \setminus (\L_I\cup\L_J)) = 2$ guarantees that
	\begin{equation}
		\tilde{\mathcal{N}} \cong \WeylalgebraA_\hbar(\tilde{\L}_I\cup \tilde{\L}_J) \otimes Fl(\tilde{\mathcal{N}})
	\end{equation}
	where $Fl(\tilde{\mathcal{N}})\subset \tilde{\mathcal{N}}$ is the Poisson centralizer of $\WeylalgebraA_\hbar(\tilde{\L}_I\cup \tilde{\L}_J)$ in $\tilde{\mathcal{N}}$. This is our analogue of \eqref{decomposition classical restriction functor}. Take $\C^*$-finite part of $Fl(\tilde{\mathcal{N}})$ and quotient by $\hbar-1$, we get an object of $\HC(\underline{\A}_{r(\lambda)})$. It turns out that we get a bijection between the set of isomorphism classes of simple Harish-Chandra bimodules over $\A_\lambda$ with minimal support, and the set of isomorphism classes of $S_n$-equivariant, finite dimensional $\underline{\A}_{r(\lambda)}$-bimodules.  
	
	\subsubsection{Finite dimensional modules over quantized slice algebra}
	The slice quiver varieties are hypertoric varieties, and their quantizations are hypertoric enveloping algebras, in the sense of \Cite{braden2012hypertoric}. There is a combinatorial description of finite dimensional modules over such algebras. We recall these results in \Cref{Hypertoric quantum slice algebras}, and classify finite dimensional simple modules over the the quantum slice algebra. This will complete our original classification problem. 
	\subsection*{Acknowledgement} 
	I am very grateful to Ivan Losev for suggesting this problem, many fruitful discussions and numerous remarks that helped improve the exposition. I would like to thank Do Kien Hoang, Mengwei Hu, Dmytro Matvieievskyi and Trung Vu for useful discussions.
	This work is partially supported by the NSF under grant DMS-2001139.
	
	\section{Preliminaries}\label{More on quiver varieties}
	\subsection{Extended quiver}\label{Framed vs nonframed quiver varieties}
	We can view any quiver variety $\M_\lambda^\theta(v,w)$ or $\M_\lambda^0(v,w)$ with a framing as one without framing. Namely, let $Q,v,w, \lambda,\theta$ be as before. Let ${Q^\infty}$ be the \textit{extended quiver} defined as follows. Its set of vertices is  ${Q^\infty_0} = Q_0 \cup \{\infty\}$. The arrows between vertices in $Q_0$ are the same as those in $Q_1$. In addition, for each $i\in Q_0$, there are $w_i$ arrows from the vertex $i$ to the vertex $\infty$. 
	
	We denote the new simple root associated to the vertex $\infty$ by $\alpha_\infty$. Define the extended dimension vector $\tilde{v}\in \ZZ_{\ge 0}^{{Q^\infty_0}}$ by $\tilde{v}_\infty = 1$ and $\tilde{v}_i = v_i$ for $i\in Q_0$. 
	
	It is clear that $R(Q,v,w) = R({Q^\infty},\tilde{v})$, so  $\M^\theta_\lambda(Q,v,w) = \M_{\tilde{\lambda}}^{\tilde{\theta}}({Q^\infty},\tilde{v})$, 
	where $\tilde{\lambda}_\infty = -v\cdot \lambda$, $\tilde{\lambda}_i = \lambda_i$ for $i\in Q_0$; $\tilde{\theta}_\infty = -\theta \cdot v, \tilde{\theta}_i = \theta_i$ for $i \in Q_0$. 
	\subsection{Symplectic leaves and local structures}\label{subsection Symplectic leaves and local structures}
	An affine quiver variety has finitely many symplectic leaves. They are classified by \textit{representation types}, which we recall below. 
	Let $Q,v,w,\lambda$ be as before, and let $x\in\M_\lambda^0$. 
	Any representative of $x$ in $T^*R$ consists of a collection of linear maps as in \eqref{linear maps listed in an element of T*R}. Let $V_x = \bigoplus_{i\in Q_0 } V_i \oplus \bigoplus_{i\in Q_0 } W_i$. 
	\begin{definition}
		Suppose $V_x =S_0 \oplus S_1^{\oplus n_1} \oplus ... \oplus S_k^{\oplus n_k}$, 
		where $S_i$'s are pairwise non-isomorphic simple representations of $Q^\infty$, and $\dim(S_0)_\infty = 1$. Write $v^i = \dim S_i \in \ZZ_{\ge 0}^{{Q^\infty_0}}$. Then we say 
		\begin{equation}\label{general form of a representation type}
			\tau = (v^0,1;v^1,n_1;v^2,n_2;...;v^k,n_k)
		\end{equation}
		is the \textit{representation type} of $x$. 
	\end{definition}
	The stabilizer of an element of $T^*R$ with representation type $\tau$ in \eqref{general form of a representation type} is 
	\begin{equation}\label{definition of stabilizer H}
		H = \prod_{i=1}^k \GL(n_i)
	\end{equation} 
	By \cite[Section 3.iv]{nakajima1998quiver}, a symplectic leaf of $\M_\lambda^0$ constitutes of all elements with a fixed representation type. In what follows, we will often say a symplectic leaf is associated to a representation type $\tau$, and denote the leaf by $\L_\tau$. 
	
	The following results of Crawley-Boevey characterizes the dimension vectors that can appear in a representation type. 
	Recall the function $p$ defined in \eqref{definition of form p}. 
	\begin{theorem}[{\cite[Theorem 1.2]{crawley2001geometry}}]{\label{CB simple dimension criteria}}
		The following conditions are equivalent.
		\begin{enumerate}
			\item There is a simple representation of $Q^\infty$ in  $\mu\inverse(\tilde{{\lambda}})$ with dimension vector $v'$.
			\item $v'$ is a positive root of ${Q^\infty}$, $\tilde{{\lambda}}\cdot v' = 0$, and for any decomposition $v' = \beta^1+ ... + \beta^n$ where $n\ge 2$ and $\beta^i$ are positive roots of ${Q^\infty}$ such that $\beta^i\cdot \tilde{{\lambda}}=0$, we have $p(v') > \sum_{i=1}^n p(\beta^i)$. 
		\end{enumerate}
	\end{theorem}
	A similar combinatorial criterion tells us about the flatness of the moment map $\mu$. 
	\begin{theorem}[{\cite[Theorem 1.1]{crawley2001geometry}}]{\label{CB flatness of moment map}}
		Fix a dimension vector $v$ and a framing $w$. The following conditions are equivalent.
		\begin{enumerate}
			\item The moment map $\mu: T^*R(Q^\infty, \tilde{v}) \to \g^*$ is flat. 
			\item For any decomposition $\tilde{v} = \beta^1+ ... + \beta^n$ where $n\ge 2$ and $\beta^i$ are positive roots of ${Q^\infty}$, we have $p(v') \ge \sum_{i=1}^n p(\beta^i)$. 
		\end{enumerate}
	\end{theorem}
	
	\begin{theorem}[{\cite[Theorem 1.3]{crawley2001geometry}}]{\label{CB strata dimension formula}}
		Let $\tau = (v^1,n_1;...;v^k,n_k)$ be a representation type such that $\sum_{i=1}^k n_iv^i = \tilde{v}$.
		The stratum of $\M_{\tilde{{\lambda}}}^0(\tilde{v})$ associated to $\tau$ has dimension
		\begin{equation}d(\tau) = 2\sum_{i=1}^k p(v^i).\end{equation}
	\end{theorem}
	
	Next, we follow \cite[Section 2.1.6]{bezrukavnikov2021etingof} to give a description of the structure of a quiver variety near a point. Similar results in the hyper-Kahler setting are first given in \cite[Section 6]{nakajima1994instantons}.
	
	Fix $\tau = (v^0,1;v^1,n_1;v^2,n_2;...;v^k,n_k)$ where $v^0_\infty = 1$, define a new quiver $\underline{Q}$, called the \textit{slice quiver}. It has $k$ vertices, in bijection with the $k$ dimension vectors $\{v^1,...,v^k\}$. The number of arrows from $i$ to $j$ is $-(v^i,v^j)$ if $i\neq j$, and we attach $1-\frac{1}{2}(v^i,v^i) = p(v^i)$ edge loops at each $i$. 
	
	Let the dimension vector $\underline{v} \in \ZZ_{\ge 0}^{\underline{Q}_0}$ be defined by $\underline{v}_i = n_i$. Let the framing $\underline{w} \in \ZZ_{\ge 0}^{\underline{Q}_0}$ be defined by $\underline{w}_i = -(v^0,v^i)$. 
	Let $H = \prod_{i=1}^k \GL(n_i) $, so that $\C^*\times H$ is identified with the subgroup of $G$ of automorphisms of a representation $V$ with type $\tau$. 
	Define the quiver varieties $\underline{\M}_0^{\theta}(\underline{v},\underline{w}) := \M_0^{\theta}(\underline{Q},\underline{v},\underline{w}),\underline{\M}_0^0(\underline{v},\underline{w}):= \M_0^0(\underline{Q},\underline{v},\underline{w})$. Here we abuse notation and write $\theta$ for the restriction of $\theta$ to $H$. 
	
	Write $\M_\lambda^0(v,w)^{\wedge_x} =\Spec (\C[\M_\lambda^0(v,w)]^{\wedge_x})$ where $\wedge_{x}$ in the right hand side denotes the completion at the maximal ideal corresponding to $x$. 
	For any vector space $R_0$, define $(\underline{\M}_0^0(\underline{v},\underline{w})\times R_0)^{\wedge_0}$ similarly. 
	
	We define  \begin{equation}\M_\lambda^\theta(v,w)^{\wedge_x} := \M_\lambda^0(v,w)^{\wedge_x} \times_{\M_\lambda^0(v,w)} \M_\lambda^\theta (v,w);
	\end{equation} 
	and define
	$\underline{\M}^{\theta}_0(\underline{v})^{\wedge_0}$ similarly. 
	\begin{proposition}[{\cite[Section 2.1.6]{bezrukavnikov2021etingof}}]{\label{structure of nbhd theorem}}
		There is an isomorphism $\M_\lambda^0(v,w)^{\wedge_x} \cong  (\underline{\M}_0^0(\underline{v},\underline{w})\times R_0)^{\wedge_0}$ as Poisson schemes, where $R_0$ is a symplectic vector space. 
		There is a commutative diagram of Poisson morphisms
		\begin{equation}
			\begin{tikzcd}
				\M_\lambda^\theta (v,w)^{\wedge_x}  \arrow[d,"\rho"] \arrow[r,"\sim"] 
				&
				(\underline{\M}^\theta_0(\underline{v},\underline{w})\times R_0)^{\wedge_0} \arrow[d,"\underline{\rho} \times \id"] \\
				\M_\lambda^0(v,w)^{\wedge_x}
				\arrow[r,"\sim"] & (\underline{\M}_0^0(\underline{v},\underline{w}) \times R_0)^{\wedge_0}
			\end{tikzcd}, 
		\end{equation}
		where $\underline\rho: \underline{\M}^\theta_0(\underline{v},\underline{w})\to \underline{\M}_0^0(\underline{v},\underline{w})$ is the natural projective map. 
	\end{proposition}
	
	In \Cref{Formal neighbourhood of minimal leaf}, we will obtain an analogue of this structure theorem, where we complete along affine open neighbourhoods of a symplectic leaf instead of completing at a point. 
	
	\subsection{More on quantizations}
	\label{more on quantizations subsection}
	In this subsection we recall more results about the quantizations $\A_\lambda, \A_\lambda^\theta$, and a quantum analogue of \Cref{structure of nbhd theorem}. 
	
	Recall for a filtered algebra $\A = \bigcup_{i\ge 0} \A_{\le i}$ the Rees algebra of $\A$ is the graded algebra $\A_\hbar = \bigoplus_{i} \A_{\le i} \hbar^i$. For example, we can define the algebra $\D(R)_\hbar$. Recall the quantum comoment map $\Phi: \g\to D(R)$ in \eqref{not symmetrized qcm}. The map $\Phi_\hbar: \g\to \A_\hbar, x\mapsto \hbar^2 x_R$ is also a quantum comoment map, in the sense that 
	\begin{equation}[\Phi_\hbar(x),d] = \hbar^2 x.d \end{equation}
	for any $d\in D(R)_\hbar$. In what follows, we will abuse notations and write $\Phi$ for $\Phi_\hbar$. 
	
	We need a \textit{symmetrized} quantum comoment map, defined as follows. Let 
	\begin{equation}\label{definition of varrho}
		\varrho(v)= -\frac{1}{2}\chi_{\bigwedge^{top} R}
	\end{equation}  
	where $\chi_{\bigwedge^{top} R}$ is the character of $\g$ given by the $\g$ action on the top exterior product of $R$.  Define the symmetrized quantum comoment map 
	\begin{equation}\label{symmetrized qcm}
		\Phi^{\sym} (x)= \Phi(x) - \hbar^2  \<\varrho(v),x\>. 
	\end{equation}
	There are two advantages of $\Phi^{\sym}$ over $\Phi$. First, it does not depend on the orientation of the arrows in the quiver. Moreover, it satisfies the following equation
	\begin{equation}\label{phisym is even}
		\sigma (\phi(x)) = \phi(x)
	\end{equation}
	where $\sigma$ is the unique antiautomorphism of $D(R)_\hbar$ satisfying that $\sigma|_R = \id_R, \sigma|_{R^*} = \id_{R^*}$, and $\sigma(\hbar) = -\hbar$. 
	
	Let us give another description of $\Phi^{\sym}$. Write the map in \eqref{not symmetrized qcm} as $\Phi_R$. The group $G$ also acts on $R^*$, so we similarly get a quantum comoment map $\Phi_{R^*}: \g \to D(R^*)$. We can identify $D(R^*)$ with $D(R)$ via an algebra isomorphism $F$, defined as follows. 
	If $\alpha \in R^* \subset D(R)$ is a linear function, the $F(\alpha) = - \d_{\alpha} \in Vect(R^*)$. If $\d_v \in D(R)$ is a vector corresponding to $v\in R$, then $F(\d_v) = v$ as a linear function on $R^*$. 
	It is then straightforward to check that 
	\begin{equation}
		\Phi^{\sym} = \frac{\Phi_R + F\circ \Phi_{R^*}}{2}.
	\end{equation} 
	
	Fix a representation type $\tau$ as in \eqref{general form of a representation type}, and recall the group $H$ in \eqref{definition of stabilizer H}. Let $\cartanh$ be the Lie algebra of $H$.
	Define the character $\underline{\varrho}(\underline{v})$ of $\cartanh$. Let $r_0$ be the natural restriction $(\g^*)^G\to (\cartanh^*)^H$. 
	Define $r: \C^{Q_0}\to \C^{\underline{Q}_0}$  by 
	\begin{equation}\label{definition of symmetrized restriction}
		r(\lambda) = r_0(\lambda -\varrho(v) ) +  \underline{\varrho}(\underline{v}).
	\end{equation} 
	
	The following is a quantum analogue of \Cref{structure of nbhd theorem}: 
	\begin{equation}
		(\A_\lambda)_\hbar^{\wedge_x} \cong (\underline{\A}_{r(\lambda)})^{\wedge_0}_\hbar \completeotimes_{\C[[\hbar]]} \WeylalgebraA_\hbar^{\wedge_0}.
	\end{equation}
	\subsection{Flower quiver varieties}\label{Flower quiver varieties}
	In this subsection we specialize previous results and constructions to the case of a flower quiver. 
	Let $Q$ be the quiver with a single vertex and $\ell \ge 2$ loops. We call $Q$ the flower quiver with $\ell$ loops.
	For example, the following picture is a flower quiver with $\ell = 7$; $w$ indicates the framing. 
	
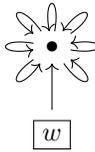
\[
	\begin{tikzcd}
		\bullet \arrow[out=-10,in=10,loop,swap]
		\arrow[out=35,in=55,loop,swap]
		\arrow[out=80,in=100,loop,swap]
		\arrow[out=125,in=145,loop,swap] 
		\arrow[out=170,in=190,loop,swap] 
		\arrow[out=215,in=235,loop,swap]
		\arrow[out=-55,in=-35,loop,swap]
		\\
		\boxed{w} \arrow[u]
	\end{tikzcd}
	\]
	
	Fix $n\ge 2, w\ge 1$. 
	The representation space is  
	\begin{equation}R =  (\sl_n)^{\oplus \ell} \oplus \Hom(\C^n ,\C^w).\end{equation}
	Its dual space is 
	\begin{equation}R^* = (\sl_n)^{\oplus \ell} \oplus \Hom(\C^w ,\C^n)\end{equation}
	where we identify $\sl_n$ with its dual via the trace pairing. 
	The group $G = \GL(n)$ acts naturally on $\C^n$ and induces a Hamiltonian action on $T^*R$. 
	If we write $x = (X_1,...,X_{\ell},Y_1,...,Y_{\ell},p,q) \in T^*R$, where $X_i\in \sl_n, Y_i \in (\sl_n)^* \cong \sl_n$ and $p: \C^w \leftrightarrows \C^n : q$, then the moment map is defined by
	\begin{equation}\mu(x) = \sum_{k=1}^{\ell} [X_k,Y_k] - pq.\end{equation}
	
	We define $\M_0^0(Q,n,w)$ and $\M^\theta_0(Q,n,w)$ as before. 
	When $\theta \neq 0$, $G$ acts on $\mu\inverse(0)^{\theta-ss}$ freely, so $\M^\theta_0$ is smooth. 
	
	As in \Cref{Framed vs nonframed quiver varieties}, we can view $\M^\theta_0$ (here $\theta$ can be 0) as a new quiver variety 
	$\M^\theta_0(Q^\infty,\tilde{v})$ without framing, where $Q^\infty$ is obtained from $Q$ by adding a vertex $\infty$ and $w$ edges from the original vertex to $\infty$. The new dimension vector $\tilde{v}$ is $1$ at $\infty$ and $n$ at the original vertex. If $\theta = k$, then the new stability condition, which we abuse notation and still denote by $\theta$, has component $-nk$ at the vertex $\infty$ and $k$ at the original vertex. 
	
	Let us denote the simple root associated to the original vertex by $\alpha$ and the root associated to $\infty$ by $\alpha_\infty$.
	We compute
	\begin{align}\label{some equation for pairing for quiver variety}
		\begin{split}
			& (\alpha,\alpha) = 2-2\ell, (\alpha,\alpha_\infty)  = - w; \\
			& p(n\alpha) = 1+n^2({\ell}-1) ;\\
			& p(\alpha_\infty+n\alpha) = n^2 ({\ell}-1) + nw .
		\end{split}
	\end{align}
	By the above computations and \Cref{CB simple dimension criteria,CB flatness of moment map}, we immediately get the following result. 
	\begin{proposition}\label{flatness of mu and existence of irrep with dimension nalpha}
		The following are true. 
		\begin{enumerate}
			\item The moment map $\mu$ is flat. 
			\item For any $0 < n'\le n $, there is an irreducible representation of $Q^\infty$ in $\mu\inverse(0)$ with dimension $n'\alpha$.
			\item For any $0 \le n'\le n $, there is an irreducible representation of $Q^\infty$ in $\mu\inverse(0)$ with dimension $\alpha_\infty + n'\alpha$. 
		\end{enumerate}
	\end{proposition}
	In particular, $\tau = (\alpha_\infty+n_0\alpha;n_1 \alpha,m_1;\cdots , n_k\alpha,m_k)$ is a valid representation type whenever $n_0+m_1n_1+ \cdots m_kn_k = n$. 
	For reasons that will be clear after \Cref{Relevant leaves}, we only care about leaves $\L_\tau$ associated to the representation type of the form
	$\tau = (\alpha_\infty+n_0\alpha,1;n_1\alpha,1;\cdots n_k\alpha,1)$. 
	By \Cref{CB strata dimension formula}, $\L_\tau$ has dimension 
	\begin{equation}\label{dimension for flower leaves}
		d(\tau) =2(p(\alpha_\infty+n_0\alpha) +  \sum p(n_i\alpha)) = 2\sum_{i=0}^k n_i^2({\ell}-1)+2n_0w+2k. 
	\end{equation}
	On the other hand, the open leaf corresponds to $\tau_O = (\alpha_\infty + n\alpha,1)$, and we see that $\dim \M_0^0 = 2(n^2 ({\ell}-1) + nw ).$
	
	\section{Relevant leaves}\label{Relevant leaves}
	Let $Q$ be the flower quiver with $\ell \ge 2$ loops. Let $n$ be a dimension vector. 
	Consider the leaf $\L_\tau$ corresponding to a representation type $\tau$
	\begin{equation}\tau = (\alpha_\infty+n_0\alpha,1; n_1\alpha,m_1; n_2\alpha,m_2;\cdots n_k\alpha,m_k).\end{equation}
	By construction, the quiver associated to the slice to $\L_\tau$ has $k+1$ vertices. The vertex corresponding to $n_i\alpha_i,m_i$ carries $n_i^2({\ell}-1)+1$ edge loops.The dimension vector for defining the slice quiver variety $\underline{\M}_0^0$ has component $m_i$ at the vertex corresponding to $n_i\alpha$. 
	
	Suppose $\lambda \in \C^{Q_0}$ and $\A_\lambda$ is as defined in \eqref{definition of quantum algebra Alambda}. 
	We only care about symplectic leaves whose closure is the support of some HC $\A _\lambda$-bimodule. The following corollary follows from (2) of \Cref{restriction functor and support property}. 
	\begin{corollary}
		If the closure of the symplectic leaf $\L_\tau$ is the support of a HC $\A_\lambda$-bimodule, then the slice algebra $\underline{\A} _{r(\lambda),\tau}$ has a nonzero finite dimensional HC bimodule. 
	\end{corollary}
	
	Such representation types are characterized by the following proposition. 
	\begin{proposition}{\label{only single loops carries finite dim modules}}
		If the algebra $\underline{\A} _{\lambda,\tau}$ has a nonzero finite dimensional HC bimodule, then all $m_i = 1$. 
	\end{proposition} 
	We remark that the condition all $m_i = 1$ is not sufficient for $\underline{\A} _{\lambda,\tau}$ to have a nonzero finite dimensional module. The existence of such modules also depends on $\lambda$. 
	
	Let $\M_0^0(Q,v,w)$ be the affine quiver variety associated to some quiver $Q$, with dimension vector $v$, and framing $w$. Let $\theta>0$ be generic and consider $\rho:
	\M^\theta_0(Q,v,w)\to \M^0_0(Q,v,w)$. \Cref{only single loops carries finite dim modules} follows from the following geometric result. 
	\begin{theorem}\label{only single edge loop has lagrangian component}
		If the fiber $\rho\inverse(0)$ has dimension $\frac{1}{2}\dim \M^\theta_0(Q,v,w)$, then for each vertex $i\in Q_0$ such that $v_i\ge 2$, it carries at most one edge loop. 
	\end{theorem}
	\begin{remarklabeled}
		By \cite[Proposition 1.2]{kaledin2000symplectic}, the map $\rho$ is semismall. Therefore, the maximal possible dimension of $\rho\inverse(0)$ is $\frac{1}{2}\dim \M^\theta_0(Q,v,w)$. 
	\end{remarklabeled}
	\begin{proof}[Proof of \Cref{only single loops carries finite dim modules}]
		Assume \Cref{only single edge loop has lagrangian component}. Suppose $\B$ is a finite dimensional HC $\underline{\A} _{\lambda,\tau}$-bimodule. 
		Then, since $0$ is the only $\C^*$-fixed point of $\M_0^0$, $\gr \B$ is supported at $0\in \underline{\M}_0^0(\underline{v},\underline{w})$. 
		For a generic stability condition $\theta$, the localization $\underline{\A}^\theta_{\lambda,\tau} \otimes_{\underline{\A}_{\lambda,\tau}} \B$ is a coherent sheaf of $\underline{\A}^\theta_{\lambda,\tau}$-modules which is supported on $\rho\inverse(0)$. 
		By Gabber's involutivity theorem, its support is coisotropic. By \Cref{only single edge loop has lagrangian component}, this is possible only if all $m_i = 1$, and \Cref{only single loops carries finite dim modules} follows. 
	\end{proof}
	
	To prove \Cref{only single edge loop has lagrangian component}, we need the following constructions from \cite{bozec2016Quiverswithloops}. 
	Let $(x,y,p,q)\in T^*R$ where $(x,p)\in R$ and $(y,q) \in R^*$. 
	Let the torus $\C^*$ act on $T^*R$ by $t.(x,y,p,q) = (tx,y,p,tq)$. It induces a $\C^*$-action on ${\M}_0^\theta(v,w) $. We denote the subset of $\C^*$-fixed points in ${\M}_0^\theta(v,w) $ by $F$. 
	
	By \cite[Section 5]{nakajima1994instantons}, $F\subset\M_0^\theta(v,w)$ is a nonsingular subvariety. Let $F = \coprod_{s\in S} F_s$
	denote its decomposition into connected components. 
	
	\begin{definition}[\cite{bozec2016Quiverswithloops}]
		Define the following subset of ${\M}_0^\theta (v,w)$: 
		\begin{equation}
			\mathfrak{L}(v,w) = \{x\in {\M}_0^\theta(v,w) | \lim_{t\to \infty} t.x \text{ exists}\}.
		\end{equation}
	\end{definition}
	The limit $\lim_{t\to \infty} t.x$ for any $x\in \mathfrak{L}(v,w)$ lies in $F$. 
	Let $\mathfrak{L}(v,w)_s = \{x\in \mathfrak{L}(v,w) |  \lim_{t\to \infty} t.x \in F_s\}$. The connected components of $\mathfrak{L}(v,w)$ are exactly $\mathfrak{L}(v,w)_s, s\in S$. When $v,w$ is clear from the context, we will write $\mathfrak{L}_s$ for $\mathfrak{L}(v,w)_s$. 
	\begin{lemma}[{\cite[Proposition 2.2]{bozec2016Quiverswithloops}}]\label{Bozec limit set is Lagrangian}
		The subvariety $\mathfrak{L}(v,w) \subset \M_0^\theta(v,w)$ is Lagrangian. \qed
	\end{lemma}
	\begin{lemma}\label{rhoinverse 0 lies in limit set}
		We have $\rho\inverse(0)\subset \mathfrak{L}(v,w)$.
	\end{lemma}
	\begin{proof}
		Since $\rho^{-1}(0)$ is projective, the image of any morphism $\C^* \to \rho\inverse(0) $ has a limit at $\infty$. 
	\end{proof}

	We now consider the existence of Lagrangian components of $\rho\inverse(0)$. Let $x \in \mathfrak{L}_s$ be a $\C^*$-fixed point. Suppose the vertex $i$ has an edge loop $a$. Let $x_a\in \sl(v_i)$ be the component of $x$ corresponding to $a$. 
	\begin{lemma}\label{lemma in proving no lagraigian component}
		If $\rho\inverse(0)$ has a Lagrangian component, then $\rk x_a = v_i-1$. 
	\end{lemma}
	
	\begin{proof}
		By \Cref{Bozec limit set is Lagrangian} and \Cref{rhoinverse 0 lies in limit set}, a Lagrangian component must be some $\mathfrak{L}_s$. In particular, $\rho(F_s) = 0$. 
		Suppose $\rk x_a <v_i-1$. 
		We will show that $\rho(F_s)$ contains a nonzero element, which leads to a contradiction.
		
		Let $y_a\in\End(V_i)$ be associated to the opposite loop of $a$. Let $r = (\cdots ,x_a,y_a,\cdots) $ be a representative of $x$ in $\mu\inverse(0)^{\theta-ss}$. 
		Since $x$ is $\C^*$-fixed, we have 
		\begin{equation}[\cdots, x_a,y_a,\cdots] = [\cdots, tx_a,y_a,\cdots]\end{equation}
		in $\M_0^\theta(v,w)$ for all $t\in \C^*$. Here we use the notation in \eqref{representative of an element of quiver variety in T*R}. 
		
		Since the action of $\GL(V_i)$ on $\mu\inverse(0)^{\theta - ss}$ is free, $\Stab_{\C^* \times \GL(V_i)} x \cong \C^*$. It follows that there is a one parameter subgroup $g(t) \in \GL(V_i)$,  such that 
		\begin{equation}\label{definition of C* invariant element}
			g(t) x_a g(t)\inverse = tx_a, g(t) y_a g(t) \inverse = y_a.
		\end{equation}
		We may fix a basis of $V_i$ in which $g(t) = \diag(t^{n_1},\cdots,t^{n_{v_i}})$, and $n_1 \ge n_2\ge\cdots n_{v_i}$. 
		It follows from \eqref{definition of C* invariant element} that $(x_a)_{ij} \neq 0$ only if $n_i-n_j = 1$. In particular, $x_a$ is strictly block upper triangular with 0 diagonal blocks. More precisely, suppose 
		\begin{equation}\label{block description of g(t)}
			\begin{split}
				&n_1=n_2=\cdots=n_{p_1} \\
				>& n_{p_1+1}=\cdots =n_{p_1+p_2} \\ 
				>& \cdots\\
				>& n_{p_1+ \cdots p_{k-1}+1} = \cdots =n_{p_1+ \cdots p_k} 
			\end{split}
		\end{equation}
		where $\sum p_j = v_i$. Then $x_a$ has $k-1$ blocks $A_{12}, A_{23}, \cdots, A_{k-1,k}$, where $A_{j,j+1}$ has size $p_j \times p_{j+1}$, and $A_{j,j+1} \neq 0$ only if $p_{j+1} = p_j +1$.

		We will construct an element $ y' \in \sl(V_i)$ satisfying the following properties:
		\begin{enumerate}[label=(\alph*)]
			\item $g(t)  y' g(t)\inverse = y'$; 
			\item $[x_a,y'] = 0$;
			\item $y'$ is not nilpotent.
		\end{enumerate}
		Suppose we have constructed such a $y'$. Then for any $z\in\C$, the element 
		\begin{equation}
			r_z : = [\cdots, x_a,y_a+zy',\cdots]
		\end{equation}
		still lies in $F_s$. By (c), $y_a+zy'$ is not nilpotent for some $z$, so $\rho(r_z) \neq 0$, contradiction. 
		
		We proceed to construct $y'$, in 2 steps. 
		\begin{enumerate}[label=\textit{Step} \arabic*.]
			\item Suppose $x_a$ has Jordan canonical form. By the assumption $\rk x_a < v_i-1$, there exists $1\le j \le v_i-1$ such that the $(j,j+1)$-entry of $x_a$ is 0. 	
			Let 
			\begin{equation}
				y' = \diag \left(1,1,...,1,-\frac{j}{v_i-j},...,-\frac{j}{v_i-j}\right) \in \sl(v_i),\end{equation}
			where the first $j$ diagonal entries are $1$. 
			Then $[x_a,y']=0$, $g(t) y' g(t)\inverse = y'$, and $y'$ is not nilpotent, as desired. 
			
			\item We construct an element $P\in \GL_{v_i}$ such that $P\inverse x_a P = J$ is in Jordan canonical form, and for any diagonal matrix $D\in \GL_{v_i}$, $P\inverse D P$ is in block diagonal form, with block size $p_1,p_2,\cdots, p_k$. 
			
			Suppose $P\inverse x_a P = J$ and  the columns of $P$ are $u_1,\cdots,u_{v_i}$. Then $x_a u_j = u_{j-1}$ or $=0$. Therefore, thanks to the block upper triangular form of $x_a$, we may assume each $u_j$ is nonzero in only one block, i.e. 
			\begin{equation}
				u_j = 
				\begin{pmatrix}
					0\in \C^{p_1} \\
					\vdots \\
					0\in \C^{p_{m-1}} \\ 
					\neq 0 \in \C^{p_m} \\
					0\in \C^{p_{m+1}} \\
					\vdots\\
					0\in \C^{p_k} 
				\end{pmatrix}
			\end{equation}
			Since each column of $P$ has the above form, we can find a permutation matrix $\sigma$ such that $P\sigma$ is block diagonal with block size $p_1,\cdots,p_k$. 
			
			Now apply Step 1 to $J$; we get a diagonal, non nilpotent $y'$ such that 
			\begin{equation}
				[y',J] = 0.
			\end{equation}
			Let $y'' = \sigma\inverse y' \sigma$, which is still diagonal, traceless and not nilpotent. Then, $P\sigma y'' (P\sigma)\inverse$ is block diagonal, and 
			\begin{equation}
				[P\sigma y'' (P\sigma)\inverse, x_a] = [y',P\inverse x_a P] = [y',J] = 0.
			\end{equation}
			The element $P\sigma y'' (P\sigma)\inverse$ has desired properties. This completes the proof.
		\end{enumerate}		
	\end{proof}
	
	We can now prove \Cref{only single edge loop has lagrangian component}. 
	\begin{proof}[Proof of \Cref{only single edge loop has lagrangian component}]
		Suppose some vertex $i$ carries edge loops $a$ and $b$, and $v_i\ge 2$. Suppose $r = [\cdots,x_a,y_a,x_b,y_b,\cdots] \in F_s$ where $x_a,x_b\in \sl(V_i)$ are associated to two different edge loops at $i$, and $y_a,y_b$ are associated to the opposite loops respectively. 
		By \Cref{lemma in proving no lagraigian component}, we may assume $x_a$ is a single nilpotent Jordan block, and the one parameter subgroup $g(t)$ takes the form 
		\begin{equation}
			g(t) = \diag(t^{v_i}, t^{v_i-1}, \cdots, t ).
		\end{equation}
		Since $g(t) x_b g(t)\inverse = tx_b$, we see $(x_b)_{ij} = 0$ unless $j-i=1$. 
		
		Let $y_a' = \diag(a_1,\cdots , a_{v_i}), y_b' = \diag(b_1,...,b_{v_i})$, so that $g(t)  y_\bullet' g(t)\inverse =  y_\bullet'$, for $\bullet = a,b$. Consider the linear system of equations
		\begin{equation}\label{system of diagonal ya yb for deformation}
			\left\{ 
			\begin{array}{cc} 
				&[x_a, y_a'] +[x_{b},y_b'] = 0  \\
				&\tr y_a' = \tr y_b' = 0
			\end{array}
			\right.
			.
		\end{equation}
		More precisely, the first equation is equivalent to 
		\begin{equation}
			a_{j+1}-a_j + (x_b)_{j,j+1}(b_{j+1}-b_j)=0, 1\le j \le v_i-1.
		\end{equation}
		Since there are more variables than equations, the system \eqref{system of diagonal ya yb for deformation} always has a nonzero solution. Fix such a pair $y_a',y_b'$, so that at least one of $y_a',y_b'$ is not nilpotent. 
		Consider the family $r_z$ obtained from $r$ by replacing $y_a$ with $y_a + z y_a'$ and $y_b$ with $y_b + z y_b'$. 
		By construction, we get a family of elements in $F_s$. But $\rho ( r_z )\neq 0$ for some $z\neq 0$, so $F_s \not\subset \rho\inverse(0)$. In particular, $\rho\inverse(0)$ cannot contain a Lagrangian component $F_s$. 
	\end{proof}
	
	\begin{definition}
		A leaf 
		$\L_\tau$ is called relevant if it corresponds to representation type of the form
		\begin{equation}\label{general form of a relevant rep type}
			\tau = (\alpha_\infty+n_0\alpha,1; n_1\alpha,1; n_2\alpha,1;\cdots n_k\alpha,1).
		\end{equation}
	\end{definition}
	The support of a Harish-Chandra bimodule must be the union of the closures of some relevant leaves. 
	In the remaining of this paper, we only care about relevant leaves.
	
	\begin{proposition}
		The dimension of the leaf corresponding to the representation type 
		$\tau = (\alpha_\infty+ n_0\alpha,1; n_1\alpha,1; n_2\alpha,1;\cdots n_k\alpha,1)$ is given by 
		\begin{equation}\label{dimension of relevant leaf}
			d_\tau = 2k + 2n_0 w + (2\ell -2) \sum_{i=0}^k n_i^2.
		\end{equation}
	\end{proposition}
	\begin{proof}
		This follows from \Cref{CB strata dimension formula} and \eqref{some equation for pairing for quiver variety}. 
	\end{proof}
	To finish this section, we describe the boundary of a relevant leaf. Let $\tau$ be as in \eqref{general form of a relevant rep type} and $\L_\tau$ be a relevant leaf; let
	$\tau' = (\alpha_\infty+ r_0\alpha,1; r_1\alpha,m_1; r_2\alpha,m_2;\cdots r_{k'}\alpha,m_{k'})$ be an arbitrary representation type. Define 
	\begin{equation}
		\tau_{rel}' = (\alpha_\infty+ r_0\alpha,1; 
		\underbrace{ r_1\alpha,1;\cdots,r_1\alpha,1}_{m_1\text{ copies} };\cdots
		\underbrace{ r_2\alpha,1;\cdots,r_2\alpha,1}_{m_2\text{ copies} };
		\underbrace{r_{k'}\alpha,1;\cdots,r_{k'}\alpha,1}_{m_{k'}\text{ copies} }).
	\end{equation}
	It is clear that $\tau'_{rel}$ defines a relevant leaf whose closure contains $\L_{\tau'}$. 
	\begin{proposition}\label{proposition leaf closure} 
		Let $\tau =(\alpha_\infty+ n_0\alpha,1; n_1\alpha,1; n_2\alpha,1;\cdots n_k\alpha,1)$.	
		The following are true. 
		\begin{enumerate}
			\item Let $\tau' = (\alpha_\infty+ r_0\alpha,1; r_1\alpha,m_1; r_2\alpha,m_2;\cdots r_{k'}\alpha,m_{k'})$. Then $\L_{\tau'} \subset \d {\L_\tau}$ if and only if $\L_{\tau'_{rel}} \subset \d {\L_\tau}$.
			\item Now let $\tau' = (\alpha_\infty+ r_0\alpha,1; r_1\alpha,1; r_2\alpha,1;\cdots r_{k'}\alpha,1)$. 
			Then $\L_{\tau'}\subset \d {\L_\tau}$ if and only if the following condition holds: upon reordering $r_1,\cdots, r_{k'}$, there are integers $0\le a_0<a_1 \cdots <a_k=k'$, such that 
			\begin{equation}\label{finer partition means boundary of leaf}
			\begin{split}
					r_0+r_1\cdots + r_{a_0} &= n_0 ;\\
					r_{a_1+1}+\cdots + r_{a_2} &= n_1;\\
					&\vdots\\
					r_{a_{k-1}+1}+\cdots + r_{a_k} &= n_k.
			\end{split}
\end{equation}
			and $\tau\neq \tau'$. 
		\end{enumerate}
	\end{proposition}
	\begin{proof}
		Let us prove (1). 
		The ``if" part is clear. Suppose $\L_{\tau'} \subset \d {\L_\tau}$. Let $H \subset G$ (resp. $H'$, $H'_{rel}$) be the stabilizer of a representative  $x \in T^*R$ (resp. $x',x'_{rel}$) of an element of $\L_\tau$ (resp $\L_{\tau'}, \L_{\tau'_{rel}}$).  
		By \cite[Proposition 3.15]{bellamy2021symplectic}, $\L_{\tau'} \subset \d {\L_\tau}$ if and only if $H$ is conjugate to a proper subgroup of $H'$. Note that $H \cong (\C^*)^k$ is a torus, $H' = \prod\GL(m_i)$, and the maximal torus of $H'$ is precisely $H'_{rel}$. It follows that any proper embedding $H \hookrightarrow H'$ factors through $H'_{rel}$ up to conjugation. Therefore, $\L_{\tau'_{rel}}\subset \d\L_\tau$. 
		
		Let us prove (2). Note that the groups $H,H'$ takes the following block diagonal forms in $G=\GL(n)$: 
		\begin{equation}
			H = \left(
			\begin{array}{cccc}
				\id_{n_0} & & & \\
				& \C^* \otimes \id_{n_1} & & \\
				& & \ddots & \\
				& & & \C^*\otimes\id_{n_k} 
			\end{array}
			\right), H'= \left(
			\begin{array}{cccc}
				\id_{r_0} & & & \\
				& \C^* \otimes \id_{r_1} & & \\
				& & \ddots & \\
				& & & \C^*\otimes\id_{r_{k'}} 
			\end{array}
			\right).
		\end{equation}
		Here, $\id_n$ denotes the identity matrix of size $n$. The group $H$ is a subgroup of $H'$ upon conjugation by $G$ if and only if \eqref{finer partition means boundary of leaf} holds. 
	\end{proof}
	
	The minimal relevant leaf $\L$ corresponding to $ (\alpha_\infty,1; \alpha,1; \alpha,1;\cdots \alpha,1)$ has dimension $2n \ell$. 
	By \Cref{proposition leaf closure}, any leaf in $\d\L$ corresponds to representation type of the form 
	$\tau = (\alpha_\infty,1; \alpha,m_1;\cdots; \alpha,m_k)$, which has dimension $2\ell k$. 
	The following corollary is then clear. 
	
	\begin{corollary}\label{boundary of minimal leaf has codim 4}
		For the minimal leaf $\L$, $\codim_{\overline{\L}}\d\L \ge 4$. 
	\end{corollary}
	\begin{remarklabeled}\label{remark cod 4 condition can fail}
		\Cref{boundary of minimal leaf has codim 4} may fail for non minimal relevant leave. For example, consider $\tau'= (\alpha_\infty,1; \alpha,1; \cdots; \alpha,1; 2\alpha,1)$. The leaf $\L_{\tau '}$ is relevant, and by \Cref{proposition leaf closure}, the minimal leaf $\L$ lies in the boundary of $\L_{\tau '}$. By \eqref{dimension of relevant leaf}, $\dim \L_{\tau'} = 2n \ell + 4 \ell-6$. When $\ell = 2$, we see $\codim_{\overline{{\L_{\tau'}}}}\d{\L_{\tau'}} = 2$. 
	\end{remarklabeled}
	
	\section{Slice to the minimal relevant leaf}\label{section slice to minimal leaf}
	In this bookkeeping section we describe the slice quiver, the slice algebra and its quantizations associated to the minimal relevant leaf. We also define certain group actions that will be useful later. 
	
	\subsection{The slice quiver}\label{subsection the slice quiver}
	
	Recall the minimal relevant leaf corresponds to the representation type $(\alpha_\infty,1;\alpha,1; \cdots;\alpha,1)$. 
	By the construction in \Cref{subsection Symplectic leaves and local structures}, the associated slice quiver $\sliceQ$ has $n$ vertices. Let us label them by $\alpha_1,\ldots,\alpha_n$. For each pair $1 \le i\neq j\le n$, there are $2\ell -2$ arrows between the vertices $\alpha_i,\alpha_j$; there are $\ell$ edge loops at $\alpha_i$ for each $1\le i \le n$. 
	
	Let us orient the arrows in the most symmetric way, i.e. for each pair $1\le i \neq j \le n$, there are ${\ell}-1$ arrows from $\alpha_i$ to $\alpha_j$ and ${\ell}-1$ arrows from $\alpha_j $ to $ \alpha_i$. 
	
	The dimension vector is $\underline{v} = (1,1,\cdots,1)$ and the framing is $\underline{w} = (w,w,\cdots,w)$. 
	Let $H =  \GL(\underline{v}) = (\C^*)^n$. 
	
	Note that, since we have assigned traceless endomorphisms to the edge loops in the original flower quiver, there are some conditions on the endomorphisms we assign to the edge loops in $\sliceQ$.
	However, since all components of the dimension vectors are 1, the edge loops will contribute to a direct product component $\C^{2n}$ in the slice quiver variety. Therefore, we may remove the edge loops from $\sliceQ$. 
	The representation space $T^*\sliceR$ has dimension $(n^2-n)(2\ell -2) + 2nw$; this dimension is precisely the number of arrows in the extended double quiver of $\sliceQ^\infty$, excluding the edge loops. 
	
	For each pair $1\le i \neq j\le n$, let $v_{ij}^k\in\C, 1\le k \le {\ell}-1 $ be the homomorphisms corresponding to the arrows from $i$ to $j$.  Also, for each $1\le i\le n$ and  $1\le k \le w$, let $p_i^k\in\C$ represent the framing. 
	
	\subsection{The slice algebra}\label{subsection the slice algebra}
	Set $\underline{A}:= \C[\underline{\M}_0^0(\underline{v},\underline{w})]$; $\underline{\M}_0^0(\underline{v},\underline{w})$ is defined in \Cref{subsection Symplectic leaves and local structures}. 
	Let us first label the linear functions on $T^*\sliceR$. For each pair $1\le i \neq j\le n$ and each $1\le k \le {\ell}-1$, let $x_{ij}^k$ be the linear function dual to $v_{ij}^k$. For each pair $1\le i \le n$ and each $1\le k \le w$, let $(p_i^k)^*$ be the corresponding linear function dual to $p_i^k$. 
	Let $y_{ij}^k$ and $(q_i^k)^*$ correspond to the opposite arrows of $x_{ij}^k$ and $(p_i^k)^*$ in the double quiver, respectively. 
	
	If we label an element of $H$ by $h = (h_1,\cdots, h_n)$ where $h_i\in \C^*$ corresponds to the vertex $\alpha_i$, then it is easy to see 
	\begin{equation}\label{H action on algebra on TstarsliceR}
		\begin{split}
			h.x_{ij}^k = h_i\inverse h_j x_{ij}^k, h.y_{ij}^k = h_i h_j\inverse y_{ij}^k; \\
			h.(p_i^k)^* = h_i (p_i^k)^*, h.(q_i^k)^* = h_i\inverse (q_i^k)^*.
		\end{split}
	\end{equation} 
	By definition, the slice algebra $\underline{A}$ is the Hamiltonian reduction 
	$\C[\mu_H\inverse(0)]^H$, where $\mu_H$ is the moment map $T^*\sliceR \to \cartanh^*$. 
	\subsection{Quantization parameters}\label{subsection morphism of Quantization parameters}
	In this subsection we compute the morphism $r$ relating quantization parameters defined in \eqref{definition of symmetrized restriction}. First, we compute the characters $\varrho$ and $\underline{\varrho}$ \eqref{definition of varrho}. 
	
	The group $G = \GL(n)$ acts on the endomorphism spaces associated to each edge loop of the flower quiver by conjugation, and acts on the framing space by left multiplication. Therefore, only the framing part contributes to $\varrho$ and 
	\begin{equation}
		\varrho = -\frac{1}{2}\chi_{\bigwedge^{top} R} = -\frac{w}{2}\tr. 
	\end{equation}
	
	Note that we have chosen the orientation of $\sliceQ$ such that for each pair $1\le i\neq j\le n$, the dimensions of $H$-weight spaces in $\sliceR$ with weights $h_ih_j\inverse$ and $h_jh_i\inverse$ are equal. 
	Therefore, the character $\underline{\varrho}$ of $\cartanh$ comes only from the framing part. For $\eta =( \eta_1,\eta_2,\cdots,\eta_n) \in \cartanh$, 
	\begin{equation}
		\underline{\varrho}(\eta) = -\frac{w}{2} (\eta_1 + \eta_2+\cdots +\eta_n). 
	\end{equation}
	In our situation, the definition \eqref{definition of symmetrized restriction} reads
	\begin{equation}
		r(\lambda) =  (\lambda,\cdots,\lambda). 
	\end{equation}
	From now on, we will write for abbreviation $\sliceA_{\lambda} := \sliceA_{r(\lambda)}$ and $\sliceA_{\lambda,\hbar} := \sliceA_{r(\lambda),\hbar}$. 
	\subsection{$S_n$- and $\C^*$-actions}
	We equip the algebras $\underline{A}$ and $\sliceA_\lambda$ with a $S_n$-action and two $\C^*$-actions.
	
	Note that $S_n$ acts on $T^*\sliceR$ by permuting the $n$ vertices of the slice quiver, and the subset $\mu_\hbar\inverse(0)$ is $S_n$-stable. This induces an $S_n$-action on $\sliceA$. 
	Moreover, since $r(\lambda) = (\lambda,\cdots,\lambda)$ is $S_n$-invariant, $S_n$ also acts on the algebras $\sliceA_\lambda$ and $\sliceA_{\lambda,\hbar}$. 
	
	The first $\C^*$-action is defined as follows. 
	\begin{definition}\label{definition of model torus action}
		Define the following $\C^*$-action on $\C[T^*\sliceR]$, by specifying the degrees of linear functions: 
		\begin{equation}
			\begin{split}
				&\deg x_{ij}^1 = 2, \deg y_{ij}^1 = 0\\
				&\deg x_{ij}^k =  \deg y_{ij}^k = 1, \text{ for $1\le i\neq j\le n, 2 \le k\le {\ell}-1$}\\
				&\deg p_i^k = \deg q_i^k = 1, \text{ for $1\le i\le n, 1<k\le w$}. 
			\end{split}
		\end{equation}

		We call this action the \textit{model} $\C^*$-action. We will show later in \Cref{subsection model varieties} that this action is constructed so that the morphism from the ``model variety" to $T^*R$ is $\C^*$-equivariant. 
	\end{definition}
	The second $\C^*$-action is more natural. 
	\begin{definition}\label{definition of conical torus action}
		Define the following $\C^*$-action on $T^*\sliceR$, by specifying the degrees of linear functions: 
		\begin{equation}
			\begin{split}
				&\deg x_{ij}^k =  \deg y_{ij}^k = 1, \text{ for $1\le i\neq j\le n, 1\le k\le {\ell}-1$}\\
				&\deg p_i^k = \deg q_i^k = 1, \text{ for $1\le i\le n, 1<k\le w$}. 
			\end{split}
		\end{equation}
		
		The degree 0 part of $\C[T^*\sliceR]$ is $\C$, and the maximal ideal of $0 \in T^*\sliceR$ is positively graded. 
		We call this action the \textit{conical} $\C^*$-action. 
	\end{definition}
	Both the model $\C^*$-action and the conical $\C^*$-action commute with the $H$-action and $S_n$-action on $T^*\sliceR$. Moreover, the symplectic form $\omega$ on $T^*R$ has degree 2 with respect to both $\C^*$-actions. 
	
	Let us describe the difference of the two $\C^*$-actions. 
	Define
	\begin{equation}
		\alpha = \sum_{i\neq j} x_{ij}^1 y_{ij}^1.
	\end{equation}
	This is an $H$-invariant element of $\C[T^*\sliceR]$ and hence descends to an element of $A$. It has degree 2 with respect to both $\C^*$-actions. 
	Similarly, define 
	\begin{equation}\label{definition of alpha relating two torus actions}
		\alpha_\hbar = \sum_{i\neq j} \hbar^2 x_{ij}^1 y_{ij}^1
	\end{equation}
	which is an $H$-invariant element of the homogenized Weyl algebra $\WeylalgebraA_{T^*\sliceR,\hbar}$ and hence descends to an element of $\sliceA_{\lambda,\hbar}$. The following proposition is straightforward to check. 
	\begin{proposition}\label{relation between model and conical Cstar action}
		Let $\eu_m$ be the Euler derivation on $\sliceA$ or $\sliceA_{\lambda,\hbar}$ induced by the model $\C^*$-action, and $\eu_c$ the Euler derivation on $\underline{A}$ or $\sliceA_{\lambda,\hbar}$ induced by the model $\C^*$-action. Then 
		\begin{enumerate}
			\item $\eu_m - \{\alpha,-\} = \eu_c $ on $\underline{A}$.
			\item $\displaystyle \eu_m - \frac{1}{\hbar^2}\ad \alpha_\hbar  = \eu_c $ on $\sliceA_{\lambda,\hbar}$.
		\end{enumerate}
	\end{proposition}
	\section{Formal neighbourhood of the minimal leaf}\label{Formal neighbourhood of minimal leaf}
	In this section we describe the formal neighbourhood of the minimal relevant symplectic leaf $\L$ of $\M$. Our goal is to obtain a decomposition similar to \eqref{formal nbhd of point in introduction}. 
	We will first work in $T^*R$, and then do Hamiltonian reduction. 
	
	\subsection{Affine opens}\label{subsection good affine open cover}
	In this subsection, we describe an affine open covering of ${\L}$. We start with the preimage of $\L$ in $T^*R$.
	
	Up to conjugation by some element of $G$, any element of $T^*R$ that lies in the preimage of $\L$ has the form 
	\begin{equation}(X_1,...,Y_{\ell},p,q),
	\end{equation}
	where all $X_i,Y_i$ are traceless diagonal matrices and $p,q=0$; moreover, for each $1\le j \neq k \le n$, the $2\ell $-tuples consisting of the $j$-th and $k$-th diagonal entry, i.e. 
	\begin{equation}\label{jth diagonal entry of s in L}
		((X_1)_{j,j},(X_2)_{j,j}, \cdots, (Y_{\ell})_{j,j})
	\end{equation}
	and 
	\begin{equation}\label{kth diagonal entry of s in L}
		((X_1)_{k,k},(X_2)_{k,k}, \cdots, (Y_{\ell})_{k,k})
	\end{equation}
	are distinct elements of $\C^{2\ell }$. 
	
	Let $\tilde{\L}$ denote the space of all diagonal $X_i,Y_i$ satisfying the above conditions. Then, $\tilde{\L}$ is a non-affine open subvariety of the vector space $\C^{2\ell (n-1)}$. There is a natural locally closed embedding $\tilde{\L} \hookrightarrow T^*R, (X_1,\cdots Y_{\ell}) \mapsto (X_1,\cdots Y_{\ell}, p=0, q=0) $. Moreover, $\tilde{\L}$ is the universal cover of $\L$ with Galois group $S_n$. 
	
	Consider a nonzero vector $K = (a_1,\cdots,a_\ell,b_1,\cdots,b_\ell)\in \C^{2\ell}$. For any $s = (X_1,\cdots,Y_\ell) \in \tilde{\L}$, set 
	\begin{equation}
		K\bullet s = \sum_{i=1}^\ell a_i X_i + \sum_{i=1}^\ell b_i Y_i
	\end{equation}
	which is a traceless diagonal matrix. 
	Define 
	\begin{equation}\label{definition of tildeLK}
		\tilde{\L}_K = \left\{ s\in \tilde{\L} | K\bullet s  \text{ has pairwise distinct diagonal entries }  \right\}
	\end{equation}
	For example, when $K = (1,0,\cdots,0)$, $\tilde{\L}_K$ consists of $s$ whose $X_1$ component has pairwise distinct diagonal entries. 
	
	If $K'=tK$ for some $t\in\C^*$, then $\tilde{\L}_K = \tilde{\L}_{K'}$. Therefore, we may assume $K  = [a_1,\cdots, b_\ell] \in \PP^{2\ell -1}$. 
	The following proposition summarizes the properties of the open subsets $\tilde{\L}_K$. 
	\begin{lemma}\label{lemma properties of good affine open cover}
		The following are true. 
		\begin{enumerate}
			\item Each $\tilde{\L}_K$ is a $\C^*$-stable affine open subset of $\tilde{\L}$, and $\codim_{\tilde{\L}} \tilde{\L}_K = 1$. 
			\item $\displaystyle \tilde{\L} = \bigcup_{K\in\PP^{2\ell -1}} \tilde{\L}_K$.
			\item All $\tilde{\L}_K$ are isomorphic as symplectic affine varieties.
			\item For $K_1\neq K_2 \in \PP^{2\ell -1}$, 
			$\codim_{\tilde{\L} } \tilde{\L} \setminus (\tilde{\L}_{K_1}\cup \tilde{\L}_{K_2}) = 2$.
			\item For pairwise distinct $K_1,K_2,K_3\in \PP^{2\ell -1}$, $\codim_{\tilde{\L}_{K_1}} \tilde{\L}_{K_1} \setminus (\tilde{\L}_{K_2}\cup \tilde{\L}_{K_3}) = 2$.
		\end{enumerate}
		
	\end{lemma}
	\begin{proof}
		Note that the complement to $\tilde{\L}_K$ in $\tilde{\L}$ is the union of $\binom{n}{2}$ hyperplanes, $\bigcup_{i\neq j} H_{K,i,j}$, where
		\begin{equation}
			H_{K,i,j} = \{s\in \tilde{\L}| (K\bullet s)_{ii} - (K\bullet s )_{jj} = 0\} .
		\end{equation} Since each $H_{K,i,j}$ is $\C^*$ stable, (1) follows. 
		
		(2) is a consequence of the fact that \eqref{jth diagonal entry of s in L} and \eqref{kth diagonal entry of s in L} are pairwise distinct for each pair $j\neq k$. 
		
		Let us prove (3). Note that the group $\Sp(2\ell )$ acts naturally on the symplectic vector space $T^*R$ by symplectomorphisms. More precisely, let $\Omega$ be the symplectic form on $\C^{2\ell}$ represented by the matrix
		$
		\begin{pmatrix}
			0 & \id_\ell\\
			-\id_\ell & 0
		\end{pmatrix}
		$, so that $\Sp(2\ell)$ consists of matrices $g\in \GL(2 \ell)$ such that $g^t \Omega g = \Omega$. Then $\Sp(2\ell )$ acts on $T^*R$ by treating $(X_1,\cdots, Y_\ell)$ as a vector. 
		This $\Sp(2\ell)$-action preserves the symplectic form on $T^*R$ by the choice of $\Omega$, and also preserves the subset $\tilde{\L}$. Now let $K_1,K_2 \in \PP^{2\ell -1}$; let $K_1', K_2'\in \C^{2\ell } \setminus 0$ be representatives of $K_1,K_2$. Choose $g\in \Sp(2\ell)$ such that $(K_1')^t g  = (K_2')^t$; then $\tilde{\L}_{K_2} = g \tilde{\L}_{K_1}$, and (3) follows. 
		
		It is easy to see that if $K_1\neq K_2\in \PP^{2\ell -1}$, then for any $i\neq j, i' \neq j'$, the hyperplanes $H_{K_1,i,j}$ and $H_{K_2,i',j'}$ are transversal. This proves (4). 
		
		(5) follows from (4). 
	\end{proof}
	For technical reasons, it will be convenient to name two special open affine subsets. 
	Let $I = [1,0,0,\cdots,0], J = [0,1,0,\cdots,0] \in \PP^{2\ell -1}$, so that $\tilde{\L}_I$ (resp. $\tilde{\L}_J$) consists of $s\in \tilde{\L}$ whose $X_1$ (resp. $X_2$) component has pairwise distinct diagonal entries. 
	\subsection{Model varieties}\label{subsection model varieties}
	In \cite{losev2006symplectic}, Losev constructed ``model varieties" for neighbourhoods of points in the Hamiltonian reduction of a symplectic variety by a reductive group. Here we adapt this idea and describe the neighbourhood of an open subset of a symplectic leaf. 
	
	Recall $G = \GL(n), \g = \gl(n)$; the stabilizer in $G$ of any element of $\tilde{\L}$ (viewed as a subset of $T^*R$) is the maximal torus $H \cong (\C^*)^n$, with Lie algebra $\cartanh = \C^n$. 
	Let $U = (\g/\cartanh)^*$; it is equipped with the coadjoint $H$-action. Furthermore, let $\C^*$ act on $U$ by $t.u = t^2 u$. 
	
	Let $V$ be a symplectic vector space of dimension $(n^2-n)(2\ell -2) + 2nw$, with symplectic form $\omega_V$. We identify $V$ with the symplectic vector space $T^*\sliceR$ in \Cref{subsection the slice quiver}. 
	Let 
	\begin{equation}\label{darboux basis for abstratc V}
		\left\{ 
		\begin{array}{c}
			v_{ij}^k,w_{ij}^k ; \\
			p_{\alpha,\beta},q_{\beta,\alpha} 
		\end{array}
		\left| 
		\begin{array}{c}
			1\le i \neq j \le n, 1\le k \le {\ell}-1 ;\\
			1\le \alpha \le n, 1\le \beta\le w
		\end{array}
		\right .
		\right\}
	\end{equation} 
	be the Darboux basis of $V = T^*\sliceR$ corresponding to the arrows of the extended double quiver of $\sliceQ$. Then 
	\begin{align}
		\begin{split}
			&\omega_V(v_{i,j}^k, v_{j',i'}^{k'} ) = \omega_V(w_{i,j}^k, w_{j',i'}^{k'} ) = \omega_V (p_{\alpha,\beta},p_{\alpha',\beta'} ) = \omega_V(q_{\beta,\alpha},q_{\beta',\alpha'}) = 0,\\
			&\omega_V(v_{i,j}^k, p_{\alpha,\beta} ) = \omega_V(v_{i,j}^k, q_{\beta,\alpha} ) =\omega_V(w_{j,i}^k, p_{\alpha,\beta} )  =\omega_V(w_{j,i}^k, q_{\beta,\alpha}) = 0,\\
			&\omega_V(v_{i,j}^k, w_{j',i'}^{k'} )=\delta_{i,i'}\delta_{j,j'}\delta_{k,k'},\\
			&\omega_V( p_{\alpha,\beta},q_{\beta',\alpha'}) = \delta_{\alpha,\alpha'}\delta_{\beta,\beta'}.
		\end{split}
	\end{align}
	
	The group $H$ acts on $V$ dual to the action \eqref{H action on algebra on TstarsliceR}. 
	We have two $\C^*$-actions on $V$, the model action and the conical action, dual to \Cref{definition of model torus action} and \Cref{definition of conical torus action}, respectively. We also have an $S_n$ action on $V$ permuting the indices $1\le i,j\le n$. 
	
	Following \cite{losev2006symplectic},  we define the homogeneous bundle
	\begin{equation}\label{definition of model homogeneous bundle}
		X = G \times^H(U\oplus V)\times \tilde{\L}.
	\end{equation}
	
	For $g\in G, u\in U,v\in V$ and $s\in \tilde{\L}$, we write $[g,u,v,s]\in X$ for the image of $(g,u,v,s)$ under the quotient by $H$. 
	
	We equip $X$ with a $G$-invariant symplectic form as follows. At $x = [e,u,v,s]$ where $e\in G$ is the identity element, $u\in U$, $v\in V, s\in \tilde{\L}$, consider the tangent space $T_xX = \g/\cartanh \oplus U \oplus V \oplus T_s \tilde{\L}$. The subspace $T_s \tilde{\L}$ is a symplectic vector space; we denote its symplectic form by $\omega_{\tilde{\L}}$. Let $\xi_i\in \g/\cartanh, u_i\in U, v_i\in V, s_i \in T_s\tilde{\L}$ for $i=1,2$. Then, $T_xX$ is the direct sum of three symplectic vector spaces $(\g/\cartanh \oplus U) \oplus V \oplus T_s \tilde{\L}$. More precisely, the symplectic form on $T_xX$ is given by  
	\begin{multline}
		\omega_x((\xi_1,u_1,v_1,s_1),(\xi_2,u_2,v_2,s_2))  = \\
		\<\xi_1,u_2\> - \<\xi_2,u_1\> + \omega_V(v_1,v_2) + \omega_{\tilde{\L}} (s_1,s_2) + \<\mu_H(v) + u, [\xi_1,\xi_2]\>
	\end{multline}
	which is skew symmetric, nondegenerate, and $H$-invariant. We then extend this form to all of $X$ by $G$-invariance. 
	The variety $X$ together with the 2-form $\omega_X$ is the direct product of the model variety in \cite[Proposition 1]{losev2006symplectic} and the symplectic variety $\tilde{\L}$. 
	The natural $G$-action on $X$ is Hamiltonian, and the moment map is 
	\begin{equation}
		\mu_G([g,u,v],s) = \Ad(g)(u+\mu_H(v)).
	\end{equation}
	
	Recall that $\C^*$ acts on $U$ by degree 2. Let $\C^*$ act on $V$ by the model action. The natural rescaling $\C^*$-action on $T^*R$ restricts to an action on $\wt{\L}$. We thus get a $\C^*$-action on $X$. Namely, for $\gamma\in \C^*$, 
	\begin{equation}\label{torus action on model variety}
		\gamma . [g,u,v,s] := [g,\gamma.u,\gamma.v, \gamma s].
	\end{equation} 
	We still call this action the \textit{model} $\C^*$-action. 
	
	It is straightforward to check this action is well defined, and that the form $\omega$ has degree 2 with respect to this action.  
	
	Now consider the affine opens $\tilde{\L}_K$ defined in \eqref{definition of tildeLK}. 
	Replace $\tilde{\L}$ with $\tilde{\L}_K$ in \eqref{definition of model homogeneous bundle}, we obtain $X_K$ which is an affine open subvariety of $X$, stable under the model $\C^*$-action. 
	We will construct a $G\times \C^*$-equivariant morphism $\phi_I: X_K \to T^*R$ extending the natural inclusion $[e,0,0, \tilde{\L}_K] = \tilde{\L}_K \hookrightarrow T^*R$. 
	Since all $\tilde{\L}_K$ are isomorphic (part (3) of \Cref{lemma properties of good affine open cover}), it suffices to consider the special open subset $\tilde{\L}_I$ at the end of \Cref{subsection good affine open cover}. 
	
	The symmetric group $S_n$ acts on $X_I$ in the following way. 
	\begin{itemize}
		\item First, view $S_n$ as a subgroup of $G$ consisting of permutation matrices, and let $S_n$ act on $G$ by $\sigma.g:=g\sigma\inverse$, where $g\in G, \sigma\in S_n$. 
		\item Let $S_n$ acts on $U = (\g/\cartanh)^*$ by adjoint action.
		\item 
		Let $S_n$ act on $V$ by permuting the Darboux basis. More precisely, for $\sigma \in S_n$ and basis element $v_{ij}^k$, we define $\sigma\cdot v_{ij}^k = v_{\sigma(i),\sigma(j)}^k$ and similarly for $w_{ij}^k$. 
		\item
		Finally, let $S_n$ act on $\tilde{\L}_I$ by permuting the diagonal entries of each of $X_1,\cdots,Y_{\ell}$. 
	\end{itemize}
	
	We get a well defined free action of $S_n$ on $X_I$ and a quotient variety $X_I/S_n := \Spec \C[X_I]^{S_n}$. 
	
	The model $\C^*$-action commutes with the $S_n$-action, so it descents to a $\C^*$-action on $X_I/S_n$. 
	Let $\omega_I$ denote the restriction of the degree 2 symplectic form $\omega_X$ to $X_I$; it is $S_n$ invariant, hence descends to a degree 2 symplectic form on $X_I/S_n$, which we still denote by $\omega_I$.
	
	The main result of this subsection is the following. 
	\begin{theorem}\label{affine open morphism before completion}
		There is a $G$-equivariant morphism of varieties
		\begin{equation}
			\phi_I: X_I/S_n \to (T^*R)_I
		\end{equation}
		such that 
		\begin{enumerate}
			\item $\phi_I$ is $\C^*$-equivariant; 
			\item for any $x\in X_I$ of the form $[e,0,0,s]$, the linear map $d_x \phi_{I}: T_xX_I \to T_s(T^*R)\cong T^*R$ is a symplectomorphism; 
			\item the restriction of $\phi_I$ to the closed subset $(G/H\times  \tilde{\L}_I)/S_n$ is an isomorphism onto its image $G\tilde{\L}_I$ in $T^*R$. Algebraically, $\phi_I$ induces an isomorphism of algebras
			\begin{equation}
				\C[(T^*R)_I]/\idealI_{G\tilde{\L}_I} \xrightarrow{\sim }
				\C[X_I/S_n]/\idealI_{(G/H \times \tilde{\L}_I)/S_n}.
			\end{equation}
			Here, $\idealI_{G\tilde{\L}_I}$ denotes the ideal of the closed subset $G\tilde{\L}_I$, and $\idealI_{(G/H \times \tilde{\L}_I)/S_n}$ denotes the ideal of $(G/H \times \tilde{\L}_I)/S_n$. 
		\end{enumerate}
	\end{theorem}
	The remaining part of this subsection will be devoted to constructing $\phi_I$. We start by defining an $S_n$-invariant morphism $X_I \to (T^*R)_I$, which we still call $\phi_I$.
	
	Define 
	\begin{equation}\label{basis of UI}
		U_I  = \left\{
		(X_1,...,Y_{\ell})\in  \sl_n^{\oplus 2\ell }\left| \begin{array}{c}
			X_1 = ... = X_{\ell} = Y_2 = ...=Y_{\ell} = 0;\\
			\text{diagonal of }Y_1 = 0
		\end{array}
		\right. \right\}. 
	\end{equation}
	It is easy to see that that $U_I$ is isomorphic to $(\g/\cartanh)^*$ as an $H$-module. 
	
	Fix $s\in \tilde{\L}_I$. Note that $\g.s \cap U_I = 0\in T^*R$. We define $V_{I,s}$ to be the skew orthogonal complement of $T_s \tilde{\L} \oplus \g.s  \oplus U_I $ in $T^*R$, which depends on $s$. 
	We will construct an $H$-equivariant linear isomorphism $\phi_{I,s}^U: U \to U_I $ and an $H$-equivariant symplectomorphism 
	$\phi_{I,s}^V: V\to V_{I,s}$, and then define 
	\begin{equation}\phi_I([g,u,v],s) = g.(s+ \phi_{I,s}^U(u) + \phi_{I,s}^V(v))\end{equation}
	which is clearly $G$-equivariant. 
	
	The maps $\phi_{I,s}^U$ and $\phi_{I,s}^V$ can be defined ``positionwise". More precisely, for each pair $1\le i\neq j \le n$, we consider the $2\ell $-dimensional symplectic vector space $R_{ij}$, spanned by the $(i,j)$ entries of $X_1,\cdots,X_{\ell}$, and its dual, $R_{ij}^*$, spanned by the $(j,i)$ entries of $Y_1,\cdots,Y_{\ell}$. 
	Also, let $R^p$ denote the space of framings $\Hom_\C (\C^w,\C^n)$.  
	We note that $T^*R = \bigoplus_{i\neq j} T^*R_{ij} \oplus T^*R^p$. 
	
	For any $i\neq j$, $U_I \cap T^*R_{ij} = (0,0,\cdots,0; \C  ,0\cdots,0)$, where the $\C$ lies in the $(l+1)$-th (corresponding to $Y_1$) position. 
	Let $e_{ij}$ be the elements of $\gl(n)$ with $1$ in the $(i,j)$-th entry and $0$ elsewhere. Let $s = (X_1(s),\cdots,Y_{\ell}(s))$. 
	Then 
	\begin{equation}
		\g.s \cap T^*R_{ij}  = \C\cdot (s_1,s_2,...,s_{\ell},t_1,t_2,...,t_{\ell}),
	\end{equation}
	where 
	\begin{align*}
		\begin{split}
			s_k&= X_k(s)_{ii} - X_k(s)_{jj} \\
			t_k & = Y_k(s)_{jj} - Y_k(s)_{ii}.
		\end{split}
	\end{align*}
	
	By our choice of $I$, $s_1 \neq 0$. Let us fix a basis element 
	\begin{equation}\label{basis of UIsij}
		u_{I,s,(ij)} = \left(0,\cdots,0;\frac{1}{s_1},0,\cdots,0\right)\in U_I\cap T^*R_{ij}
	\end{equation}
	so that when we view it as an element of $T^*R$, 
	\begin{equation}\omega_{T^*R}(e_{ij}.s, u_{I,s,(ij)}) = 1. \end{equation}
	As the position ranges over all $1\le i\neq j \le n$, the vectors $u_{I,s,(i,j)}$ form a basis of $U_I$ \textit{depending on} $s$. 
	
	Now we consider $V_{I,s}$, the skew orthogonal complement of $\g.s \oplus U_I$ in $T^*R$. Note that $T^*R^p \subset V_{I,s}$. We have the following ``positionwise decomposition" 
	\begin{equation}
		V_{I,s} = \left( \bigoplus_{1\le i \neq j \le n} (T^*R_{ij}\cap V_{I,s} )\right)\oplus T^*R^p.
	\end{equation}
	Therefore, we will describe its ``positionwise" components 
	\begin{equation}V_{I,s,(i,j)} := T^*R_{ij}\cap V_{I,s}.\end{equation}
	By definition, $V_{I,s,(i,j)}$ is exactly the skew orthogonal complement of $e_{ij}.s \oplus U_{I,(i,j)}$ in $T^*R_{ij} \cong \C^{2\ell }$. 
	Consider the following elements of $V_{I,s,(i,j)}$. 
	Note the difference in degrees of entries of $v^1,w^1$ with other vectors; this difference indicates that the model $\C^*$-action is the correct action to use.
	\begin{align}\label{naive basis for I = X1X1X1}
		\begin{split}
			v_{I,s,(ij)}^1 & = (0,\frac{1}{s_1},0,\cdots,0; \frac{t_2}{s_1^2},0,\cdots,0) \\
			v_{I,s,(ij)}^2 & =(0,0,1,\cdots,0; \frac{t_3}{s_1},0,\cdots,0) \\
			& \vdots \\
			v_{I,s,(ij)}^{{\ell}-1} & = (0,0,\cdots,1; \frac{t_{\ell}}{s_1},0,\cdots,0)  \\
			w_{I,s,(ij)}^1 & = (0,\cdots,0;-{s_2},{s_1},0,\cdots,0)\\
			w_{I,s,(ij)}^{2} & = (0,\cdots,0;-\frac{s_3}{s_1},0,1,0,\cdots,0)\\
			& \vdots \\
			w_{I,s,(ij)}^{{\ell}-1} & = (0,\cdots,0;-\frac{s_{\ell}}{s_1},0,0,\cdots,1).
		\end{split}
	\end{align}
	In the above notation, there are $\ell$ entries before ``;" corresponding to $X_1,\cdots X_{\ell}$ and $\ell$ entries after ``;" corresponding to $Y_1,\cdots Y_{\ell}$.
	
	Moreover, let us write $p_{s, (i,\beta)}, q_{s,(\beta,i)}$, $1\le i\le n,1\le \beta\le w$, for the entries of $p,q\in T^*R^p$. Then, 
	\begin{equation}\label{darboux basis for VIs}
		\left\{
		\begin{array}{c}
			v^k_{I,s,{ij}}, w^k_{I,s,{ij}},\\ 
			p_{s, (i,\beta)}, q_{s,(\beta,i)}
		\end{array} 
		\left| 
		\begin{array}{c}
			1\le i \neq j \le n, 1\le k\le {\ell}-1  \\
			1\le i \le n, 1\le \beta\le w
		\end{array}
		\right. 
		\right\}
	\end{equation}
	form a Darboux basis for $V_{I,s}$. The $H = (\C^*)^n$-action in terms of this basis is as follows. For $h=(t_1,\cdots, t_n)\in H$, we have 
	\begin{align}
		\begin{split}
			& h.v^k_{I,s,{ij}} = t_j t_i \inverse v^k_{I,s,{ij}}; \\ 
			&h.w^k_{I,s,{ij}} = t_jt_i\inverse w^k_{I,s,{ij}};\\
			&h.p_{s, (i,\beta)} = t_i p_{s, (i,\beta)}, \\ 
			&h.  q_{s,(\beta,i)}=t_i\inverse q_{s, (\beta,i)}. 
		\end{split}
	\end{align}
	Let us proceed to defining the morphism $\phi_I: X_I \to T^*R$. First, identify $U = (\g/\cartanh)^*$ with the subspace of $n\times n$ matrices with 0 diagonal. 
	Define the linear isomorphism 
	\begin{equation}
		\phi_{I,s}^U: U \to U_I, e_{ij} \mapsto u_{I,s,(i,j)},
	\end{equation}
	so that for any $u\in U,\xi\in \g$, we have
	\begin{equation}\<u,\xi\> = \omega_{T^*R}( \xi.s, \phi_{I,s}^U(u)).\end{equation}
	This morphism is $H$-equivariant. 
	
	Next, recall $V$ has Darboux basis \eqref{darboux basis for abstratc V}. We define a linear isomorphism
	\begin{align}\label{definition of phiIVs}
		\begin{split}
			\phi_{I,s}^V: V &\to V_{I,s}, \\
			v_{ij}^k &\mapsto v_{I,s,(i,j)}^k, \\
			w_{ij}^k &\mapsto w_{I,s,(i,j)}^k,\\
			p_{i,\beta} &\mapsto p_{s,(i,\beta)},\\
			q_{\beta,i}&\mapsto q_{s,(\beta,i)}.
		\end{split}
	\end{align}
	Since $\phi_{I,s}^V$ sends the Darboux basis \eqref{darboux basis for abstratc V} to the Darboux basis \eqref{darboux basis for VIs}, it is a symplectomorphism. Moreover, $\phi_{I,s}^V$ is $H$-equivariant with respect to the $H$-action \eqref{H action on algebra on TstarsliceR} on the source, and the restriction of the $G$-action to $H$ on the target.  
	
	Finally, we define the $G$-equivariant morphism
	\begin{align}\label{affine open morphism before completion formula}
		\begin{split}
			\phi_I:X_I &\to T^*R, \\
			[g,u,v,s] &\mapsto g.(s+\phi_{I,s}^U(u) + \phi_{I,s}^V(v)).
		\end{split}
	\end{align}
	The image of $\phi_I$ lies in $(T^*R)_I$ and $\phi_I$ is $S_n$-invariant, so it induces a morphism (by an abuse of notation) $\phi_I:X_I/S_n \to (T^*R)_I$. 
	We now show $\phi_I$ is our desired morphism. 
	\begin{proof}[Proof of \Cref{affine open morphism before completion}]
		For (1), note that the model $\C^*$-action is constructed so that $\phi_I$ is equivariant. For (2), note that when $x = ([e,0,0],s) \in X_I$, the linear map $d_x \phi_{I}: T_xX_I \to T_s(T^*R)\cong T^*R$ sends a Darboux basis to a Darboux basis. For (3), just restrict $\phi_I$ to the closed subset $u=0, v=0$. 
	\end{proof}
	
	Let $(X_I/S_n)^\wedge$ be the affine formal scheme obtained by completing $X_I/S_n$ along the closed subset $(G/H \times \tilde{\L}_I)/S_n$, and let $(T^*R)_I^\wedge$ be the affine formal scheme obtained by completing $(T^*R)_I$ along the closed subset $G\tilde{\L}_I$. 
	The morphism $\phi_{I}$ defined in \Cref{affine open morphism before completion} induces a morphism of formal schemes 
	\begin{equation}\label{definition of phiI complete}
		\phi_I^\wedge: X_I^\wedge/S_n\to (T^*R)_I^\wedge. 
	\end{equation}
	The morphism $\phi_I^\wedge$ is $\C^*$-equivariant and restricts to an isomorphism of the varieties of closed points on both sides, i.e. $(G/H \times \tilde{\L}_I)/S_n\cong G\tilde{\L}_I$. 
	
	For any $s\in G\tilde{\L}_I$, let $\mf{m}_s$ denote the ideal of $s$ in $A_I$ and $\mf{n}_s$ denote the ideal of $s$ in $\C[(T^*R)_I^\wedge]$. The morphism $\phi_I^\wedge$ induces an isomorphism between the completed rings 
	$A_I^{\wedge_{\mf{m}_s}}$ and $\C[(T^*R)_I^\wedge]^{\wedge_{\mf{n}_s}}$; 
	in fact, both are isomorphic to $\C[G\tilde{\L}_I]^{\wedge \mf{p}_s}\completeotimes \C[U\oplus V]^{\wedge_{\mf{m}_0}}$ where $\mf{p}_s$ is the ideal of $s$ in $\C[G\tilde{\L}_I]$.
	Therefore, $\phi_I^\wedge$ is an isomorphism of formal schemes. 
	
	We finish this subsection by the following description of generators of $\C[X_I]$. It will be used in \Cref{inner derivation}. 
	\begin{lemma}\label{lemma CXI generated by linear functions}
		The algebras $\C[X_I]$ is generated by $\C[G]^H, (\C[G] \otimes V^*)^H, (\C[G]\otimes U^*)^H$ and $\C[\tilde{\L}_I]$ as an algebra. Similar statement holds for $\C[X_J]$ and $\C[X_{I\cap J}]$. 
	\end{lemma}
	\begin{proof}
		Recall $\C[X_I] = (\C[G]\otimes (\C[U] \oplus \C[V] ))^H \otimes \C[\tilde{\L}_I] =((\C[G]\otimes \C[U])^H \oplus (\C[G]\otimes \C[V])^H) \otimes \C[\tilde{\L}_I]$. Therefore, it suffices to show that $(\C[G]\otimes \C[U])^H$ is generated by $\C[G]^H$ and $ (\C[G] \otimes U^*)^H$, and that $(\C[G]\otimes \C[V])^H$ is generated by $\C[G]^H$ and $ (\C[G] \otimes V^*)^H$. We show the first claim and the second is proved analogously. 
		
		Since $H$ is a torus, we have a weight space decomposition 
		\begin{equation}
			(\C[G]\otimes \C[U])^H = \bigoplus_{\chi \in \mathfrak{X}(H)} \C[G]_\chi\otimes \C[U]_{-\chi}
		\end{equation}
		where $\mathfrak{X}(H) \cong \ZZ^n$ denotes the character lattice of $H$, and $\bullet_\chi$ denotes the $\chi$-weight subspace. It suffices to show each $\C[G]_\chi\otimes \C[U]_{-\chi}$ is generated by $\C[G]$ and $(\C[G]\otimes U^*)^H$. 
		Denote by $u_{ij}$ the coordinate function of $U$ corresponding to the $(i,j)$-th entry. Let $F = f\otimes u_{i_1,j_1}u_{i_2,j_2}\cdots u_{i_k,j_k} \in \C[G]_\chi\otimes \C[U]_{-\chi}$, where $f\in \C[G]_\chi$.
		
		Note that the weight $wt(u_{ij})$ is 1 at position $i$ and $-1$ at position $j$, and 0 elsewhere. Note also that $wt(\det)=(1,1,\cdots,1)$. Therefore, $wt(\det^k F)\in\ZZ_{\ge 0}^n$. We may rewrite $\det^kF$ as a finite sum
		\begin{equation}
			{\det }^k F = \sum_{a\ge 1} \prod_{b=1}^{k} F_{ab}
		\end{equation}
		such that 
		\begin{equation}
			wt(F_{ab}) = -wt(u_{i_b,j_b}) + (1,1,\cdots,1). 
		\end{equation}
		Then, the decomposition
		\begin{equation}
			F =  \sum_{a\ge 1} \prod_{b=1}^{k} \frac{F_{ab} u_{i_b,j_b}}{\det} 
		\end{equation}
		shows that $F$ is generated by $ (\C[G] \otimes U^*)^H$. 
	\end{proof}
	\subsection{Moser's theorem on formal schemes}\label{subsection Moser's theorem}
	Recall the symplectic form $\omega_I$ on $X_I$ and also the closed 2-form $\phi_I^*\omega_{T^*R}$ on $X_I$. These two form are different in general. In this section, we prove an analogue of Moser's theorem on formal schemes, and then compare the restrictions of $\omega_I$ and $\phi_I^*\omega_{T^*R}$ to $X_I^\wedge/S_n$. 
	
	Let us recall the notion of differential forms on formal schemes. For our purpose, let $X$ be the total space of a vector bundle over a smooth variety $Y$, so that $X$ is also smooth. Let $\idealI\subset \sheafO_X$ be the ideal sheaf defining $Y$. 
	Let $\sheafO_{\hat{X}}$ denote the sheaf $\varprojlim_n \sheafO_X/\idealI^n$ on $Y$. 
	The next lemma follows from \cite[Corollary II.9.8]{hartshorne1977algebraic}. 
	\begin{lemma}
		The following are true. 
		\begin{enumerate}
			\item For any $k\ge 1$, there is an isomorphism of sheaves of $\sheafO_{\hat{X}}$-modules 
			\begin{equation}\varprojlim \Omega^k_X/\idealI^n \cong \sheafO_{\hat{X}}\otimes_{\sheafO_X} \Omega^k_X.\end{equation}
			In particular, $\varprojlim \Omega^k_X/\idealI^n$ is a locally free sheaf of $\sheafO_{\hat{X}}$-modules.
			\item Let $\hat{\idealI}$ be the $\idealI$-adic completion of $\idealI$. Then $\sheafO_X/\idealI^n \cong \sheafO_{\hat{X}}/(\hat{\idealI})^n$. 
		\end{enumerate}
	\end{lemma}
	
	Denote $\Omega_{\hat{X}}^k  = \varprojlim_n \Omega^k_X/\idealI^n \Omega^k_X$. The differentials $d: \Omega^k_X \to \Omega^{k+1}_X$ is continuous in the $\idealI$-adic topology, and therefore induces differentials $d: \Omega_{\hat{X}}^k \to \Omega_{\hat{X}}^{k+1}$. 
	\begin{lemma}
		Suppose $\phi$ is an automorphism of sheaf of algebras of $\sheafO_{\hat{X}}$ that preserves $ \hat{\idealI}$. Then $\phi$ induces an $\sheafO_{\hat{X}}$-linear automorphism of $\Omega_{\hat{X}}^k$. 
	\end{lemma}
	
	Since $X$ is the total space of a vector bundle on $Y$, there is a natural $\C^*$-action on $X$ given by dilation on the fibers, so that the linear functions along fibers have degree 1. 
	Let $\eu_F$ denote its Euler vector field ($F$ indicates the $\C^*$-action is fiberwise). Since The fixed points of this $\C^*$-action is precisely $Y$, we see $\eu_F \in H^0( X, \idealI  Vect(X))$. 
	
	\begin{lemma}\label{lemma for moser theorem difference of symp form is exact}
		Suppose $\omega \in H^0(X, \Omega_{{X}}^k)$ is a closed $k$-form such that $\omega_y=0$ for all $y\in Y$. 
		Then there is $\mu \in H^0(X, \idealI^2\Omega_{X}^{k-1})$ such that $\omega = d\mu$. 
	\end{lemma}
	\begin{proof}
		The $\C^*$-action on $X$ induces a $\ZZ_{\ge 0}$-grading on $\Omega_X^k$. Write $\omega = \sum_{i=0}^N \omega_i$ where $\omega_i$ has weight $i$. Since $\omega|_Y =0$, we see $\omega_0 = 0$ and each $\omega_i\in H^0(X,\idealI \Omega_X^k)$.
		Define 
		\begin{equation}
			\mu = \sum_{i=0}^{N} \frac{1}{i+1} \iota_{\eu_F}\omega_i. 
		\end{equation}
		This is an element of $H^0(X,\idealI^2\Omega_{X}^{k-1})$ since $\eu_F \in H^0( X, \idealI  Vect(X))$. 
		By Cartan magic formula, we see $d\mu = \omega$. 
	\end{proof}
	\begin{theorem}[Moser's theorem]\label{theorem exist auto intertwining two symp forms}
		Suppose $\omega_0,\omega_1$ are two symplectic forms on $X$ such that $(\omega_0-\omega_1)_y = 0 $ for all $y\in Y$. Then there is an automorphism $\Phi$ of $\sheafO_{\hat{X}}$ that preserves $\hat{\idealI}$, induces the identity on $\sheafO_X/\idealI$, and satisfies $\Phi^*\omega_1 = \omega_0$, viewing $\omega_0,\omega_1$ as elements of $\Omega_{\hat{X}}^2$. 
	\end{theorem}
	
	\begin{proof}
		By \Cref{lemma for moser theorem difference of symp form is exact}, there is $\mu\in H^0(X,\idealI^2\Omega^1_X)$ such that
		\begin{equation}\label{definition of 1form mu}
			d\mu = \omega:= \omega_1-\omega_0.
		\end{equation}
		Let $\omega_t:= \omega_0 + t \omega  $ for $0\le t\le 1$. Since $\omega_t$ is constant on $Y$, the degeneracy locus of $\omega_t$ is either empty or a divisor disjoint from $Y$. Let $X_t$ be an open subset of $X$ such that $\omega_t$ is nondegenerate on $X_t$; then $Y$ is a closed subvariety of $X_t$. For each $t$, there is $\eta^t\in Vect(X_t)$ satisfying
		\begin{equation}\label{definition of field etat}
			\iota_{\eta^t} \omega_t = -\mu|_{X_t}.
		\end{equation}
		Since $\mu\in \idealI^2\Omega^1$, we have 
		\begin{equation}\label{etat maps to Isquare}
			\eta^t \sheafO_{X_t}\subset \idealI^2 \sheafO_{X_t}. 
		\end{equation}
		In particular, $\eta^t|_Y = 0$. The vector field $\eta^t$ can be extended to a derivation of $\sheafO_{\hat{X}}$, and $\eta^t \sheafO_{\hat{X}} \subset \hat{\idealI}^2$. 
		
		We are going to define automorphisms $\Phi_t^*$ of $\sheafO_{\hat{X}}$ for $0\le t \le 1$, satisfying
		\begin{equation}\label{ode for flow}
			\frac{d}{dt} (\Phi_t^*f) = \eta^t (\Phi_t^*f)
		\end{equation} 
		and $\Phi_0^* = \id$. These conditions are motivated by the definition of flow of $\eta^t$ in the differential geometric setting. 
		
		To define $\Phi_t$, we first define a sequence of vector fields $\eta_n$ inductively. Set $\eta_0 = \eta^0$, i.e. $\iota_{\eta_0}\omega_0 = -\mu$. If $\eta_{n}$ is defined, define $\eta_{n+1}$ by the equation 
		\begin{equation}\label{inductively define eta n+1}
			\iota_{\eta_{n+1}}\omega_0 + \iota_{\eta_{n}}\omega = 0.\end{equation}
		Since $\omega \in \idealI\Omega^2_{\hat{X}}$, $\omega_0$ is nondegenerate, and $\eta_0  \sheafO_{\hat{X}}  \subset \hat{\idealI}^{2}$, it follows inductively that for any $n$, 
		\begin{equation}\label{etan converges}
			\eta_n \sheafO_{\hat{X}} \subset \hat{\idealI}^{n+2}.
		\end{equation}
		Next, we fix a section $f$ of $ \sheafO_{\hat{X}}$ and define a sequence $f_n \in \sheafO_{\hat{X}}$ inductively. Set $f_0 = f$; suppose $f_n$ is defined; define 
		\begin{equation}\label{inductively define fn+1 from previous fn}
			f_{n+1} = \frac{1}{n+1} \sum_{k=0}^n \eta_kf_{n-k}.
		\end{equation}
		Since $\eta_k\sheafO_{\hat{X}} \subset \hat{\idealI}^{k+2}$, we see by Leibniz rule $\eta_k\hat{\idealI}^{m} \subset \hat{\idealI}^{m+k+1}$. By induction, we conclude 
		\begin{equation}\label{f sub n in Taylor series of flow is convergent}
			f_n \in \hat{\idealI}^{n}.
		\end{equation}
		
		For $0\le t\le 1$, define 
		\begin{equation}\label{taylor series of Phitf}
			\Phi_t^*f = \sum_{n=0}^\infty t^nf_n.
		\end{equation}
		By \eqref{f sub n in Taylor series of flow is convergent}, $\Phi_t^*f $ is a well defined section of $\sheafO_{\hat{X}}$, and $\Phi_0^* = \id$. Let us show that  $\Phi_t^*$ satisfies \eqref{ode for flow}. 
		
		First, we show that for any local section $f$ of $\sheafO_{\hat{X}}$,
		\begin{equation}\label{taylor expansion of etat}
			\eta^t f= \sum_{n=0}^\infty t^n \eta_n f;
		\end{equation} the right hand side converges by \eqref{etan converges}. Since all $\omega_t$ are nondegenerate, it suffices to check that for any local vector field $u$ on $\hat{X}$, $\iota_{\eta^t}\omega_t (u) = \sum_{n=0}^\infty t^n \iota_{\eta_n}\omega_t (u)$; this follows directly from \eqref{inductively define eta n+1}. 
		
		Now by \eqref{taylor expansion of etat} and \eqref{inductively define fn+1 from previous fn}, we see $\Phi^*_t$ defined by \eqref{taylor series of Phitf} satisfies \eqref{ode for flow}. Note that \eqref{ode for flow} together with $\Phi_0^* = \id$ is precisely the definition of the flow of $\eta^t$ in the analytic setting, so $\Phi_t^*\omega_t = \omega_0$. In particular, $\Phi_1^*\omega_1 = \omega_0$, as desired. 
	\end{proof}
	In our situation, let $X = X_I/S_n, Y = (G/H \times \tilde{\L}_I)/S_n$, so $X$ is the total space of a vector bundle of rank $\dim U + \dim V$ over $Y$. Let $\omega_0 = \omega_I$ and $\omega_1 = \phi_I^*\omega_{T^*R}$. 
	By (2) of \Cref{affine open morphism before completion}, $\omega_1-\omega_0$ vanishes on $Y$. 
	Therefore, by \Cref{theorem exist auto intertwining two symp forms}, we obtain an automorphism $\Phi_I$ of $X^\wedge_I/S_n$ such that $\Phi_I^* (\phi_I^\wedge)^*\omega_{T^*R} = \omega_I$. 
	
	Moreover, note that $\omega_0,\omega_1$ both have degree 2 with respect to the model $\C^*$-action on $X$ (different from the fiberwise dilation action in the proof of \Cref{lemma for moser theorem difference of symp form is exact}). Therefore, we can take the 1-form $\mu$ in \eqref{definition of 1form mu} to be homogeneous of degree 2, and hence take the vector fields $\eta^t$ in \eqref{definition of field etat} to be $\C^*$-invariant. Since $G$ is reductive, we can further take $G$-invariant parts of $\mu$ and $\eta^t$. It follows that $\Phi_I$ is $G\times \C^*$-equivariant. 
	
	In summary, we obtain the following result. 
	\begin{proposition}\label{the final affine symplectic isomorphism from XI to TstarRI}
		The morphism 
		\begin{equation}\label{symplectic isomorphism between formal schemes before twisting}
			\phi_I^\wedge
			\circ \Phi_I : X_I^\wedge/S_n\to (T^*R)_I^\wedge
		\end{equation}
		is a $G\times\C^*$-equivariant symplectic isomorphism of formal schemes. It induces an isomorphism on the varieties of closed points. 
	\end{proposition}
	
	\subsection{Unipotency of transition map}\label{Unipotency of transition map}
	In this subsection we compare 
	the symplectomorphisms $\phi_I^\wedge
	\circ \Phi_I $ and $\phi_J^\wedge
	\circ \Phi_J $ on the open subsets ${\L}_I$ and ${\L}_J$. 
	
	Let us first define a morphism $\phi_J$ with analogous properties to $\phi_I$ in \Cref{affine open morphism before completion}. 
	Recall the Darboux basis \eqref{naive basis for I = X1X1X1} for $V_{I,s,(1,2)}$; we define similarly 
	\begin{align}\label{naive basis for J = X2X2X2}
		\begin{split}
			v_{J,s,(1,2)}^1 & = (-\frac{1}{s_2},0,0,\cdots,0; 0,-\frac{t_1}{s_2^2},\cdots,0) \\
			v_{J,s,(1,2)}^2 & =(0,0,1,\cdots,0;0, \frac{t_3}{s_2},\cdots,0) \\
			& \vdots \\
			v_{J,s,(1,2)}^{{\ell}-1} & = (0,0,\cdots,1;0, \frac{t_{\ell}}{s_2},\cdots,0)  \\
			w_{J,s,(1,2)}^1 & = (0,\cdots,0;-s_2, {s_1},0,\cdots,0)\\
			& \vdots \\
			w_{J,s,(1,2)}^{{\ell}-1} & = (0,\cdots,0;0,-\frac{s_{\ell}}{s_2},0,\cdots,1)
		\end{split}
	\end{align}
	where there are $\ell$ entries before ``;" corresponding to $X_1,\cdots,X_{\ell}$ and $\ell$ entries after ``;" for $Y_1,\cdots Y_{\ell}$. 
	This is a Darboux basis for $V_{J,s,(1,2)}$. Note that we made a choice so that $w_{J,s,(1,2)}^1 = w_{I,s,(1,2)}^1$, so there is an extra minus sign in the definition of $v_{J,s,(1,2)}^1 $ compared to $v_{I,s,(1,2)}^1 $. 
	
	Similar to \eqref{basis of UIsij}, let 
	\begin{equation}
		u_{J,s,(1,2)} = (0,\cdots,0;0, \frac{1}{s_2},0,\cdots,0).
	\end{equation} 
	It is dual to $e_{12}.s$ in $U_{J,(1,2)}$. The morphism $\phi_J$ and $\phi_J^\wedge$ are defined similarly to \eqref{affine open morphism before completion formula}, but using the Darboux basis \eqref{naive basis for J = X2X2X2} for $\phi_{J,s}^V$ and using $u_{J,s,(i,j)} $ for $\phi_{J,s}^U$. 
	
	In the remainder of this paper, we often use the notation $\bullet_{I\cap J} = \bullet_{I}\cap \bullet_J, \bullet_{I\cup J} = \bullet_{I}\cup \bullet_J$. For example, $\tilde{\L}_{I\cap J} = \tilde{\L}_I\cap \tilde{\L}_J$.
	Let $X_{I\cap J}^\wedge$ denote the affine formal scheme obtained by completing $X_{I\cap J}$ along the closed subvariety $G/H \times \tilde{\L}_{I\cap J}$. 	
	We can restrict the morphisms $\phi_I,\phi_J$ to (different) isomorphisms $X_{I\cap J}^\wedge/S_n \to (T^*R)_{I\cap J}$, and restrict the automorphisms $\Phi_I$ and $ \Phi_J$ to $X_{I\cap J}^\wedge/S_n$. 
	
	By the constructions in the previous subsection, we have the following diagram 
	\begin{equation}\label{comm diagram defining theta}
		\begin{tikzcd}
			X_{I\cap J}^\wedge/S_n \arrow[r, "\Phi_I"] \arrow[dd, dashed ]{}{\theta_{IJ}}  & X_{I\cap J}^\wedge/S_n \arrow[dr, "\phi_I^\wedge"] \arrow[dd, dashed ]{}{\phi_{IJ}} \\
			& & (T^*R)^\wedge_{I\cap J} \arrow[dl, "(\phi_J^\wedge)\inverse"]\\
			X_{I\cap J}^\wedge/S_n & X_{I\cap J}^\wedge/S_n \arrow[l, "\Phi_J\inverse"]
		\end{tikzcd}
	\end{equation}
	where we define the automorphisms $\theta_{IJ}$ and $\phi_{IJ} $ so that the diagram commutes. By \Cref{the final affine symplectic isomorphism from XI to TstarRI}, $\theta_{IJ}$ is $\C^*$-equivariant and preserves the symplectic form $ \omega_{{I\cap J}}$. 
	By (3) of \Cref{affine open morphism before completion}, $\phi_{IJ}$ restricts to the identity on $A_{I\cap J} / \idealI_{G\tilde{\L}_{I\cap J}}$. 
	
	Fix $x = [e,0,0,s]$ and consider the induced linear automorphism $d_x\phi_{IJ}$ on the tangent space $T_x(X_{I\cap J}^\wedge/S_n) \cong \g/\cartanh \oplus U \oplus V \oplus T_s\tilde{\L}$. 
	Let us write down the matrix of $d_x\phi_{IJ}$. Note first that $d_x\phi_{IJ}|_{T_x\tilde{\L}}$ is the identity. Moreover, fix $(i,j)$ with $1\le i \neq j\le n$, we see by the construction of $\phi_I$ and $\phi_J$ that $d_x\phi_{IJ}$ stabilizes the subspace 
	\begin{equation}\label{ij subspace in X}
		\C.(e_{ij}) \oplus \C. (e_{ij}^*) \oplus \Span \{v_{ij}^k,w_{ij}^k\}_{k=1}^{{\ell}-1}
	\end{equation}
	\begin{lemma}
		For any $i\neq j$, $d_x\phi_{IJ}$ is unipotent on the subspace \eqref{ij subspace in X}.
	\end{lemma}
	\begin{proof}
		Without loss of generality let $(i,j) = (1,2)$. Fix the ordered basis 
		\begin{equation}\label{ordered basis of 12position}
			\{e_{12}, w_{1,2}^1,w_{1,2}^2,\cdots,w_{1,2}^{{\ell}-1},v_{1,2}^{{\ell}-1},\cdots, v_{1,2}^1,(e_{12}^*)\}\end{equation}
		of the subspace \eqref{ij subspace in X}. 
		By the definition of $\phi_I,\phi_J$,
		\begin{equation}\label{image of basis under dxphiI}
			\begin{split}
				d_x\phi_I (e_{12}) &= (s_1,s_2,\cdots,t_\ell)\\
				d_x\phi_I (v_{1,2}^k) &= v_{I,s,(1,2)}^k  \\
				d_x\phi_I (w_{1,2}^k) &= w_{I,s,(1,2)}^k \\
				d_x\phi_I(e_{12}^*) &= u_{I,s,(1,2)}
			\end{split}
		\end{equation}
		and 
		\begin{equation}\label{image of basis under dxphiJ}
			\begin{split}
				d_x\phi_J (e_{12}) &= (s_1,s_2,\cdots,t_\ell)\\
				d_x\phi_J (v_{1,2}^k) &= v_{J,s,(1,2)}^k \\
				d_x\phi_J (w_{1,2}^k) &= w_{J,s,(1,2)}^k \\
				d_x\phi_J(e_{12}^*)& = u_{J,s,(1,2)}.
			\end{split}
		\end{equation}
		To compute the matrix of the restriction of $d_x\phi_{IJ}$ to the subspace \eqref{ij subspace in X} in the basis \eqref{ordered basis of 12position}, 
		it suffices to express the vectors in the right hand side of \eqref{image of basis under dxphiJ} as linear combinations of those of \eqref{image of basis under dxphiI}. One computes that in the basis \eqref{ordered basis of 12position}, $d_x\phi_{IJ}$ takes the following matrix form: 
		\begin{equation}\label{matrix for differential of phiIJ}
			\left(
			\begin{array}{c|c|cccc|cccc|c|c}	
				1 &   &                     &                     &        &                        &                       &                           &        &                    & -\frac{1}{s_1s_2}                                   &                   \\ \hline
				& 1 & -\frac{s_3}{s_1s_2} & -\frac{s_4}{s_1s_2} & \cdots & -\frac{s_\ell}{s_1s_2} & \frac{t_\ell}{s_1s_2} & \frac{t_{\ell-1}}{s_1s_2} & \cdots & \frac{t_3}{s_1s_2} & \frac{1}{s_1s_2}(\frac{t_2}{s_1} - \frac{t_1}{s_2}) & -\frac{1}{s_1s_2} \\ \hline
				&   & 1                   &                     &        &                        &                       &                           &        &                    & \frac{t_3}{s_1s_2}                                  &                   \\
				&   &                     & 1                   &        &                        &                       &                           &        &                    & \frac{t_3}{s_1s_2}                                  &                   \\
				&   &                     &                     & \ddots &                        &                       &                           &        &                    & \vdots                                              &                   \\
				&   &                     &                     &        & 1                      &                       &                           &        &                    & \frac{t_\ell}{s_1s_2}                               &                   \\ \hline
				&   &                     &                     &        &                        & 1                     &                           &        &                    & \frac{s_\ell}{s_1s_2}                               &                   \\
				&   &                     &                     &        &                        &                       & 1                         &        &                    & \frac{s_{\ell-1}}{s_1s_2}                           &                   \\
				&   &                     &                     &        &                        &                       &                           & \ddots &                    & \vdots                                              &                   \\
				&   &                     &                     &        &                        &                       &                           &        & 1                  & \frac{s_3}{s_1s_2}                                  &                   \\ \hline
				&   &                     &                     &        &                        &                       &                           &        &                    & 1                                                   &                   \\ \hline
				&   &                     &                     &        &                        &                       &                           &        &                    &                                                     & 1
			\end{array}
			\right)
		\end{equation}
		where all empty entries are 0. The lemma follows since the matrix is strictly upper triangular. 
	\end{proof}
	Note that $d_x\phi_{IJ}$ is the identity on the subspace of $V$ spanned by $p_{\alpha,\beta}, q_{\alpha,\beta}$, we obtain the following result. 
	\begin{proposition}\label{differeitnal of phiIJ is unipotent}
		For $x = [g,0,0,s]$, the differential $d_x\phi_{IJ}$ is unipotent on the whole space $T_x(X_{I\cap J}^\wedge/S_n)$. 
	\end{proposition}
	
	\subsection{Convergence of logarithm}
	In this subsection we establish the existence of logrithm the automorphism $\theta_{IJ}$ on $X_{I\cap J}^\wedge/S_n$ defined in \eqref{comm diagram defining theta}. We will show the existence of $\log \Phi_I^*$, $\log \Phi_J^*$ and $\log \phi_{IJ}^*$ separately, and then apply the BCH formula. 
	
	Consider the algebra morphism $\Phi_I^*$ on  $A_I = \C[X_I^\wedge/S_n]$, and the ideal $\idealI_I \subset A_I$  corresponding to $G\tilde{\L}_I$. 
	We show $\log \Phi_I^*$ is well defined. More precisely, 
	\begin{proposition}
		For any $i\ge 0$, 
		\begin{equation}\label{pronilpotency of PhiI-Id}
			(\Phi_I^*-\id)\idealI_I^i \subset \idealI_I^{i+1} .
		\end{equation}
		In particular, the series
		\begin{equation}\label{definition of log PhiI}
			\log \Phi_I^* : =\sum_{i=1}^\infty (-1)^{i+1} \frac{(\Phi_I^* - \id)^i}{i}
		\end{equation}
		is a well defined derivation of $\C[X_I^\wedge/S_n]$. 
		
	\end{proposition} 
	\begin{proof}
		We prove the following slightly stronger statement: for all $t\in [0,1]$, 
		\begin{equation}\label{condition for taking log}
			(\Phi_{I,t}^*-\id) \idealI_I^i \subset \idealI_I^{i+1}
		\end{equation}
		Recall $\Phi_{I,t}$ is the integral flow of the vector fields $\eta^t$, defined in the proof of \Cref{theorem exist auto intertwining two symp forms}. 
		
		When $t = 0$, $\Phi_{I,0} = \id$, so \eqref{condition for taking log} holds for $t=0$ and arbitrary $i$. By the definition of flow, 
		\begin{equation}\frac{d}{dt}(\Phi_{I,t}^*-1)f =\frac{d}{dt}\Phi_{I,t}^*f= \eta^t(\Phi_{I,t}^*f).\end{equation}
		By \eqref{etat maps to Isquare}, the right hand side lies in $\idealI_I^2$. Therefore, \eqref{condition for taking log} holds for $i=0,1$ and $t\in [0,1]$. Now suppose \eqref{condition for taking log} holds for $i$. Suppose $f \in \idealI_I^{i+1}$, and assume without loss of generality that $f = f_ig$ where $f_i\in \idealI_I^i$ and $g\in\idealI_I$. Then 
		\begin{equation}    (\Phi_{I,t}^*-1)f =((\Phi_{I,t}^*-1)f_i)g  +  (\Phi_{I,t}^*f_i) (\Phi_{I,t}^*-1)(g)\end{equation}
		where both terms on the right hand side lies in $\idealI_I^{i+2}$ by induction hypothesis. Condition \eqref{condition for taking log}, and hence the proposition, is proved. 
	\end{proof}
	The same argument also shows $\log \Phi_J^*$ is well defined. Moreover, we see from \eqref{pronilpotency of PhiI-Id} that for any $i\ge 0$, 
	\begin{equation}\label{pronilpotency of log PhiI}
		\log \Phi_I^* (\idealI_I^i) \subset \idealI_I^{i+1}, \log \Phi_J^* (\idealI_J^i) \subset \idealI_J^{i+1}.
	\end{equation}
	
	Now let us turn to $\phi_{IJ}$. Let $\phi_{IJ}^*$ denote the correponding automorphism on $A_{I\cap J} 
	$. 
	Thanks to \Cref{differeitnal of phiIJ is unipotent}, we see there is a large enough $N$ such that for any $k\ge 0$,
	\begin{equation}\label{explicit nilpotency of phiIJ}
		(\phi_{IJ}^*-1)^N \idealI_{I\cap J}^k \subset \idealI_{I\cap J}^{k+1}
	\end{equation}
	so the linear map $\phi_{IJ}^*-1$ is pro-nilpotent. Therefore, $\log \phi_{IJ}^*$ is well defined.

	Finally, let us consider $\log\theta_{IJ}^* = \log( \Phi_I^*  \circ \phi_{IJ}^*\circ (\Phi_J^*)\inverse) $. Here, we abuse notation and write $\Phi_I,\Phi_J$ for their restrictions on $X_{I\cap J}^\wedge/S_n$. We will show $\log\theta_{IJ}^*$ is a well defined derivation on $A_{I\cap J}$ via the Baker–Campbell–Hausdorff (BCH) formula. 
	
	Let us first examine $\log( \Phi_I^*  \circ \phi_{IJ}^*)$. By BCH formula, 
	\begin{equation}\label{BCH for two operators}
		\log( \Phi_I^*  \circ \phi_{IJ}^*) = \log  \Phi_I^* + \log \phi_{IJ}^* + \frac{1}{2}[\log  \Phi_I^*, \log \phi_{IJ}^*] + \text{higher order brackets} .
	\end{equation}
	By \eqref{explicit nilpotency of phiIJ}, we see there is a large enough $N$ such that $(\log \phi_{IJ}^*)^N \idealI_{I\cap J}^i \subset \idealI_{I\cap J}^{i+1} $. Also, by \eqref{condition for taking log} we see $(\log \Phi_I^*) \idealI_{I\cap J}^i\subset \idealI_{I\cap J}^{i+1}$. These observations imply that, for any $k\ge 1$, if $\psi$ is an order $kN$ bracket between  $\log \phi_{IJ}^*$ and $\log \Phi_I^*$, then $\psi \idealI_{I\cap J}^i\subset \idealI_{I\cap J}^{i+k}$. 
	This implies the operator defined by \eqref{BCH for two operators} is convergent on $A_{I\cap J}$, and moreover, 
	\begin{equation}
		(\log( \Phi_I^*  \circ \phi_{IJ}^*))^N \idealI_{I\cap J}^i\subset \idealI_{I\cap J}^{i+1}.
	\end{equation}
	Now apply the same argument to the pair of operators $(\Phi_J^*)\inverse $ and $\Phi_I^*\circ\phi_{IJ}^*$, we conclude that $\log\theta_{IJ}^*$ is a derivation on $A_{I\cap J}$. 
	Moreover, 
	there exists a large enough $N$, such that for all $i\ge 0$, 
	\begin{equation}\label{log of theta is pro nilpotent}
		(\log\theta_{IJ}^*)^N \idealI_{I\cap J}^i = \log( \Phi_I^*  \circ \phi_{IJ}^*\circ (\Phi_J^*)\inverse)^N \idealI_{I\cap J}^i
		\subset\idealI_{I\cap J}^{i+1}.
	\end{equation}
	Since $\theta_{IJ}$ is a symplectic automorphism, it follows that $\log\theta_{IJ}^*$ is a symplectic vector field. 
	
	\subsection{Pro-nilpotency of logarithm}
	In this subsection we give a more precise description for the pro-nilpotency of $\log \phi_{IJ}^*$ and $\log \theta_{IJ}^*$. 
	Define the following subspaces of $V$: 
	\begin{align}
		&V^k =
		\begin{cases}
			\Span_\C\{w^k_{ij}| 1\le i\neq j \le n, \}, &1\le k \le {\ell}-1\\
			\Span_\C\{v^{k-\ell+1}_{ij}| 1\le i\neq j \le n\},  & l\le k \le 2\ell -2
		\end{cases}\\
		&V^\mathbf{p}= \Span_\C\{p_{\alpha,\beta}|1\le \alpha\le n, \ell \le \beta\le w\},\\
		&V^\mathbf{q} =  \Span_\C\{q_{\alpha,\beta}|1\le \alpha\le n, 1\le \beta\le w\}.
	\end{align}
	Each of the subspaces $V^k,V^\mathbf{p},V^\mathbf{q}$ and $U$ is stable under the $H $ and the $ S_n$-actions. 
	
	The $S_n$-action on $X_{I\cap J}$ is free, so the derivations $\log \phi_{IJ}^*$ and $\log \theta_{IJ}^*$ lift to derivations of $\C[X_{I\cap J}]$. 
	Set $B: = G\times \tilde{\L}_{I\cap J}$. 
	We can rewrite
	\begin{equation}
		\C[X_{I\cap J}^\wedge ] = \left(\C[B] \completeotimes \C[U]^{\wedge 0} \completeotimes
		\left( \widehat{\bigotimes} _{k=1}^{2\ell -2}\C[V^k]^{\wedge 0} \right) 
		\completeotimes \C[V^\mathbf{p}]^{\wedge 0}\completeotimes \C[V^\mathbf{q}]^{\wedge 0} \right)^{H}. 
	\end{equation}
	
	The matrix form \eqref{matrix for differential of phiIJ} of $d\phi_{IJ}$ implies the following property of the derivation $\log \phi_{IJ}^*$. Recall the generators of $\C[X_{I\cap J}]$ in \Cref{lemma CXI generated by linear functions}. 
	\begin{proposition}\label{pronilpotency for derivation log phiIJ}
		The following are true. 
		\begin{enumerate}
			\item If $f \in \C[B]$, then $\log \phi_{IJ}^* (f) \in  \idealI_{I\cap J}$.
			\item If
			$f\in (\C[G]\otimes (V^k)^*)^H$, then, 
			\begin{multline}
				\log \phi_{IJ}^* (f) \in \C[B] \otimes (V^{k+1})^* + \C[B] \otimes (V^{k+2})^* \\ + \cdots +\C[B] \otimes (U)^* + \idealI_{I\cap J}^2.
			\end{multline}
			\item If 
			$f\in \C[B]\otimes (U^* \oplus  (V^\mathbf{p})^* \oplus (V^\mathbf{q})^*)$. Then 
			\begin{equation}
				\log \phi_{IJ}^* (f) \in \idealI_{I\cap J}^2.
			\end{equation}
		\end{enumerate}
	\end{proposition}
	Recall the pro-nilpotency condition \eqref{pronilpotency of log PhiI} for $\log\Phi_I^*$ and $\log \Phi_J^*$. By the BCH formula, $\log \theta^*_{IJ}$ is an infinite series in $\log\Phi_I^*, \log\Phi_J^*$, $\log \phi_{IJ}^*$, and brackets among them. 
	We conclude the following result. 
	\begin{proposition}\label{pronilpotency for delta}
		The statements (1)-(3) of \Cref{pronilpotency for derivation log phiIJ} holds with $\log \theta^*_{IJ}$ replacing $\log \phi^*_{IJ}$. 
	\end{proposition}
	\subsection{Inner derivation}\label{inner derivation}
	In this subsection we show that, up to some modification, the symplectic vector field $\log\theta_{IJ}^*$ is in fact Hamiltonian. More precisely, we will prove the following result. 
	\begin{proposition}\label{proposition trasition is inner}
		There are $G\times \C^*$-equivariant symplectomorphisms $\Psi_I$ of $X_I^\wedge/S_n$ and $\Psi_J$ of $X_J^\wedge/S_n$, and a $G$-invariant element $f\in A_{I\cap J}$ of degree 2, 
		such that 
		\begin{enumerate}
			\item $\Psi_I,\Psi_J$ restrict to the identity morphism on $(G\times^H \tilde{\L}_I)/S_n$ and $(G\times^H \tilde{\L}_J)/S_n$ respectively. 
			\item $\Psi_I,\Psi_J$ restrict to automorphisms of $X_{I\cap J}^\wedge/S_n$; 
			\item as an automorphism of $A_{I\cap J}$, $ (\Psi_J\inverse\theta_{IJ}\Psi_I )^*= \exp\{f,\bullet\}.$
		\end{enumerate}
	\end{proposition}
	\begin{remarklabeled}
		Most of the rest of this subsection is devoted to the proof of \Cref{proposition trasition is inner}, which is motivated by \cite[Section 4.1]{losev2021harish}. Let us briefly explain the idea of the proof. 
		
		First, we will use a de Rham cohomology vanishing argument to show that (\eqref{derivation version of coboundary condition}) $\delta = \delta_I - \delta_J + \{f,-\}$ for Poisson derivations $\delta_I$ of $\C[X_I^\wedge]$, $\delta_J$ of $\C[X_J^\wedge]$, and $f\in \C[X_{I\cap J}^\wedge]$. We would like to let $\Phi_I = \exp(\delta_I),\Phi_J = \exp(\delta_J)$, but these operators do not a priori converge. 
		
		Therefore, we need to modify $\delta_I, \delta_J$ and $f$ so that the above operators are well defined. This is done by giving a multigrading to the algebra $\C[X_{I\cap J}^\wedge]$, and apply \Cref{pronilpotency for delta}. The most technical part of the proof, \Cref{lemma truncated derivations are poisson}, is to show that the modified versions of $\delta_I, \delta_J$ are still Poisson derivations. 
	\end{remarklabeled}
	
	We start by defining the de Rham cohomology of the formal scheme $(X_{I\cup J}/S_n)^\wedge$, similar to the definition of de Rham cohomology of an algebraic variety. 
	
	Namely, consider the affine formal schemes $(X_I/S_n)^\wedge$ and $(X_J/S_n)^\wedge$ that cover $(X_{I\cup J}/S_n)^\wedge$. 
	Denote by 
	$\Omega^k_I = \Omega^k_{(X_I/S_n) /\C}$, which is a module over $A_I = \C[X_I/S_n]$. Denote by $A_I^\wedge$ the completion of $\C[X_I/S_n]$ along the ideal $\idealI_I$ of the closed subset $(G/H \times \tilde{\L}_I)/S_n$. Define $(\Omega^k_I)^\wedge$ to be the $\idealI_I$-adic completion of $\Omega^k_I$. Similarly, define $\Omega^k_J, A_J, A_J^\wedge, (\Omega^k_J)^\wedge$. 
	We can write down the following \v{C}ech-de Rham complex 
	\begin{equation}\label{cech derham complex formal version}
		\begin{tikzcd}[sep=small]
			\vdots  &  \vdots &   \\ 
			(\Omega^2_{I})^\wedge \oplus (\Omega^2_{J})^\wedge \arrow[r] \arrow[u]  &  (\Omega^2_{I\cap J})^\wedge \arrow[r] \arrow[u] & 0 \\ 
			(\Omega^1_{I})^\wedge \oplus (\Omega^1_{J})^\wedge \arrow[r] \arrow[u]&  (\Omega^1_{I\cap J})^\wedge \arrow[r] \arrow[u]& 0 \\ 
			(\Omega^0_{I})^\wedge \oplus (\Omega^0_{J})^\wedge \arrow[r] \arrow[u]&  (\Omega^0_{I\cap J,n})^\wedge \arrow[r] \arrow[u]& 0 \\ 
		\end{tikzcd}
	\end{equation}
	\begin{definition}
		The de Rham cohomology of $(X_{I\cup J}/S_n)^\wedge$ is defined to be the cohomology of the total complex associated to the double complex \eqref{cech derham complex formal version}. The $i$-th cohomology is denoted by $H^i_{\dR}((X_{I\cup J}/S_n)^\wedge)$. 
	\end{definition}
	Note that $X_{I\cup J}/S_n$ is a vector bundle over $(G/H \times \L_{I\cup J}) / S_n$. We will therefore use the following algebraic version of Poincare lemma. Suppose 
	$\pi: E \to B$ is a vector bundle over a smooth variety. Let $E^\wedge$ denote the completion of $E$ along the zero section. 
	For any affine open cover $\{B_\alpha\}$ of $B$ and 
	$E_\alpha = \pi\inverse B_\alpha$ of $E$, we can write down a \v{C}ech-de Rham complex of $E^\wedge$ similar to \eqref{cech derham complex formal version} and define $H^i_{\dR}(E^\wedge)$. 
	\begin{lemma}\label{algebraic poincare lemma}
		We have $H^i_{\dR}(E^\wedge) = H^i_{\dR}(B)$. In particular, $H^i_{\dR}(E^\wedge)$ is independent of the choice of the cover $B_\alpha$. 
	\end{lemma}
	\begin{proof}
		This is an analogue of the analytic version of Poincare lemma. We first show 
		$H^i_{\dR}(E^\wedge) \cong H^i_{\dR}((E \times \C)^\wedge)$ (the latter is completed along the closed subset $B\times\C$) using the homotopy $K$ defined completely analogously in the beginning of \cite[Chapter 4.1]{bott1982differential}. Note that the natural maps 
		\begin{align}
			\begin{split}
				\pi: E^\wedge \to E \to B\\
				\iota: B \to E^\wedge
			\end{split}
		\end{align}
		satisfies that $\pi \circ\iota = \id_B$ and $\iota\circ\pi$ is homotopic to $\id_E$. An analogue of \cite[Corollary 4.1.2]{bott1982differential}, replacing $\RR$ by $\C$, shows that $\pi^*$ and $\iota^*$ are isomorphisms between $H^i_{\dR}(E^\wedge)$ and $H^i_{\dR}(B)$. 
	\end{proof}
	\begin{lemma}\label{lemma H2 of base is 0}
		$H^2_{\dR}((G/H \times \L_{I\cup J}) / S_n) = 0. $
	\end{lemma}
	\begin{proof}
		We first show 
		$H^i_{dR}(\tilde{\L}_{I\cup J}) = 0$ for $i=1,2$. 
		Note $\tilde{\L}_{I\cup J}$ is an open subset of $\C^d$, where $d =2\ell (n-1) $. Its complement $Z$ is a (singular) closed subset of $\C^d$ with pure dimension $d -2$. 
		There is a long exact sequence of Borel-Moore homology groups (\cite[(2.6.10)]{chriss1997representation})
		\begin{align}
			\begin{split}
				\cdots \to H_{2d-1}(\C^d) \to H_{2d-1}(\tilde{\L}_{I\cup J}) \to H_{2d-2}(Z) \\
				\to H_{2d-2}(\C^d) \to H_{2d-2}(\tilde{\L}_{I\cup J}) \to H_{2d-3}(Z) \to \cdots.
			\end{split}
		\end{align}
		By Poincare duality, the above sequence can be written as 
		\begin{align}
			\begin{split}
				\cdots \to H^1(\C^d) \to H^1(\tilde{\L}_{I\cup J}) \to H_{2d-2}(Z) \\
				\to H^2(\C^d) \to H^2(\tilde{\L}_{I\cup J}) \to H_{2d-3}(Z) \to \cdots,
			\end{split}
		\end{align}
		where $H^\bullet$ stands for singular cohomology with complex coefficients. Note that $H^1(\C^d) = H^2(\C^d) = 0$, and $ H_{2d-3}(Z) = H_{2d-2}(Z) = 0$ since $Z$ has real dimension $2d-4$. We conclude $H^1(\tilde{\L}_{I\cup J}) = H^2(\tilde{\L}_{I\cup J}) =0$. The algebraic de Rham theorem implies $H^i_{dR}(\tilde{\L}_{I\cup J}) = 0$ for $i=1,2$.

		It is classical that $H^2(G/H,\C) = H^2(G/B,\C) = \C^{n-1}$. Note that $G/H \times \tilde{\L}_{I\cup J} \to (G/H \times \tilde{\L}_{I\cup J})/S_n$ is a covering map. By Kunneth formula, 
		\begin{equation}
			H^2_{dR}((G/H \times \tilde{\L}_{I\cup J})/S_n) =  (\C^{n-1})^{S_n} = 0.
		\end{equation}
		This finishes the proof of the lemma.
	\end{proof}
	By \Cref{algebraic poincare lemma} and \Cref{lemma H2 of base is 0}, we conclude
	\begin{equation}\label{second derham formal vanish}
		H^2_{dR}(X_{I\cup J }^\wedge/S_n) =0.
	\end{equation}
	Denote $\delta = \log\theta_{IJ}^*$. 
	Let $\beta = \iota_{\delta}\omega_{IJ}$ be the 1-form on $X^\wedge_{I\cap J}/S_n$, where $\iota$ denotes contraction. It is a $G$-invariant degree 2 element in $(\Omega_{I\cap J}^1)^\wedge$. 
	Since $\delta$ is a symplectic vector field, the 2-cochain $(\beta,0)$ is a cocycle. By \eqref{second derham formal vanish}, $\beta$ is in fact a coboundary, which means (recall the complex \eqref{cech derham complex formal version}) there are closed 1-forms 
	\begin{equation}\beta_I \in (\Omega_I^1)^\wedge,
		\beta_J \in (\Omega_J^1)^\wedge,
	\end{equation}
	and a function $f\in A_{I\cap J}$ such that
	\begin{equation}\label{beta is a coboundary}
		\beta = \beta_I|_{X_{I\cap J}^\wedge/S_n}-\beta_J|_{X_{I\cap J}^\wedge/S_n}+df.
	\end{equation}
	Since $(\Omega_I^1)^\wedge, (\Omega_J^1)^\wedge$ and $A_{I\cap J}$ are all limits of rational $G\times \C^*$-modules, we can assume that $\beta_I, \beta_J$ and $f$ are $G$-invariant and have degree 2. 
	Let $\delta_I$ and $\delta_J$ be derivations of $\C[X_I^\wedge/S_n]$ and $\C[X_J^\wedge/S_n]$ respectively, defined by 
	\begin{equation}
		\iota_{\delta_I} {\omega_I} = \beta_I, 	\iota_{\delta_J} {\omega_J} = \beta_J. 
	\end{equation}
	Since $\beta_I,\beta_J$ are closed, $\delta_I$ and $\delta_J$ are Poisson derivations. 
	
	The $S_n$-action on $X_I^\wedge,X_J^\wedge$ and on $X^\wedge_{I\cap J}$ are free. Therefore, the derivations $\delta,\delta_I$ and $\delta_J$ lift to derivations of $ \C[X_{I}^\wedge]$, $\C[X_{ J}^\wedge]$ and $\C[X_{I\cap J}^\wedge]$ respectively. 
	Let us still denote the lifts by $\delta, \delta_I,\delta_J$. Moreover, we view $f$ as an element of $\C[X_{IJ}^\wedge]$. 
	We obtain the following equality of derivations of $\C[X_{I\cap J}^\wedge]$: 
	\begin{equation}\label{derivation version of coboundary condition}
		\delta = \delta_I - \delta_J + \{f,-\}.
	\end{equation}
	Recall $\delta$ is $G$-, $S_n$- and model $\C^*$-equivariant. Therefore, take the suitable invariant parts of $\delta_I,\delta_J$, we may assume $\delta_I,\delta_J$ are $G$-, $S_n$- and model $\C^*$-equivariant as well. 
	
	We would like to take $\exp$ of both sides of \eqref{derivation version of coboundary condition}, and take $\Psi_I=\exp\delta_I,\Psi_J = \exp\delta_J$ in \Cref{proposition trasition is inner}, but it is not clear that $\delta_I,\delta_J$ are pro-nilpotent; therefore, $\exp \delta_I$ may not be well defined.
	The remainder of this subsection is devoted to modifying $\delta_I,\delta_J$ so that they become pro-nilpotent and Poisson. 
	First, we define a multigrading on $\C[X_{I\cap J}]$ so that $\delta$ ``increases" the multigrading. 
	
	\begin{definition}
		We give the algebra $\C[X_{I\cap J}]$ an algebra multigrading by the lattice 
		\begin{equation}
			\Lambda=  \ZZ \oplus \ZZ^{2\ell}  
		\end{equation} as follows. For the generators in \Cref{lemma CXI generated by linear functions}, we let 
		\begin{itemize}
			\item $\C[G]^H$ and $\C[\tilde{\L}_{I\cap J}]$ have degree 0; 
			\item $(\C[G]\otimes U^*)^H$ has degree $2 \oplus 0  $;  
			\item $(\C[G] \otimes (V^\mathbf{p})^*)^H$ have degree $0  \oplus e_1 $, where $e_i \in \ZZ^{2\ell}$ is 1 at the $i$-th component and 0 elsewhere;   
			\item for $1\le k \le 2\ell -2$, $(\C[G]\otimes(V^k)^*)^H$ has degree $ 0 \oplus e_{k+1}  $; 
			\item  $(\C[G] \otimes (V^\mathbf{q})^*)^H$ has degree $0  \oplus e_{2\ell}$. 
		\end{itemize}
	\end{definition}
	
	The following properties of this multigrading are straightforward to check. 
	\begin{lemma}\label{assorted facts about multigrading}
		The following statements are true. 
		\begin{enumerate}
			\item Each multigraded piece $\C[X_{I\cap J}]_d$ is stable under the $G$-, $S_n$- and model $\C^*$-action.
			\item If $\{\C[X_{I\cap J}]_{0\oplus e_i},\C[X_{I\cap J}]_{0\oplus e_j} \} \neq 0$, then $i+j = 2\ell +1$. 
			\item If a multigraded piece $\C[X_{I\cap J}]_d$ is nonzero, then all components of $d$ are nonnegative. 
			\item For $d\in \Lambda$, let $|d|$ denote the sum of components of $d$. For any $n\ge 0$, 
			\begin{equation}
				\bigoplus_{|d| \ge 2n} \C[X_{I\cap J}]_d \subset I_{I\cap J} ^{n}
			\end{equation}
			\item For any $d_1,d_2\in \Lambda$, we have
			\begin{equation}
				\{\C[X_{I\cap J}]_{d_1},\C[X_{I\cap J}]_{d_2} \} \subset \bigoplus_{|d| = |d_1|+|d_2|-2} \C[X_{I\cap J}]_d
			\end{equation}
		\end{enumerate}
	\end{lemma}
	The following definition is motivated by \Cref{pronilpotency for delta}. 
	\begin{definition}
		Let $\nil$ be the subset of $\Lambda$ defined by 
		\begin{equation}
			\begin{split}
				\nil = &  \{ 0\oplus (-e_k +e_{k'} )  | 1\le k<k'\le 2\ell  \} \\
				\cup &\{-2 \oplus (e_k + e_k') | k+k'> 2\ell +1 \} .
			\end{split}
		\end{equation}
		Note that $|d|=0$ for all $d\in\nil$. 
	\end{definition}
	\begin{lemma}\label{rewrite delta is a pronilpotent derivation}
		We have
		\begin{equation}
			\delta = \sum_{d\in \nil} \delta_d + \sum_{|d'| \ge 1} \delta_{d'}
		\end{equation} 
		where $\delta_d$ is the degree $d$ component of $\delta$. Note the second sum is convergent since $\C[X_{I\cap J}^\wedge]$ is $\idealI_{I\cap J}$-adic complete. 
	\end{lemma}
	\begin{proof}
		We check the lemma on the generators of $\C[X_{I\cap J}^\wedge]$, i.e. we show 
		\begin{equation}
			\delta (a)= \sum_{d\in \nil} \delta_d (a)+ \sum_{|d'| \ge 1} \delta_{d'}(a)
		\end{equation} 
		for $a \in \C[G]^H, (\C[G] \otimes V^*)^H, (\C[G]\otimes U^*)^H$ or  $\C[\tilde{\L}_{I\cap J}]$. The cases for $a\in \C[G]^H, (\C[G] \otimes V^*)^H$ or  $\C[\tilde{\L}_{I\cap J}]$ follows easily from \Cref{pronilpotency for delta}. Let us assume $a= f\otimes u\in (\C[G]\otimes U^*)^H$, where $f\in \C[G], u\in U^*$. By \Cref{pronilpotency for delta}, we have
		\begin{equation}
			\begin{split}
				\delta (f\otimes u) &\in \sum_{i\le j} f_{ij} \otimes v_iv_j + \bigoplus_{|d'|\ge 3} \C[X_{I\cap J}^\wedge]_{d'}, \\ & \text{ where } f_{ij}\in \C[G], v_i\in (V^i)^*,v_j\in (V^j)^*. 
			\end{split}
		\end{equation}
		It suffices to show that $f_{ij} = 0$ unless $i+j > 2\ell+1$. Fix $i$ and let $f'\otimes v_{2\ell +1-i}$ be an arbitrary element of $(\C[G]\otimes (V^{2\ell +1-i})^*)^H$. Since $\delta$ is Poisson, 
		\begin{equation}
			\delta \{f\otimes u, f'\otimes v_{2\ell +1-i}\} = \{\delta f\otimes u, f'\otimes v_{2\ell +1-i}\} + \{f\otimes u, \delta f'\otimes v_{2\ell +1-i}\}.
		\end{equation}
		Let us compare the degree $0 \oplus e_j$ parts of both sides. 
		Since $\{f\otimes u, f'\otimes v_{2\ell +1-i}\}$ has degree $0 \oplus e_{2\ell +1-i}$, 
		we see $(\delta \{f\otimes u, f'\otimes v_{2\ell +1-i}\})_{0\oplus e_j} = 0$ unless $i+j>2\ell+1$. 
		Similarly, $(\{ f\otimes u, \delta f'\otimes v_{2\ell +1-i}\})_{0\oplus e_j} = 0$ unless $i+j>2\ell+1$. 
		The $0\oplus e_j$ component of the first summand of the right hand side is 
		\begin{equation}
			\begin{split}
				\{\delta f\otimes u, f'\otimes v_{2\ell +1-i}\}_{0\oplus e_j} & = \left( \left\{\sum_{i\le j} f_{ij} \otimes v_iv_j , f'\otimes v_{2\ell +1-i} \right\}\right)_{0\oplus e_j} \\
				&= f_{ij} f' v_j \{v_i,v_{2\ell +1-i}\}
			\end{split}
		\end{equation}
		which must also be 0 for any choice of $f'\otimes v_{2\ell +1-i}$ unless $i+j >2\ell+1$. 
		Since the Poisson bracket is nondegenerate on $V$, it follows that $f_{ij} = 0$ when $i+j \le 2\ell+1$. This finishes the proof of the lemma.  
	\end{proof}
	\begin{remarklabeled}
		\Cref{pronilpotency for delta} is slightly stronger than \Cref{rewrite delta is a pronilpotent derivation}. For example, by \Cref{pronilpotency for delta} we see $\delta_{0 \oplus (-e_1+e_2)} = 0$, so we could have excluded $0 \oplus (-e_1+e_2)$ from $\nil$. However, the current definition of $\nil$ is more symmetric and easier to use. 
	\end{remarklabeled}
	Now let us consider $\delta_I$ and $\delta_J$. Define 
	\begin{equation}\label{definition of deltaI0}
		\delta_{I,0} := \sum_{d\in \nil} \delta_{I, d} + \sum_{|d'| \ge 1} \delta_{I,d'} 
	\end{equation}
	and define $\delta_{J,0}$ similarly. 
	\begin{lemma}\label{lemma truncated derivations are poisson}
		The operators $\delta_{I,0}$ and $\delta_{J,0}$ are $G$-, $S_n$- and model $\C^*$-equivariant Poisson derivations of $\C[X_{I\cap J}^\wedge]$. 
	\end{lemma}
	\begin{proof}
		We prove the statement for $\delta_{I,0}$; the proof for $\delta_{J,0}$ is completely analogous. 
		We claim that both $\sum_{d\in \nil} \delta_{I, d} $ and $ \sum_{|d'| \ge 1} \delta_{I,d'} $ are Poisson derivations. 
		
		First, $\sum_{d\in \nil} \delta_{I, d}$ is a finite sum of graded pieces of a derivation on a multigraded algebra, hence a derivation. Since each multigraded piece of $\C[X_{I\cap J}^\wedge]$ is stable under the $G$-, $S_n$- and model $\C^*$-action, $\sum_{d\in \nil} \delta_{I, d}$ is $G$-, $S_n$- and model $\C^*$-equivariant. 
		
		We check $\sum_{d\in \nil} \delta_{I, d}$ is Poisson on the generators in \Cref{lemma CXI generated by linear functions}. Namely, let $a$ and $b$ be elements of $\C[G]^H, (\C[G] \otimes V^*)^H, (\C[G]\otimes U^*)^H$ or  $\C[\tilde{\L}_{I\cap J}]$. We need to show 
		\begin{equation}\label{goal check truncated poisson der}
			\sum_{d\in \nil} \delta_{I, d} \left\{a,b\right\} = \left\{\sum_{d\in \nil} \delta_{I, d} a, b\right\} + \left\{a,\sum_{d\in \nil} \delta_{I, d} b \right\}.  
		\end{equation}
		We check case by case. First, we observe that, for any $d\in \nil$,
		\begin{equation}\label{nil derivation kills terms without V}
			\delta_{I, d} a \in 
			\begin{cases}
				\{0\} &\text{ if $a\in \C[G]^H$ or  $\C[\tilde{\L}_{I\cap J}]$ } \\
				(\C[G] \otimes V^*)^H &\text{ if $a\in(\C[G] \otimes V^*)^H$}\\
				(\C[G] \otimes S^2(V^*))^H &\text{ if $a\in(\C[G] \otimes U^*)^H$}
			\end{cases}.
		\end{equation}
		\begin{enumerate}[labelsep=0cm,align=left,leftmargin=0cm, label=\textit{Case} \arabic*. \  ,itemindent=0.5cm]
			\item Suppose both $a,b$ lie in $\C[G]^H $ or  $\C[\tilde{\L}_{I\cap J}]$. By \eqref{nil derivation kills terms without V}, the right hand side of \eqref{goal check truncated poisson der} is $0$. Since $\{a,b\}$ is homogeneous and its degree lies in $\ZZ \oplus 0$,the degree of $\delta_{I, d}\{a,b\}$ must have a negative component. It follows that the left hand side of \eqref{goal check truncated poisson der} is also $0$. 
			
			\item Suppose $a\in \C[G]^H$ or $\C[\tilde{\L}_{I\cap J}]$ while $b\in (\C[G] \otimes V^*)^H$. 
			In this case, $\{a,b\}=0$ so the left hand side of \eqref{goal check truncated poisson der} is 0. By \eqref{nil derivation kills terms without V} both summands of the right hand side of \eqref{goal check truncated poisson der} are 0 as well. 
			\item Suppose $a\in \C[G]^H$ or $\C[\tilde{\L}_{I\cap J}]$ while $b\in (\C[G] \otimes U^*)^H$. In this case $\{a,b\}$ has degree 0 so the left hand side of \eqref{goal check truncated poisson der} is 0. The first summand of the right hand side is 0, and the second summand is also 0 since $\{\C[G]^H, (\C[G] \otimes S^2(V^*))^H\} = 0$ and 
			$\{\C[\tilde{\L}_{I\cap J}], (\C[G] \otimes S^2(V^*))^H\} = 0$. 
			\item Suppose $a\in (\C[G] \otimes U^*)^H$, $b\in (\C[G] \otimes V^*)^H$; assume $b$ has degree $0\oplus e_k$, $1\le k \le 2\ell$. 
			Then $\{a,b\}$ has degree $0\oplus e_k$. 
			Since $\delta_I$ is Poisson we have 
			\begin{equation}\label{equation for poisson derivation deltaI on a and b}
				\delta_I \{a,b\} = \{a,\delta_I b\} + \{\delta_I a,b\}
			\end{equation}
			and the left hand side of \eqref{goal check truncated poisson der} takes the sum of degree $0\oplus e_{k'}$ components of \eqref{equation for poisson derivation deltaI on a and b} over $k'>k$. 
			
			The sum of degree $0\oplus e_{k'}$ ($k'>k$) parts of $\{ \delta_{I } a, b\}$ is 
			\begin{equation}
				\sum_{k' >k} \{(\delta_I a)_{0\oplus e_{2\ell+1-k } + e_{k'}}, b\}
			\end{equation}
			which, by the definition of $\nil$, is precisely 
			$\{\sum_{d\in \nil} \delta_{I,d} a ,b\}$, since $(2\ell+1-k)+k'>2\ell+1$. 
			
			The sum of degree $0\oplus e_{k'}$ ($k'>k$) parts of $\{  a, \delta_{I }b\}$ is  
			\begin{equation}
				\sum_{k'>k} \{a,(\delta_Ib )_{0\oplus e_{k'}}\}
			\end{equation}
			which is precisely $\left\{a,\sum_{d\in \nil} \delta_{I, d} b \right\}$. This proves \eqref{goal check truncated poisson der}. 
			\item Suppose $a, b\in (\C[G] \otimes V^*)^H$; assume $a$ has degree $0\oplus e_{k_1}$ and $b$ has degree $0 \oplus e_{k_2}$. In this case, the left hand side of \eqref{goal check truncated poisson der} is 0. 
			
			Suppose first $k_1+k_2 \ge 2\ell +1$. Then for any $d \in \nil$, $\delta_{I, d}a $ has degree $0\oplus e_{k'}$ for some $k'>k_1\ge 2\ell +1-k_2$, so $\{\delta_{I, d} a, b\} = 0$; similarly, $\{ a, \delta_{I, d} b\} =0$, so \eqref{goal check truncated poisson der} holds. 
			
			Suppose $k_1+k_2 < 2\ell +1$. Then $\{a,b\}=0$. By \eqref{equation for poisson derivation deltaI on a and b}, $\{\delta_Ia,b\} + \{a,\delta b\} = 0$. In particular, the degree 0 part of the right hand side is 0. Note that
			\begin{equation}
				\{\delta_Ia,b\}_0 = \{(\delta_Ia)_{0\oplus e_{2\ell +1-k_2}},b\} =  \left\{\sum_{d\in \nil} \delta_{I, d} a, b\right\} 
			\end{equation}
			since $0\oplus (-e_{k_1} + e_{2\ell +1-k_2}) \in \nil$. Similarly, 
			\begin{equation}
				\{a,\delta_Ib\}_0 = \{a,(\delta_Ib)_{0\oplus e_{2\ell +1-k_1}}\} =  \left\{ a, \sum_{d\in \nil} \delta_{I, d}b\right\} = 0.
			\end{equation}
			Therefore, the right hand side of \eqref{goal check truncated poisson der} is 0 as well. 
			\item Suppose $a\in (\C[G] \otimes U^*)^H$, $b\in (\C[G] \otimes U^*)^H$. In this case, the two sides of \eqref{goal check truncated poisson der} are the sum of degree $0\oplus (e_k + e_k')$ (with $k+k' > 2\ell +1$) parts of the two sides of \eqref{equation for poisson derivation deltaI on a and b}. 
		\end{enumerate}
		
		This finishes the proof that $\sum_{d\in \nil} \delta_{I, d}$ is a Poisson derivation. 
		
		Next, note that $\sum_{|d'| \ge 1} \delta_{I,d'}$ is a convergent sum of graded pieces of $\delta$, hence it is a $G$-, $S_n$- and model $\C^*$-equivariant derivation. Since $\delta$ is Poisson and the Poisson bracket has degree $-2$ in the sense of (5) of \Cref{assorted facts about multigrading}, it is easy to check that $\sum_{|d'| \ge 1} \delta_{I,d'}$ is Poisson as well. This finishes the proof of \Cref{lemma truncated derivations are poisson}.
	\end{proof}
	
	\begin{lemma}\label{final pronilpotency lemma for the coboundary equation}
		There are Poisson derivations $\delta_{I,0}, \delta_{J,0}$ of $\C[X_I^\wedge], \C[X_J^\wedge]$ respectively, and a function $f_0 \in \C[X_{I\cap J}^\wedge]$, satisfying the following conditions. 
		\begin{enumerate}
			\item $\delta_{I,0}, \delta_{J,0}$ are $G$-, $S_n$- and model $\C^*$-equivariant Poisson derivations, and for any $d\in \Lambda$, 
			\begin{equation}
				\begin{split}
					&\delta_{I,0}\C[X_{I\cap J}^\wedge]_d \subset \bigoplus_{d_1 \in \nil } \C[X_{I\cap J}^\wedge]_{d+ d_1} \oplus \bigoplus_{|d_2|\ge 1} \C[X_{I\cap J}^\wedge]_{d+ d_2} \\
					&\delta_{J,0}\C[X_{I\cap J}^\wedge]_d \subset \bigoplus_{d_1 \in \nil } \C[X_{I\cap J}^\wedge]_{d+ d_1} \oplus \bigoplus_{|d_2|\ge 1} \C[X_{I\cap J}^\wedge]_{d+ d_2} 
				\end{split}
			\end{equation}
			\item $f_0$ is $G$- and $S_n$-invariant, and has degree 2 with respect to the model $\C^*$-action. For any $d\in \Lambda$, 
			\begin{equation}
				\{f_0,\C[X_{I\cap J}^\wedge]_d \}\subset \bigoplus_{d_1 \in \nil } \C[X_{I\cap J}^\wedge]_{d+ d_1} \oplus \bigoplus_{|d_2|\ge 1} \C[X_{I\cap J}^\wedge]_{d+ d_2}. 
			\end{equation}
			\item $\delta = \delta_{I,0} - \delta_{J,0} + \{f_0,-\}$.
		\end{enumerate}
	\end{lemma}
	\begin{proof}
		Let $\delta_I,\delta_J$ and $f$ be as in \eqref{derivation version of coboundary condition}. Let $\delta_{I,0}$ and $\delta_{J,0}$ be as in \eqref{definition of deltaI0}, so that they satisfy (1) by \Cref{lemma truncated derivations are poisson}. 
		
		Write $f = \sum_{d'\in \Lambda} f_{d'}$ where $f_{d'} \in \C[X_{I\cap J}^\wedge]_{d'}$, which is a convergent sum. 
		If $f_{d'} \neq 0$, then $\{f_{d'},\C[X_{I\cap J}^\wedge]_d \}\subset \bigoplus_{d_1 \in \nil } \C[X_{I\cap J}^\wedge]_{d+ d_1} \oplus \bigoplus_{|d_2|\ge 1} \C[X_{I\cap J}^\wedge]_{d+ d_2}$ if and only if $|d'|>2$ or $d' = 0\oplus (e_k+e_{k'})$ with $k+k' > 2\ell +1$. 
		Take 
		\begin{equation}
			f_0:= \sum_{k+k'>2\ell +1}f_{0\oplus (e_k+e_{k'})} + \sum_{|d'|>2} f_{d'}.
		\end{equation}
		Each summand is $G$-, $S_n$-invariant and has degree 2 with respect to the model $\C^*$-action, so that $f_0$ satisfies (2). 
		
		Thanks to \Cref{rewrite delta is a pronilpotent derivation} and \eqref{derivation version of coboundary condition}, (3) is true by the construction of $\delta_{I,0},\delta_{J,0}$ and $f_0$. 
	\end{proof}
	
	\begin{proof}[Proof of \Cref{proposition trasition is inner}]
		Let $\delta_{I,0},\delta_{J,0}$ and $f_0$ be as in \Cref{final pronilpotency lemma for the coboundary equation}. 
		By the definition of $\nil$ and (4) of \Cref{assorted facts about multigrading}, we see $\Psi_I := \exp (-\delta_{I,0}) $ and $\Psi_J := \exp(-\delta_{J,0})$ are well defined symplectic automorphisms on $\C[X_I^\wedge], \C^[X_J^\wedge]$ respectively. They are $S_n$-equivariant, and therefore descends to symplectic automorphisms on $\C[X_I^\wedge/S_n], \C^[X_J^\wedge/S_n]$ respectively. 
		Since $\delta_{I,0}$ and $\delta_{J,0}$ stabilizes the ideals $\idealI_I$ and $\idealI_J$ respectively, $\Psi_I$ and $\Psi_J$ satisfy (1) of \Cref{proposition trasition is inner}. Since $X_{I\cap J} \subset X_I$ is given by inverting a function on $\tilde{\L_I}$, $\Psi_I$ and $\Psi_J$ satisfy (2) of \Cref{proposition trasition is inner} by (1). 
		
		For linear combinations of commutators among $\delta_{I,0}, \delta_{J,0}$ and $\{f_{0},-\}$, the conditions in (1) \Cref{final pronilpotency lemma for the coboundary equation} also hold. Therefore, by the BCH formula, there exists an $f' \in \C[X_{I\cap J}]$ such that 
		\begin{equation}
			(\Psi_J\inverse\theta_{IJ}\Psi_I )^*= \exp\{f',\bullet\}.
		\end{equation}
		The function $f'$ is a convergent sum of functions of the form 
		$\delta_{I,0}^{m_1}\delta_{J,0}^{n_1} \cdots \delta_{I,0}^{m_k}\delta_{J,0}^{n_k} f_0$, where $m_1,n_1,\cdots, m_k,n_k \ge 0$. Therefore, $f'$ is $G$-, $S_n$-invariant and has degree 2 with respect to the model $\C^*$-action, and descends to an element of $\C[X_{I\cap J}^\wedge/S_n]$ satisfying (3) of \Cref{proposition trasition is inner}. This finishes the proof of \Cref{proposition trasition is inner}. 
	\end{proof}
	
	Let us record the following form of \Cref{proposition trasition is inner} and some consequences of its proof. 
	\begin{corollary}\label{structure of nbhd before reduction }
		The following are true.
		\begin{enumerate}
			\item There are $G \times \C^*$-equivariant symplectomorphisms
			\begin{equation}\label{definition of classical thetaI and theta_J}
				\begin{split}
					\theta_I:= \phi_I^\wedge\circ  \Phi_I\circ \Psi_I: X_I^\wedge \to (T^*R)_I^\wedge
					\\
					\theta_J:= \phi_J^\wedge \circ\Phi_J\circ \Psi_J:X_J^\wedge \to (T^*R)_J^\wedge
				\end{split}
			\end{equation}
			where we use the model $\C^*$-action on the sources. 
			\item $\theta_I$ and $\theta_J$ induce isomorphisms of varieties $(G\times^H \tilde{\L}_I)/S_n \cong G\tilde{\L}_I$ and $(G\times^H \tilde{\L}_J)/S_n \cong G\tilde{\L}_J$.
			\item The isomorphisms $\theta_I$ and $\theta_J$ restrict to (different) $G$-equivariant symplectomorphisms $\theta_I|_{I\cap J}, \theta_J|_{I\cap J}: X_{I\cap J}^\wedge/S_n\xrightarrow{\sim} (T^*R)_{I\cap J}^\wedge$. 
			\item There is a function $f\in \C[X_{I\cap J}^\wedge/S_n]$, $G$-invariant and homogeneous with degree $2$ with respect to the $\C^*$-action, and there exists an $N>0$ such that for any $i\ge 0$,
			\begin{equation}\label{pronil condition for bracket with f}
				(\{f, -\})^N \idealI_{I\cap J}^i \subset \idealI_{I\cap J}^{i+1}. 
			\end{equation} 
			Moreover, as algebra automorphisms of $\C[X_{I\cap J}^\wedge/S_n]$,
			\begin{equation}
				(\theta_J|_{I\cap J}\inverse \circ  \theta_I|_{I\cap J})^* = \exp\{f,-\}.
			\end{equation}
		\end{enumerate}
	\end{corollary}
	
	We now turn to formal neighbourhoods of $\L_{I \cup J}$ in $\M_0^0$, by taking suitable Hamiltonian reductions.  
	The formal neighbourhoods $X_I^\wedge,X_J^\wedge,(T^*R)_I^\wedge$ and $(T^*R)_J^\wedge$ are all $G$-invariant, and the $G$-action and $S_n$-action on them commute. 
	Therefore, 
	\begin{align}
		\begin{split}
			(X_I^\wedge/S_n)\hamiltonianreduction_0 G &\cong  (\tilde{\L}_I \times (V\hamiltonianreduction_0H)^{\wedge_0})/S_n \\
			&\cong (\tilde{\L}_I \times \underline{\M}^{\wedge_0})/S_n;\\
			(T^*R)_I^\wedge \hamiltonianreduction_0 G &\cong \M^\wedge_{I};\\
		\end{split}
	\end{align}
	and similar isomorphisms hold when $I$ is replaced by $J, I\cap J$ and $I\cup J$. Here $\M_I^\wedge$ is the completion along $\L_I$ of an open subset of $\M$ that contains $\L_I$ as a closed subset. 
	
	Take Hamiltonian reductions $\hamiltonianreduction_0 G$ for $\theta_I$, $\theta_J$ in \Cref{structure of nbhd before reduction }, we get the following result. 
	We abuse notations and write $\idealI_{I\cap J}$ for the ideal of $\L_{I\cap J}$ in $(\C[\tilde{\L}_{I\cap J}]\completeotimes \C[\underline{\M}]^{\wedge_0})^{S_n}$. 
	\begin{corollary}\label{structure of nbhd after reduction}
		The following are true. 
		\begin{enumerate}
			\item There are $\C^*$-equivariant isomorphisms of affine  Poisson schemes
			\begin{equation}\label{definition of reduced isomorphism thetaredI and thetaredJ}
				\begin{split}
					\theta^{\red}_I: (\tilde{\L}_I \times \underline{\M}^{\wedge_0})/S_n \to \M^\wedge_I\\
					\theta^{\red}_J: (\tilde{\L}_J \times \underline{\M}^{\wedge_0})/S_n \to \M^\wedge_J
				\end{split}
			\end{equation}
			\item The isomorphisms $\theta^{\red}_I$ and $\theta^{\red}_J$ restrict to Poisson isomorphisms $\theta^{\red}_I|_{I\cap J}$ and $\theta^{\red}_J|_{I\cap J}: (\tilde{\L}_{I\cap J}\times\underline{\M}^{\wedge _0})/S_n \xrightarrow{\sim } \M^{\wedge}_{{I\cap J}}$.
			\item There is a function $f\in (\C[\tilde{\L}_{I\cap J}]\completeotimes \C[\underline{\M}]^{\wedge_0})^{S_n}$ such that, as algebra automorphisms of $(\C[\tilde{\L}_{I\cap J}]\completeotimes \C[\underline{\M}]^{\wedge_0})^{S_n}$,
			\begin{equation}(\theta^{\red}_J|_{I\cap J}\inverse\circ\theta^{\red}_I|_{I\cap J})^* = \exp(\{f,\bullet\}).\end{equation}
			Moreover, there exists an $N>0$ such that for any $i\ge 0$,
			\begin{equation}
				(\{f, -\})^N \idealI_{I\cap J}^i \subset \idealI_{I\cap J}^{i+1}. 
			\end{equation} 
		\end{enumerate}
	\end{corollary}
	
\section{Minimally-supported HC bimodules}\label{section Minimally-supported HC bimodules}
	In this section, we obtain quantum analogues of the results of \Cref{Formal neighbourhood of minimal leaf}, and reduce the problem of classifying HC bimodules over $\A_\lambda$ to that over the quantum slice algebra. 
	\subsection{Quantizations of affine opens}\label{subsection Quantizations of affine opens}
	Recall that we have constructed a symplecticomorphism in \Cref{the final affine symplectic isomorphism from XI to TstarRI} of formal schemes. In this subsection we upgrade it to a quantum version. 
	
	Recall the model variety $X_I = G*_H((\g/\cartanh)^*\oplus V)\times \tilde{\L}_I$. We have a natural identification 
	\begin{equation}\label{XI as a hamiltonian reduction}
		X_I = (T^*G\times V\times \tilde{\L}_I)\hamiltonianreduction_0 H
	\end{equation}
	where $H$ acts on $T^*G\times V\times \tilde{\L}_I$ diagonally (with right multiplication on $G$ and trivial action on $\tilde{\L}_I$).
	
	We can microlocalize the algebra of differential operators $D(G)$ on $G$ and get a quantization of $T^*G$, which we still denote by $D(G)$. 
	Equip $D(G)$ with a filtration where the degree of vector fields is 2 and the degree of functions is 0. Take the Rees algebra of $D(G)$ with respect to this filtration, and take its $\hbar$-adic completion, we get a graded formal quantization $D_\hbar(G)$ of $T^*G$. 
	
	Let $\WeylalgebraA_\hbar(V)$ be the homogeneous Weyl algebra of the symplectic vector space $V$. 
	Note that $\tilde{\L}_I$ is a principal open subset of a symplectic vector space $\C^{2\ell (n-1)}$. 
	Let $\WeylalgebraA_\hbar(\tilde{\L}_I) $ denote the microlocalization of the homogeneous Weyl algebra $\WeylalgebraA_\hbar(\C^{2\ell (n-1)})$ to $\tilde{\L}_I$. 
	Similarly, let $\WeylalgebraA_\hbar(\tilde{\L}_{I\cup J}) $ denote the restriction of $\WeylalgebraA_\hbar(\C^{2\ell (n-1)})$ to $\tilde{\L}_{I\cup J}$. 
	
	The right hand side of \eqref{XI as a hamiltonian reduction} has a quantization 
	\begin{equation}
		\tilde{\D}_{I,\hbar} = (D_\hbar(G) \otimes \WeylalgebraA_\hbar(V) \otimes \WeylalgebraA_\hbar(\tilde{\L}_I) ) \hamiltonianreduction_0 H.
	\end{equation}
	Define 
	\begin{equation}\label{Sn invariant quantization before completion}
		\D_{I,\hbar} = (\tilde{\D}_{I,\hbar})^{S_n}
	\end{equation}
	It is a graded formal quantization of $X_I/S_n$. Finally, taking the completion of $\D_{I,\hbar}$ along the closed subset $G\tilde{\L}_I$ of $X_I/S_n$, we obtain a graded formal quantization 
	$
	\D_{I,\hbar}^\wedge
	$
	of $X_I^\wedge/S_n$. 
	
	On the other hand, let $\WeylalgebraA_\hbar(T^*R)$ be the completed homogeneous Weyl algebra of $T^*R$. Localize it to the open subset $(T^*R)_I$ and take its completion along the closed subset $G\tilde{\L}_I$. We get a graded formal quantization of $(T^*R)_I^\wedge$. Denote it by $\mathcal{E}_{I,\hbar}^\wedge$.
	
	Recall the isomorphism $\theta_I$ defined in \Cref{structure of nbhd before reduction }. The quantizations $\mathcal{D}_{I,\hbar}^\wedge$ and $\mathcal{E}_{I,\hbar}^\wedge$ are even and graded quantizations of $X_I^\wedge/S_n$ and $(T^*R)_I^\wedge$ respectively, with respect to the involutory antiautomorphism $\hbar\to -\hbar$ and the the model $\C^*$-action. For more details of even graded quantizations, see \cite[Section 2.2]{losev2012isomorphisms}. 
	Therefore, there is a $\C^*$-equivariant $\C[\hbar]$-algebra isomorphism
	$\vartheta_{I,\hbar}: \mathcal{E}_{I,\hbar}^\wedge \to \mathcal{D}_{I,\hbar}^\wedge $
	that fits in the following commutative diagram 
	\begin{equation}\label{commutative diagram for iso quantization lifting phiI}
		\begin{tikzcd}
			\mathcal{E}_{I,\hbar}^\wedge \arrow[r, "\vartheta_{I,\hbar}"] \arrow[d, "/\hbar"]  & \mathcal{D}_{I,\hbar}^\wedge \arrow[d, "/\hbar"] \\	\C[(T^*R)_I^\wedge] \arrow[r, "\theta_I"]& \C[X_{I}^\wedge/S_n]
		\end{tikzcd}
	\end{equation}
	A completely analogous construction gives an algebra isomorphism $\vartheta_{J,\hbar}: \mathcal{E}_{J,\hbar}^\wedge \to \mathcal{D}_{J,\hbar}^\wedge $
	satisfying a commutative diagram similar to \eqref{commutative diagram for iso quantization lifting phiI}. 
	
	Note that the $S_n$- and $G$-actions on $\tilde{\D}_{I,\hbar}$ commute, so 
	\begin{equation}
		{\D}_{I,\hbar}\hamiltonianreduction_\lambda G \cong (\tilde{\D}_{I,\hbar} \hamiltonianreduction_\lambda G )^{S_n} =( \WeylalgebraA_\hbar(\tilde{\L}_I) \otimes \underline{\A}_{\lambda, \hbar} )^{S_n}. 
	\end{equation}
	Similar isomorphisms hold if $I$ is replaced by $J$ or $ I\cap J$. We thus obtain an isomorphism by applying quantum Hamiltonian reduction $\hamiltonianreduction_\lambda G$ to the isomorphism $\vartheta_{I,\hbar}$, i.e. we have an isomorphism of algebras 
	\begin{equation}\label{bad isomorphism of completed quantum algebra over affine piece}
		\vartheta_{I,\hbar}^{\red} : \A_{\lambda,\hbar}^{\wedge \L_I} \xrightarrow{\sim} ( \WeylalgebraA_\hbar(\tilde{\L}_I) \completeotimes \underline{\A}_{\lambda, \hbar}^{\wedge_0})^{S_n}.
	\end{equation}
	Similarly, we have an isomorphism $\vartheta_{J,\hbar}^{\red} : \A_{\lambda,\hbar}^{\wedge \L_J} \cong ( \WeylalgebraA_\hbar(\tilde{\L}_J) \completeotimes \underline{\A}_{\lambda, \hbar}^{\wedge_0})^{S_n}.$
	\subsection{Lifting automorphisms}\label{subsection Lifting automorphisms}
	We can microlocalize the algebras $\mathcal{D}_{I,\hbar}^\wedge, \mathcal{D}_{J,\hbar}^\wedge$ and obtain an algebra
	$\mathcal{D}_{I\cap J,\hbar}^\wedge$ quantizing the affine formal scheme $X_{I\cap J}^\wedge$, 
	and similarly an algebra $\mathcal{E}_{I\cap J,\hbar}^\wedge$ quantizing $(T^*R)_{I\cap J}^\wedge$. 
	Let $\idealI_{I\cap J,\hbar
	}\subset D_{I\cap J,\hbar}^\wedge$ be the preimage of $\idealI_{I\cap J}$ under the projection $D_{I\cap J,\hbar}^\wedge \twoheadrightarrow A_{I\cap J} $. By construction, $D_{I\cap J,\hbar}^\wedge$ is complete in the $\idealI_{I\cap J,\hbar
	}$-adic topology. 
	
	The first main result of this subsection is the following quantum version of \Cref{structure of nbhd before reduction }. 
	
	\begin{proposition}\label{quantum inner automorphism}
		There is an element $F \in D_{I\cap J,\hbar}^\wedge$ such that 
		\begin{itemize}
			\item $F$ is $G$-invariant; 
			\item $F$ is homogeneous of degree $2$ with respect to the model $\C^*$-action, and
			\item $F$ lifts the element $f\in A_{I J}$ in (3) of \Cref{structure of nbhd before reduction }. 
			
		\end{itemize}
		There are $G\times \C^*$-equivariant automorphims $\vartheta'_{I,\hbar}$ (resp. $\vartheta'_{J,\hbar}$) of $\D_{I,\hbar}^\wedge$ (resp. $\D_{J,\hbar}^\wedge$). 
		These objects satisfy:
		\begin{enumerate}
			\item  $\vartheta'_{I,\hbar}$ and $\vartheta'_{J,\hbar}$ are identity modulo $\hbar$, and they restrict to automorphisms of $\D_{I\cap J,\hbar}^\wedge$; 
			\item 
			there is some $N>0$ such that, for any $i\ge 0$, 
			\begin{equation}\label{prounipotency for quantum ad F}
				(\frac{1}{\hbar^2}\ad F)^N \idealI_{I\cap J,\hbar}^i \subset \idealI_{I\cap J,\hbar}^{i+1}.
			\end{equation}
			\item 
			when restricted to $\D_{I\cap J,\hbar}^\wedge$,
			\begin{equation} \label{inner transition automorphism before reduction}
				\vartheta'_{I,\hbar} 
				\circ \vartheta_{I,\hbar}  
				\circ \vartheta_{J,\hbar}\inverse \circ   (\vartheta'_{J,\hbar})\inverse 
				= \exp \left( \frac{1}{\hbar^2}\ad F\right).
			\end{equation}
			
		\end{enumerate}
		
	\end{proposition}
	\begin{proof}
		To simplify notation let us write $\D_\hbar$ for $\D_{I\cap J,\hbar}^\wedge$, and $A$ for $A_{I\cap J}$ in this proof. 
		
		\textit{Step 1.}
		We first show that for any lift $\tilde{f}\in \D_\hbar$ of $f$, the operator $\exp ( \frac{1}{\hbar^2}\ad \tilde{f})$ converges. 
		Denote $\partial := \frac{1}{\hbar^2}\ad \tilde{f}$. 
		It is a lift of the derivation $\{f,\bullet\}$ of $A$. 
		Note that $\partial \idealI_{I\cap J, \hbar} + (\hbar) \subset \idealI\subset A$, and $\hbar \in \idealI_{I\cap J, \hbar}$. 
		We conculde $\partial\idealI_{I\cap J, \hbar} \subset \idealI_{I\cap J, \hbar}$. 
		
		Therefore, it suffices to show that, for any $n\ge 0$, there is a large enough $N$ such that 
		\begin{equation}\label{temp partial N is pronilpotent}
			\partial^N \D_\hbar\subset \idealI_{I\cap J, \hbar}^n.
		\end{equation}
		By \eqref{pronil condition for bracket with f}, for any $n\ge0$, there is some $N>0$ such that 
		$\{f,\bullet\}^N A_{I\cap J}\subset \idealI_{I\cap J}^{n}.$
		This implies 
		\begin{equation}
			\partial^N \D_\hbar^\wedge \subset \idealI_{I\cap J, \hbar}^n + (\hbar).
		\end{equation}
		Apply $\partial^N$ to the above inclusion, we see 
		\begin{equation}
			\partial^{2N} \D_\hbar \subset \idealI_{I\cap J, \hbar}^n + \hbar(\idealI_{I\cap J, \hbar}^n + (\hbar))  \subset \idealI_{I\cap J, \hbar}^n + (\hbar^2).
		\end{equation}
		Repeating this argument, we see
		\begin{equation}\partial^{nN} \D_\hbar \subset (\idealI_{I\cap J, \hbar})^n + (\hbar^n).\end{equation}
		But since $\hbar^n\subset \idealI_{I\cap J, \hbar}^n$, we see $\partial^{nN} \D_\hbar \subset (\idealI_{I\cap J, \hbar})^n$, as desired. 
		
		\textit{Step 2.} Denote $\Theta = \vartheta_{I,\hbar}  
		\circ \vartheta_{J,\hbar}\inverse$, viewed as an algebra automorphism of $\D_\hbar$. We show $\log\Theta$ is a well defined derivation of $\D_\hbar$. 			
		Let $\tilde{f}$ be any lift of $f$, so that 
		\begin{equation}\label{big Theta = exp adftilde mod hbar}
			\Theta \equiv \exp \frac{1}{\hbar^2} \ad \tilde{f} \mod \hbar. 
		\end{equation}
		As in Step 1, let $\partial =  \frac{1}{\hbar^2}\ad \tilde{f}$. Then $\Theta - 1 = \partial(1 + \frac{1}{2}\partial^2+ \cdots) + \hbar\psi$ for some endomorphism $\psi$ of $\D_\hbar$. Therefore, thanks to \eqref{temp partial N is pronilpotent}, for each $n\ge 0$ there is some $N$ large enough such that
		\begin{equation}
			(\Theta-1)^N \D_\hbar \subset \idealI_{I\cap J, \hbar}^n.
		\end{equation} 
		It follows that $\log \Theta$ converges. 
		
		\textit{Step 3.} We construct the automorphisms $\vartheta'_{I,\hbar}$ and $\vartheta'_{J,\hbar}$, and prove part (2). 
		
		Let $\tilde{f}$ be a degree 2 $G$-invariant lift of $f$. 
		Denote $\Delta = \log\Theta, \delta = \Delta - \partial$, so that $\delta $ is a $G \times \C^*$-equivariant derivation on $\D_\hbar$, and by \eqref{big Theta = exp adftilde mod hbar}, $\delta \equiv 0 \mod \hbar$. The derivation $\frac{1}{\hbar} \delta$ descends to a derivation on $A$. 
		Since the Poisson bracket on $A$ descents from the commutator on $\D_\hbar$, and $\Theta, \exp(\frac{1}{\hbar^2} \ad \tilde{f}) $ are both algebra automorphisms, it is clear that $\frac{1}{\hbar}\delta$ preserves the Poisson bracket. 
		In other words, the vector field on $X_{I\cap J}$ corresponding to $\frac{1}{\hbar}\delta$ is symplectic. 
		The same arguments as in \Cref{inner derivation} shows there are $G$-invariant, degree 1 closed 1-forms 
		\begin{equation}\beta_I \in (\Omega_I^1)^\wedge,
			\beta_J \in (\Omega_J^1)^\wedge,
		\end{equation}
		and a $G$-invariant, degree 1 function $f_0\in A$ such that
		\begin{equation}\label{betafor delta over hbar square is a coboundary}
			\iota_{\frac{1}{\hbar}\delta} \omega= \beta_I|_{X_{I\cap J}^\wedge/S_n}-\beta_J|_{X_{I\cap J}^\wedge/S_n}+df_0.
		\end{equation}
		Define derivations $\delta_I$ on $A_I$ by $\iota_{\delta_I}\omega = \beta_I$, and similarly, $\delta_J$ on $A_J$ by $\iota_{\delta_J}\omega = \beta_J$. They are Poisson derivations since $\beta_I,\beta_J$ are closed forms. 		
		We remark that we do not require $\delta_I,\delta_J$ and $f$ to satisfy any nilpotency condition. 
		
		Let $\tilde{\delta}_I, \tilde{\delta_J} $ be derivations on $\D_\hbar$ lifting  $\delta_I,\delta_J$ respectively; such lifts exist thanks to \cite[Lemma 2.4]{losev2021harish}. Let $\tilde{f_0}\in \D_\hbar$ be a degree 1 lift of $f_0$. 
		Note that $\Der(\D_\hbar)$ is a pro-rational $G\times \C^*$-module. Hence, by taking invariant components, we may assume $\tilde{\delta}_I, \tilde{\delta_J} $ are $G$-equivariant and degree $-1$.
		We see on $\D_\hbar$,
		\begin{equation}
			\Delta \equiv \hbar \tilde{\delta}_I -\hbar \tilde{\delta}_J + \frac{1}{\hbar^2}\ad(\tilde{f}+ \hbar\tilde{f_0}) \mod \hbar^2.
		\end{equation}
		
		Let $\vartheta'_{I,1} = \exp(\hbar\tilde{\delta}_I)$ and $\vartheta'_{J,1} = \exp(\hbar\tilde{\delta}_J)$.
		By the BCH formula, we see there is a $G$-invariant, degree 2 element $F_1\in \D_\hbar$ satisfying $F_1 \equiv \tilde{f}\mod \hbar^2$, and 
		\begin{equation}
			\Theta \equiv \vartheta'_{I,1}\circ \exp(\frac{1}{\hbar^2} F_1) \circ (\vartheta'_{J,1})\inverse \mod \hbar^2.
		\end{equation}
		Comparing with \eqref{big Theta = exp adftilde mod hbar} ,we have increased the order of $\hbar$ we mod out. 
		Repeat the above argument and use $\hbar$-completeness, we conclude the existence of $\vartheta'_{I,\hbar},\vartheta'_{J,\hbar}$ and $F$ satisfying (1) and (2) of the lemma. 
		The pro-nilpotency condition \eqref{prounipotency for quantum ad F} is a consequence of \eqref{pronil condition for bracket with f} and the fact $\hbar\in \idealI_{I\cap J, \hbar}$. 
	\end{proof}
	
	The $G\times \C^*$-equivariant automorphisms $\vartheta'_{I,\hbar} $ and $ \vartheta'_{J,\hbar} $  descend to automorphisms $(\vartheta'_{I,\hbar})^{\red} $ of $ {\D}_{I,\hbar}\hamiltonianreduction_\lambda G = ( \WeylalgebraA_\hbar(\tilde{\L}_I) \otimes \underline{\A}_{\lambda, \hbar} )^{S_n}$ and  $ {\D}_{J,\hbar}\hamiltonianreduction_\lambda G$ respectively. The $G$-invariant degree 2 element $F$ descents to a degree 2 element of ${\D}_{I\cap J,\hbar}\hamiltonianreduction_\lambda G$, which we still denote by $F$. The identity \eqref{inner transition automorphism before reduction} then descends to 
	\begin{equation}
		\Big{( } (\vartheta'_{I,\hbar})^{\red} \circ \vartheta_{I,\hbar}^{\red}\Big{) }
		\circ \Big{( }(\vartheta'_{J,\hbar})^{\red} \circ \vartheta_{J,\hbar}^{\red}\Big{) }\inverse = \exp \left( \frac{1}{\hbar^2}\ad F\right).
	\end{equation}
	This is an identity of automorphisms of ${\D}_{I\cap J,\hbar}\hamiltonianreduction_\lambda G$. 
	
	In what follows we will mostly use the composition isomorphism  $(\vartheta'_{I,\hbar})^{\red} \circ \vartheta_{I,\hbar}^{\red}$ instead of its two components. Thus, by renaming the isomorphisms, we get the following result. 
	\begin{corollary}\label{key}
		There are $ \C^*$-equivariant algebra isomorphisms
		\begin{equation}\label{isomorphism of completed quantum algebra over affine piece}		
			\begin{split}
				\vartheta_{I,\hbar}^{\red}:
				\A_{\lambda,\hbar}^{\wedge \L_I} \xrightarrow{\sim} ( \WeylalgebraA_\hbar(\tilde{\L}_I) \completeotimes \underline{\A}_{\lambda, \hbar}^{\wedge_0})^{S_n}\\
				\vartheta_{J,\hbar}^{\red}:
				\A_{\lambda,\hbar}^{\wedge \L_J} \xrightarrow{\sim} ( \WeylalgebraA_\hbar(\tilde{\L}_J) \completeotimes \underline{\A}_{\lambda, \hbar}^{\wedge_0})^{S_n}.
			\end{split}		
		\end{equation}
		There is a degree 2 element $F\in {\D}_{I\cap J,\hbar}\hamiltonianreduction_\lambda G$ that lifts the element $f\in A_{I J}$ in (3) of \Cref{structure of nbhd after reduction}.
		They satisfy 
		\begin{enumerate}
			\item $\vartheta_{I,\hbar}^{\red}$ and $\vartheta_{J,\hbar}^{\red}$ descend to the classical isomorphisms $\theta_I$ and $\theta_J$ modulo $\hbar$, respectively. 
			\item there is some $N>0$ such that, for any $i\ge 0$,  
			\begin{equation} 
				(\frac{1}{\hbar^2}\ad F)^N \idealI_{I\cap J,\hbar}^i \subset \idealI_{I\cap J,\hbar}^{i+1}.
			\end{equation}
			\item \begin{equation}
				\vartheta_{I,\hbar}^{\red} \circ  (\vartheta_{J,\hbar}^{\red} )\inverse = \exp \left( \frac{1}{\hbar^2}\ad F\right).
			\end{equation}
		\end{enumerate}
	\end{corollary}

	\subsection{Poisson bimodules}\label{subsection on coherent poisson bimodule over algebra} 
	We need the notion of Poisson bimodules over quantizations of Poisson algebras, following \cite[Section 3.3]{losev2012completions}. 
	
	Let $\A = \bigcup_{i=0}^\infty \A_{\le i}$ be a filtered $\C$-algebra such that $\gr \A$ is a Noetherian Poisson algebra $A$, where the Poisson bracket induced by the commutator on $\A$ has degree $-2$. 
	Let $\A_\hbar$ be the Rees algebra of $\A$, so that
	$[\A_\hbar,\A_\hbar]\subset \hbar^2 \A_\hbar$. 
	We define a $\C[\hbar]$-linear Poisson bracket 
	\begin{equation}\{\bullet,\bullet\} : \A_\hbar\otimes_{\C[\hbar]} \A_\hbar \to \A_\hbar, \{a_1,a_2\}:= \frac{1}{\hbar^2}[a_1,a_2].\end{equation}
	Clearly, the bracket is skew-symmetric,satisfies the Leibniz rule in each component, and satisfies the Jacobi identity.

	\begin{definition}\label{definition of Poisson bimodules over quantization}
		Let $\A_\hbar$ be as above and let $\B_\hbar$ be an $\A_\hbar$-bimodule such that the left and right $\C[\hbar]$ actions coincide. We say $\B_\hbar$ is a \textit{Poisson} bimodule over $\A_\hbar$ if $\B_\hbar$ is equipped with a $\C[\hbar]$-linear bracket map 
		\begin{equation}\{\bullet,\bullet\}: \A_\hbar\otimes \B_\hbar \to \B_\hbar,\end{equation}
		satisfying the following conditions: 
		\begin{align}
			&[a,b] = \hbar^2 \{a,b\} \label{poisson bracket vs commutator} \\
			&	\{a_1,a_2 b \} = \{a_1,a_2\}b + a_2\{a_1,b\} \\
			& \{a_1a_2,b\} = a_1\{a_2,b\} + \{a_1,b\}a_2 \\
			& \{\{a_1,a_2\},b\} = \{a_1,\{a_2,b\}\}- \{a_2,\{a_1,b\}\}
		\end{align}
		for any $a,a_1,a_2\in \A_\hbar$ and $b\in\B_\hbar$. 
	\end{definition}
	The following lemma is standard. 
	\begin{lemma}
		A finitely generated graded Poisson $\A_\hbar$-bimodule is finitely generated as a left $\A_\hbar$-module and as a right $\A_\hbar$-module. 
	\end{lemma}
	
	\begin{proof}
		Let $V\subset \B_\hbar$ be a finite dimensional $\C^*$-stable subspace such that $\B_\hbar = \A_\hbar V \A_\hbar$. By \eqref{poisson bracket vs commutator}, $\B_\hbar = \A_\hbar V + \hbar^2 \B_\hbar$; hence, $\B_\hbar = \A_\hbar V + \hbar^{2n} \B_\hbar$ for all $n\ge 1$. Since $\B_\hbar$ is $\C^*$-locally finite and the weights of elements of $\B_\hbar$ are bounded from below, we conclude $\B_\hbar = \A_\hbar V$. Right generation is proved analogously. 
	\end{proof}
	Let $\B_\hbar$ be a coherent Poisson $\A_\hbar$-bimodule. Then $\B_\hbar/\hbar \B_\hbar$ is a finitely generated $A$-module.
	\begin{definition}\label{definition of support of poisson bimodule}
		The support of $\B_\hbar$ is defined to be the support of $\B_\hbar/\hbar \B_\hbar$ in $\Spec A$, which is a closed Poisson subvariety. 
	\end{definition}
	The relation between $\HC$ and coherent Poisson bimodules is as follows, c.f. \cite[Secion 2.4]{losev2021harish}. 
	Let $\B\in \HC(\A )$. 
	Fix a good filtration of $\B$ and let $\B_\hbar= \bigoplus_{i=0}^\infty \B_{\le i} \hbar^i$ be its Rees bimodule. Then $\B_\hbar$ is a finitely generated Poisson bimodule. Conversely, if $\B_\hbar$ is a finitely generated Poisson $\A_\hbar$ bimodule, then $ \B_\hbar/(\hbar-1)$ is an object of $ \HC(\A)$, together with a good filtration. 
	
	Let us consider completions. Let $\A_\hbar,A$ be as before, and let $I\subset A$ be a graded ideal. Let $\idealI_\hbar$ be the preimage of $I$ under the map $\A_\hbar \to A$, which is $\C^*$-stable. Define 
	\begin{equation}\label{definition of Ihbar adic completion algebra}
		\A_\hbar^\wedge : = \varprojlim_k \A_\hbar / \idealI_\hbar^k.
	\end{equation} 
	The Poisson bracket on $\A_\hbar$ extends to a Poisson bracket on $\A_\hbar^\wedge$ by continuity. 
	We define $\C^*$-equivariant Poisson $\A_\hbar^\wedge$-bimodules $\B_\hbar'$  analogously to \Cref{definition of Poisson bimodules over quantization}. 
	
	For example, let $\B_\hbar$ be a finitely generated graded Poisson $\A_\hbar$-bimodule. 
	Define 
	\begin{equation}\label{definition of Ihbar adic completion module}
		\B_\hbar^\wedge : = \varprojlim_k \B_\hbar / \idealI_\hbar^k\B_\hbar.
	\end{equation} 
	This is a left $\A_\hbar^\wedge $-module. 
	Note that $[\idealI_\hbar, \B_\hbar] \subset \hbar\B_\hbar$, so $\B_\hbar'/ \idealI_\hbar^k \B_\hbar'$ is actually a $\A_\hbar^\wedge$-bimodule. 
	We can extend the Poisson bracket on $\B_\hbar$ to a Poisson bracket $\A_{\hbar}^{\wedge} \times \B_\hbar^{\wedge} \to \B_\hbar^{\wedge}$ by continuity, similar to \cite[Lemma 3.6.1]{losev2012completions}. If $\B_\hbar$ is graded, we can extend it to a $\C^*$-action on $\B_\hbar'$ by continuity as well, and each $\B_\hbar/\idealI_\hbar^k \B_\hbar$ is $\C^*$-locally finite. 
	
	The algebra $\A_\hbar^\wedge$ satisfies the conditions of \cite[Lemma 2.4.2]{losev2011finite}, so if $\B_\hbar'$ is a Poisson $\A_\hbar^\wedge$-bimodule that is finitely generated as a left $\A_\hbar^\wedge$-module. An analogue of \cite[Lemma 2.4.4]{losev2011finite} implies $\B_\hbar'$ is complete and separated in the $\idealI_\hbar$-adic topology. 

	\subsection{Sheaves and microlocalization}\label{subsection sheaves and microlocalization}
	Let us proceed to sheaves. We first recall the notion of microlocalization, c.f. \Cite[Section 2.4]{losev2012completions}. We first microlocalize the algebra $\A_\hbar^\wedge$. 
	Let $\A_\hbar^\wedge$ be as in \eqref{definition of Ihbar adic completion algebra}. 
	Let $U$ be a principal open subset of $\Spec A/I$ defined by an element $f\in A$. Let $\hat{f}$ be a lift of $f$ in $\A_\hbar^\wedge$. Let $S$ denote the set $\{\hat{f}^n|n\ge 0\}$. Thanks to \eqref{poisson bracket vs commutator} and the observation $\hbar\in \idealI_\hbar$, the set $S$ satisfies the Ore condition for $\A_\hbar^\wedge / (\idealI_\hbar)^k $, $k\ge 1$. Therefore, the localization $\A_\hbar^\wedge /\idealI_\hbar^k[\hat{f}\inverse]$ is well defined. Take principal open subsets corresponding to all elements of $A$, we see the algebras $\A_\hbar^\wedge /\idealI_\hbar^k[\hat{f}\inverse]$ glue to a sheaf on $\Spec A/I$. 
	\begin{definition}
		We define the microlocalization of the algebra $ \A_\hbar^\wedge$
		by 
		\begin{equation}
			\Loc(\A_\hbar^\wedge): = \varprojlim_k \Loc(\A_\hbar^\wedge/ \idealI_\hbar^k ).
		\end{equation}
		This is a sheaf of algebras on $\Spec A/I$ in the Zariski topology. 
	\end{definition}
	Note that the global sections of $\Loc(\A_\hbar^\wedge)$ is precisely $\A_\hbar^\wedge$. 
	
	We proceed to microlocalizing modules. 
	Let $\B_\hbar'$ be a Poisson bimodule over $\A_\hbar^\wedge$, $f,\hat{f},S$ be as above. Then for each $k>0$, $S$ satisfies the Ore condition for $\B_\hbar' / (\idealI_\hbar)^k \B_\hbar'$. 
	More precisely, for any $b\in \B_\hbar' / (\idealI_\hbar)^k \B_\hbar'$, we can express $b \hat{f}\inverse$ in terms of linear combinations of elements of the form $\hat{f}^{-n}b_n$, as follows. Write $[\hat{f},b] =\hbar b_1$; it is easy to check that 
	\begin{equation}
		b\hat{f}\inverse = \hat{f}\inverse b + \hbar\hat{f}\inverse (b_1\hat{f}\inverse).
	\end{equation}
	Repeat this process and note $\hbar^k\in \idealI^k$, we see there are elements $b_1,\cdots,b_{k-1} \in  \B_\hbar' / \idealI_\hbar^k \B_\hbar'$ such that 
	\begin{equation}
		b\hat{f}\inverse = \hat{f}\inverse b + \sum_{i=1}^{k-1}\hbar^i \hat{f}^{-i-1}b_i.
	\end{equation}
	Therefore, the localization $\B_\hbar' / (\idealI_\hbar)^k \B_\hbar' [\hat{f}\inverse]$ is well defined as a bimodule over $\A_\hbar^\wedge /\idealI_\hbar^k[\hat{f}\inverse]$. 
	They glue to a sheaf $\Loc(\B_\hbar'/ \idealI_\hbar^k \B_\hbar')$ on $\Spec A/I$ in the Zariski topology. 
	\begin{definition}
		We define the microlocalization of $\B_\hbar'$ by 
		\begin{equation}
			\Loc(\B_\hbar'): = \varprojlim_k \Loc(\B_\hbar'/ \idealI_\hbar^k \B_\hbar').
		\end{equation}
		This is a sheaf of $\Loc(\A_\hbar^\wedge)$-bimodules on $\Spec A/I$ in the Zariski topology. If $Y$ is the closed subvariety defined by $I$, we also write $\Loc_Y(\B_\hbar')$ for the above microlocalization. 
	\end{definition} 
	
	Assume the Poisson variety $\Spec A$ has finitely many symplectic leaves and let $\L$ be one of them. Let $I \subset A$ denote the ideal corresponding to $\overline{\L}$, and define $\A_{\hbar}^{\wedge} $ as in \eqref{definition of Ihbar adic completion algebra}. 
	Denote by $\A_{\hbar}^{\wedge\L} $ the restriction of $\Loc (\A_\hbar^\wedge)$ to $\L\subset \overline{\L}$, and by $\idealI_\hbar^{\wedge\L}$ the restriction of the ideal sheaf $\Loc(\idealI_\hbar^\wedge)$ to $\L$.
	
	Let $\B_\hbar$ be a finitely generated Poisson $\A_\hbar$-bimodule, and let $\B_\hbar^\wedge$ be as before. 
	Denote by $\B_\hbar^{\wedge \L} $ the restriction of $\Loc(\B_\hbar^\wedge)$ to $\L \subset\overline{\L}$. 
	This is a sheaf of Poisson $\A_{\hbar}^{\wedge\L} $-bimodules over $\L$ that are finitely generated as left modules, complete and separated in the $\idealI_\hbar^{\wedge \L}$-adic topology. 
	
	In particular, in our situation of quantized flower quiver varieties (and quantized slice quiver varieties) we can take $\A_\hbar$ to be $\A_{\lambda,\hbar}$ or $\underline{\A} _{\lambda,
		\hbar}$ in \Cref{introduction section quantization}, 
	take $\L$ to be the minimal relevant leaf of $\M_0^0$ to construct 
	$\A_{\lambda,\hbar}^ {\wedge \L}$, and take the conical point $0$ to construct 
	$\underline{\A}_{\lambda,\hbar} ^ {\wedge 0}$; the latter is a sheaf over a point and is considered as an algebra. 
	We can further restrict $\A_{\lambda,\hbar}^ {\wedge \L }$ to open subsets $\L_{I}, \L_{I\cup J}$ etc, and get sheaves of Poisson algebras $\A_{\lambda,\hbar}^ {\wedge \L_{I}}$, $\A_{\lambda,\hbar}^ {\wedge \L_{I\cup J}}$, etc. 
	
	\begin{definition}\label{definition of sheaf of coherent poisson bimodules}
		Let $A,\A_\hbar, \L$ be as above, and let $\L_0$ be an open subset of $\L$ (we allow $\L_0 = \L$). 
		We assume $\B_\hbar'$ is a sheaf of $\C^*$-equivariant Poisson $\A_\hbar^{\wedge \L_0}$-bimodules that is complete and separated in the $\idealI_{\L_0,\hbar}$-adic topology. We say $\B_\hbar'$ is
		\begin{enumerate}
			\item \textit{coherent}, if there is an affine open cover $\{U_\alpha\}$ of $\L_0$, such that for all $\alpha$, $H^0(U_\alpha,\B_\hbar')$ is finitely generated as a left module over $H^0(U_\alpha,\A_\hbar^{\wedge \L_0})$, and 
			\begin{equation}
				\B_\hbar|_{U_\alpha}  = \Loc(H^0(U_\alpha,\B_\hbar')).
			\end{equation}
			\item \textit{$\C^*$-prorational}, if the $\C^*$-action on each $\B_\hbar'/\idealI_{\L_0,\hbar}^k\B_\hbar'$ is locally finite.
		\end{enumerate} 
		
	\end{definition}
	
	\subsection{Various categories}\label{subsection various categories}
	In this subsection we return to the setup of quantized flower quiver varieties, and define several categories of interest. As is mentioned in \Cref{sbsbsection introduction Bimodules over quantized algebras}, our goal is to produce functors between these categories that lead to a bijection of isomorphism classes of simple Harish-Chandra bimodules. 
	
	For simplicity of notations, we write
	\begin{equation}\label{notation of Dredhbar}
		\D^{\red}_{\bullet,\hbar}: =   (\WeylalgebraA_\hbar(\tilde{\L}_{\bullet }) \completeotimes_{\C[\hbar]} \underline{\A}_{\lambda,\hbar} ^ {\wedge 0})^{S_n}
	\end{equation} and 
	\begin{equation}\label{notation of tildeDredhbar}
		\widetilde{\D}^{\red}_{\bullet,\hbar}: =   \WeylalgebraA_\hbar(\tilde{\L}_{\bullet}) \completeotimes_{\C[\hbar]} \underline{\A}_{\lambda,\hbar} ^ {\wedge 0} .
	\end{equation}
	where $\bullet$ can be $I,J,I\cap J, I\cup J$. 
	Let $\idealm\subset \underline{\A}_{\lambda,\hbar} ^ {\wedge 0}$ be the maximal ideal. 
	By an abuse of notation, write $\idealI_{\bullet,\hbar}$ for the ideal $(\WeylalgebraA_\hbar(\tilde{\L}_{\bullet }) \completeotimes_{\C[\hbar]}\idealm)^{S_n} \subset \D^{\red}_{\bullet,\hbar}$ 
	and also for the ideal 
	$\WeylalgebraA_\hbar(\tilde{\L}_{\bullet }) \completeotimes_{\C[\hbar]} \idealm\subset \widetilde{\D}^{\red}_{\bullet,\hbar}$. 
	
	In the remainder of this section, we will omit $\lambda$ and write $\A_\hbar$ for $\A_{\lambda,\hbar}$, $\underline{\A}_\hbar$ for $\underline{\A}_{\lambda,\hbar}$. Recall that $\L$ denotes the minimal relevant symplectic leaf. 
	\begin{definition}
		We define the following categories of Harish-Chandra and Poisson bimodules. 
		\begin{enumerate}
			\item The category $\HC(\A)_F$ consists of Harish-Chandra bimodules over $\A$ equipped with a specified good filtration; morphisms are filtered bimodule morphisms. Let $\HC_{\L}(\A)_F$ be its full subcategory consisting of objects supported on $\overline{\L}$. 
			\item The category $\CP(\A_\hbar)$ consists of graded coherent Poisson $\A_\hbar$-bimodules, \Cref{definition of Poisson bimodules over quantization}; $\CP_{\overline{\L}}(\A_\hbar)$ is its full subcategory  consisting of objects supported on $\overline{\L}$. 
			\item The category $\CP(\A_\hbar^{\wedge\overline{\L}})$ consists of $\C^*$-equivariant coherent Poisson $\A_\hbar^{\wedge\overline{\L}}$-bimodules complete and separated in the $\idealI_{\overline{\L},\hbar}$-topology;  $\CP_{\overline{\L}}(\A_\hbar^{\wedge\overline{\L}})$ is its full subcategory consisting of objects supported on $\overline{\L}$. 
			\item The category $\CP(\A_\hbar^{\wedge{\L}})$ consists of sheaves of $\C^*$-equivariant Poisson $\A_\hbar^{\wedge{\L}}$-bimodules complete and separated in the $\idealI_{ \L,\hbar}$-topology that are coherent and $\C^*$-prorational;  $\CP_{\L}(\A_\hbar^{\wedge{\L}})$ is its full subcategory consisting of objects supported on $\L$. 
			\item The category $\CP(\A_\hbar^{\wedge{\L_{\bullet}}})$, where $\bullet$ can be $I,J,I\cup J$, consists of sheaves of $\C^*$-equivariant Poisson $\A_\hbar^{\wedge{\L_{\bullet}}}$-bimodules complete and separated in the $\idealI_{{\bullet},\hbar}$-topology that are coherent and $\C^*$-prorational. 
			The category $\CP_{\L_\bullet}(\A_\hbar^{\wedge{\L_{\bullet}}})$ is its full subcategory consisting of objects supported on $\L_\bullet$.
			\item The category $\CP(\D^{\red}_{\bullet,\hbar})$ consists of sheaves of $\C^*$-equivariant (with respect to the model $\C^*$-action) Poisson $\D^{\red}_{\bullet,\hbar}$-bimodules complete and separated in the $\idealI_{{ \bullet},\hbar}$-topology that are coherent and $\C^*$-prorational;  $\CP_{\bullet}(\D^{\red}_{\bullet,\hbar})$ is its full subcategory consisting of objects supported on $\L_\bullet$.  
			\item 
			The category $\CP (\widetilde{\D}^{\red}_{\bullet,\hbar})$ consists of sheaves of $\C^*$-equivariant (with respect to the model $\C^*$-action) Poisson $\D^{\red}_{\bullet,\hbar}$-bimodules complete and separated in the $\idealI_{{ \bullet},\hbar}$-topology that are coherent and $\C^*$-prorational; $\CP_{\widetilde{\L}_\bullet}(\widetilde{\D}^{\red}_{\bullet,\hbar})$ is its full subcategory consisting of objects supported on $\widetilde{\L}_\bullet$. 
			The symmetric group $S_n$ acts on $\widetilde{\D}^{\red}_{\bullet,\hbar}$. Let  $\CP^{S_n}(\widetilde{\D}^{\red}_{\bullet,\hbar})$,  $\CP^{S_n}_{\widetilde{\L}_\bullet}(\widetilde{\D}^{\red}_{\bullet,\hbar})$ denote the subcategories consisting of $S_n$-equivariant objects and morphisms. 
			
			\item 
			The category $\CP (\underline{\A}_{\hbar}^{\wedge 0})$ consists of finitely generated Poisson $\underline{\A}_{\hbar}^{\wedge 0}$-bimodules $\B_\hbar'$ that are 
			\begin{itemize}
				\item $\C^*$-equivariant with respect to the \textit{model} $\C^*$-action;
				\item complete and separated in the $\idealm$-adic topology;
				\item $\C^*$-prorational.
			\end{itemize} 
			Let $\CP _0(\underline{\A}_{\hbar}^{\wedge 0})$ be its full subcategory consisting of objects supported at $0$. 
			The symmetric group $S_n$ acts on $\underline{\A}_{\hbar}^{\wedge 0}$. Let $\CP^{S_n}(\underline{\A}_{\hbar}^{\wedge 0})$, $\CP^{S_n}_0(\underline{\A}_{\hbar}^{\wedge 0})$ denote the subcategories consisting of $S_n$-equivariant objects and morphisms. 		
			\item  
			The category ${^c \CP }(\underline{\A}_{\hbar}^{\wedge 0})$ consists of finitely generated Poisson $\underline{\A}_{\hbar}^{\wedge 0}$-bimodules $\B_\hbar'$ that are 
			\begin{itemize}
				\item $\C^*$-equivariant with respect to the \textit{conical} $\C^*$-action;
				\item complete and separated in the $\idealm$-adic topology; 
				\item $\C^*$-prorational.
			\end{itemize} 
			Let ${^c \CP_0}(\underline{\A}_{\hbar}^{\wedge 0})$ be its full subcategory consisting of objects supported at $0$. 
			The symmetric group $S_n$ acts on $\underline{\A}_{\hbar}$. 
			Let ${^c \CP^{S_n}}(\underline{\A}_{\hbar}^{\wedge 0})$ and ${^c \CP^{S_n}_0}(\underline{\A}_{\hbar}^{\wedge 0})$ be the subcategories consisting of $S_n$-equivariant objects and morphisms. 
			
			\item 
			The category $^c \CP(\underline{\A}_\hbar)$ consists of finitely generated graded (with respect to the conical $\C^*$-action) Poisson $\underline{\A}_{\hbar}$-bimodules, $^c\CP_0(\underline{\A}_\hbar)$ is its full subcategory consisting of finite dimensional objects. 
			Let $^c \CP^{S_n}(\underline{\A})$ and $^c\CP^{S_n}_0(\underline{\A}_\hbar)$ be the subcategories consisting of $S_n$-equivariant objects and morphisms. 
			\item The category $\HC(\sliceA)_F$ consists of Harish-Chandra bimodules over $\sliceA$ equipped with a specified good filtration. The group $S_n$ acts on $\sliceA$; the category $\HC^{S_n}(\sliceA)_F$ consists of $S_n$-equivariant Harish-Chandra bimodules over $\sliceA$; $\HC^{S_n}_0(\sliceA)_F$ is its full subcategory consisting of finite dimensional objects.
		\end{enumerate}
	\end{definition}
	We will produce functors $\F_i$ from the category ($i$) to the category ($i+1$), and partial inverse functors $\F'_i$ from the category ($i+1$) to the category ($i$), for $1\le i \le 10$, 
	\subsection{$\D$-module structure}
	In this subsection we collect a few preliminary results that comes from constructing $\D$-module structures from Poisson bimodule structures. They will be essential for constructing the functors between the categories  $\CP(\A_\hbar^{\wedge\overline{\L}}),  \CP(\A_\hbar^{\wedge\L})$, $\CP(\A_\hbar^{\wedge{\L_{I\cup J}}})$ and between $\CP (\underline{\A}_{ \hbar} ^ {\wedge 0}),  \CP^{S_n}(\widetilde{\D}^{\red}_{I\cup J,\hbar}) $. 
	
	Let $U$ be a symplectic vector space, and let $\C^*$ act on $U$ by scaling. Let $\WeylalgebraA_\hbar$ be the Weyl algebra of $U$, and $\WeylalgebraA_\hbar|_{U^0}$ be its restriction to $U^0$ as a sheaf of algebras equipped with a $\C^*$-action. Let $U^0$ be a $\C^*$-stable open subset of $U$. 
	\begin{lemma}[{\cite[Section 3.5]{losev2012completions}}]\label{structure of D modules from weyl algebra}
		Let $\B_\hbar$ be a weakly $\C^*$-equivariant Poisson $\WeylalgebraA_\hbar|_{U^0}$-bimodule, and suppose $\B_\hbar$ is equipped with a decreasing filtration $F_i \B_\hbar$, such that 
		\begin{itemize} 
			\item $\B_\hbar$ is complete with respect to this filtration; 
			\item $\hbar F_i \B_\hbar \subset F_{i+1}\B_\hbar$;
			\item each $F_i\B_\hbar / F_{i+1}\B_\hbar$ is a coherent $\sheafO_{U^0}$-module. 
		\end{itemize}
		Then 
		\begin{enumerate}
			\item $\B_\hbar$ is a $\D_{U^0}$-module; each $F_i\B_\hbar$ is a $\D_{U^0}$-submodule; each $\B_\hbar/F_i\B_\hbar$ is a coherent $\sheafO_{U^0}$ module;  
			\item If $\codim_{U}U\setminus U^0 \ge 2$, then the natural morphism
			\begin{equation}
				\WeylalgebraA_\hbar|_{U^0} \completeotimes_{\C[[\hbar]]} (\B_\hbar)^{\ad U} \to \B_\hbar
			\end{equation} 
			is an isomorphism. Here $(\B_\hbar)^{\ad U}$ denotes the kernel of $\{U,-\}$ in $\B_\hbar$. 
		\end{enumerate}		
	\end{lemma}
	
	\begin{remarklabeled}\label{construction of D module structure from weyl algebra structure}
		We recall from \Cite[Section 3.5]{losev2012completions} the construction of the $\D_{U^0}$-module structure in (1) of \Cref{structure of D modules from weyl algebra} since it will be useful later. 
		
		Denote the bimodule product map $\WeylalgebraA_\hbar|_{U^0} \otimes \B_\hbar \otimes \WeylalgebraA_\hbar|_{U^0} \to \B_\hbar$ by $1\otimes b \otimes a_2 \mapsto a_1 * b * a_2$. Let $u,v\in U$ and a local section $b$ of $\B_\hbar$. Let $\check{u} = \omega(u,-) \in U^*$. Let $\partial_v$ be the vector field on $U$ corresponding to $v$. Define
		\begin{equation}\label{definition of D module structure on U0}
			\begin{split}
				&\check{u}\cdot b = \frac{1}{2} (u*b + b*u); \\ 
				&\partial_v b = \{v,b\}
			\end{split}
		\end{equation}
		where we view $v$ as an element of $\WeylalgebraA_\hbar|_{U^0}$. 
		
		It is straightforward to check that the operators $u\cdot, v\cdot, \partial_u \cdot, \partial_v\cdot$ satisfy the relations 
		\begin{equation}
			\begin{split}
				&[\check{u}\cdot, \check{v}\cdot] = [\partial_u , \partial_v ] = 0, \\ 
				&[\partial_v,\check{u}\cdot] = \<\check{u},v\>\id.
			\end{split}
		\end{equation}
		Therefore, we get a $D_U(U)$-module structure on $\B_\hbar$. If $U_f\subset U^0$ is a distinguished open set of $U$ defined by the polynomial $f\in \C[U]$, then by the definition of $f\cdot$ we have $f\cdot \equiv f* \mod \hbar$. Therefore, since $f*$ is invertible on $\B_\hbar$ and $\B_\hbar$ is $\hbar$-adic complete, we see $f\cdot$ is also invertible on $\B_\hbar$. It follows that $\B_\hbar$ is a $D_{U^0}$-module. 
	\end{remarklabeled}

	To apply the above lemma in our situation, we need the following result on equivariant descent. 
	Let $\B_\hbar$ be an object of $\CP( {\D}^{\red}_{\bullet,\hbar})$, and define
	\begin{equation}
		\F(\B_\hbar) = \F( - ) := \widetilde{\D}^{\red}_{\bullet,\hbar} \otimes_{{\D}^{\red}_{\bullet,\hbar}} \B_\hbar. 
	\end{equation}
	A priori, $\F(\B_\hbar)$ is a finitely generated {left} $\widetilde{\D}^{\red}_{I\cup J,\hbar}$-module. Note that the $S_n$-action on $\tilde{\L}$ is free. 
	Therefore, by \cite[Lemma 3.6.3 (1)]{losev2012completions}, there is a Poisson bracket satisfying \Cref{definition of Poisson bimodules over quantization}
	\begin{equation}\label{want poisson bracket for descent}
		\{\bullet,\bullet\}: \widetilde{\D}^{\red}_{\bullet,\hbar} \otimes_{\C[[\hbar]]} \F(\B_\hbar') \to \F(\B_\hbar')
	\end{equation}
	that extends the Poisson bracket $ {\D}^{\red}_{\bullet,\hbar} \otimes_{\C[[\hbar]]} \B'_\hbar  \to \B'_\hbar $. 
	We can thus define functors 
	\begin{equation} 
		\F: \CP( {\D}^{\red}_{\bullet,\hbar}) \leftrightarrows \CP^{S_n}(\widetilde{\D}^{\red}_{\bullet,\hbar}) : \F',
	\end{equation}
	given by
	\begin{equation}\label{inverse to descend functor}
		\begin{split}
			\F( - ) &:= \widetilde{\D}^{\red}_{\bullet,\hbar} \otimes_{{\D}^{\red}_{\bullet,\hbar}} -,\\
			\F'(-) & := (-)^{S_n}. 
		\end{split}
	\end{equation}
	\cite[Lemma 3.6.3 (1)]{losev2012completions} also shows the following result. 
	\begin{proposition}\label{equivariant descent on D is equivalence}
		The functors $\F$ and $\F'$ are mutual quasi-inverses.
	\end{proposition}
	Let $\L_K$ be an affine open subset of $\L$ defined in \Cref{subsection good affine open cover}. 
	\begin{proposition}\label{Bhbar locally pro O coherent D module}
		Let $\B_\hbar'$ be an object of $\CP(\A_\hbar^{\wedge\L_K})$. Then
		\begin{enumerate}
			\item $\B_\hbar'$ is equipped with a pro-$\sheafO_{\L_K}$-coherent $\D_{\L_K}$-module structure, and for each $n\ge 1$, $\B_\hbar'/\idealI_{\hbar,K}^n\B_\hbar'$ is an $\sheafO_{\L_K}$-coherent $\D_{\L_K}$-module. 
			\item The $\sheafO_{\L_K}$-module structures on $\B_\hbar'/ \idealI_{I,\hbar}\B_\hbar' $ coming from $\A_\hbar^{\wedge\L_I}$ and the $\D_{\L_K}$-module structure in (1) agree. 
			More precisely, let $f \in \C[\L_K]$ and $\tilde{f}$ be a lift of $f$ in $\A_\hbar^{\wedge\L_K}$. Then, for any $b\in \B_\hbar'$, $\tilde{f}b + \idealI_{K,\hbar}\B_\hbar' = f\cdot b +\idealI_{K,\hbar}\B_\hbar'$, where $f\cdot $ denotes the $\D$-module multiplication. 
			\item The sheaf $\B_\hbar'/ \idealI_{K,\hbar}\B_\hbar' $, equipped with the $\sheafO_{\L_K}$-module structure induced from the $\A_{\hbar}^{\L_K}$-module structure, is a vector bundle over $\L_K$. 
		\end{enumerate}
	\end{proposition}
	\begin{proof}
		It suffices to prove the statements for $K=I$. 
		Thanks to the isomorphism $\vartheta_{I,\hbar}^{\red}$ from \eqref{isomorphism of completed quantum algebra over affine piece}, we may view $\B'_{\hbar}$ as a $\C^*$-equivariant (with respect to the model $\C^*$-action) coherent Poisson $\D^{\red}_{I,\hbar}$-bimodule. 
		Apply the functor $\F$ in \eqref{inverse to descend functor}. 
		Since $ \WeylalgebraA_\hbar|_{\tilde{\L}_I}$  embeds $\C^*$-equivariantly into $\widetilde{\D}^{\red}_{\bullet,\hbar}$, we see by \Cref{structure of D modules from weyl algebra} that $\F(\B_\hbar')$ is a $\D_{\tilde{\L}_I}$-module and each $\F(\B_\hbar') / \tilde{\idealI}_{I,\hbar}^n \F(\B_\hbar')$ is $\sheafO_{\tilde{\L}_I}$-coherent.
		Take $S_n$-invariants (i.e. apply the functor $\F'$), it follows that $\B_\hbar'$ is a  $\D_{\L_I}$-module such that each $\B_\hbar'/\idealI_{\hbar,I}^n\B_\hbar'$ is $\sheafO_{\L_I}$-coherent. This proves (1). 
		
		Let us prove (2). Let $f\in\C[\L_I]\subset \C[\tilde{\L}_I]$. Let $\tilde{f}$ be an $S_n$-invariant lift of $f$ in $\WeylalgebraA_\hbar|_{\tilde{\L}_I}$, so $\tilde{f}\otimes 1 \in \D^{\red}_{I,\hbar}$. 
		By the construction \Cref{construction of D module structure from weyl algebra structure}, we see that for any $b\in \B_\hbar'$, $f\cdot (b+\idealI_{\hbar,I}\B_\hbar') = \tilde{f} b + \idealI_{\hbar,I}\B_\hbar'$, and 
		(2) follows. 
		
		(3) follows from (2) and the fact that an $\sheafO$-coherent $\D$-module is a vector bundle. 
	\end{proof}
	
	The following result is a corollary of {\cite[Cor. VIII.2.3]{grothendieck1968sga}}. 
	\begin{theorem}\label{grothendieck finiteness}
		Let $X$ be a variety, $Y\subset X$ be a closed subvariety, and $U = X\setminus Y$. Let $\iota: U \to X$ denote the inclusion. Let $\mathcal{E}$ be a coherent $\sheafO_X$-module. 
		If $U$ is smooth, $\mathcal{E}|_U $ is a vector bundle, and $\codim_X Y \ge n$, then for any $i \le n-2$, $R^i \iota_* (\mathcal{E}|_U)$ is a coherent $\sheafO_X$-module. 
	\end{theorem}
	We conclude this subsection with the following result about sheaves annihilated by some power of $\idealI_\hbar$. 
	\begin{corollary}\label{finiteness of H1 and H2 for vector bundles}
		Let $\B_\hbar'$ be an object of $\CP(\A_\hbar^{\wedge\L})$. Suppose that $\idealI_{\hbar}^k\B_\hbar' = 0$ for some $k\ge 1$. Then the following are true. 
		\begin{enumerate}
			\item For $i=0,1,2$, $H^i (\L, \B_\hbar')$  is finitely generated as a left $\A_\hbar^{\wedge\overline{\L}}$-module. 
			\item Suppose $f\in \C[\overline{\L}]$ and $\hat{f}$ is a lift of $f$ in $\A_\hbar^{\wedge\overline{\L}}$. Suppose the distinguished affine open subset $U_f \subset \L$. Then $H^i (\L, \B_\hbar')[\hat{f}\inverse] = 0$ for $i=1,2$. 
			\item For $f$ in (2), $H^0(U_f,\B_\hbar') = H^0(\L,\B_\hbar') [\hat{f}\inverse]$. 
		\end{enumerate}
	\end{corollary}
	\begin{proof}
		We prove all statements by induction on $k$. 
		
		We prove (1).
		Assume $k=1$. 
		Consider the affine open cover $\L_K$ of $ \L$, and apply \Cref{Bhbar locally pro O coherent D module} to each $\B_\hbar'|_{\L_K}$. Thanks to (3) of that proposition, we see $\B_\hbar' = \B_\hbar' / \idealI_{ \hbar}\B_\hbar'$ is a vector bundle on $\L$. 
		By \Cref{boundary of minimal leaf has codim 4}, $\codim_{\overline{\L}}\d\L \ge 4$;
		by \Cref{grothendieck finiteness}, (1) holds for $k=1$. 
		
		Now assume (1) holds for $k$ and $\B_\hbar'$ is annihilated by $\idealI_{ \hbar}^{k+1}$.Then the short exact sequence 
		\begin{equation}
			0 \to \idealI_{\hbar}^{k}\B_\hbar' / \idealI_{\hbar}^{k+1}\B_\hbar'  \to \B_\hbar' / \idealI_{\hbar}^{k+1}\B_\hbar' \to \B_\hbar' / \idealI_{\hbar}^{k}\B_\hbar' \to 0
		\end{equation}
		is just 
		\begin{equation}\label{simplified ses when B killed by I to power n+1}
			0 \to \idealI_{\hbar}^{k}\B_\hbar'  \to \B_\hbar'   \to \B_\hbar' / \idealI_{\hbar}^{k}\B_\hbar' \to 0
		\end{equation}
		and the first and third terms are annihilated by $\idealI_{\hbar}^{k}$. 
		Take the long exact sequence for $R\Gamma({\L},-)$ associated to the above short exact sequence, we get 
		\begin{equation}\label{part of les of H1 H2}
			\cdots H^i(\L, \idealI_{\hbar}^{k}\B_\hbar') \to H^i(\L, \B_\hbar') 
			\to H^i(\L,\B_\hbar' /  \idealI_{\hbar}^{k}\B_\hbar') \to \cdots
		\end{equation}
		The first and third terms are finitely generated for $i=0,1,2$ by the induction hypothesis, hence so is the middle term. 
		We proved (1) for $k+1$. 
		
		We now prove (2). When $k=1$, we view $\B_\hbar' $ as an $\sheafO_{\L}$-module, and view $H^i(\L,\B_\hbar' ) $ as $\C[\overline{\L}]$-modules. Since $\overline{\L}$ is affine, the long exact sequence for cohomology with support shows 
		\begin{equation}
			H^i(\L,\B_\hbar') \cong H^{i+1}_{\partial\L} (\L,\B_\hbar')
		\end{equation}
		for $i\ge 1$, and the right hand side is clearly supported on $\partial\L$. When $i=1,2$, $H^i(\L,\B_\hbar')$ is finitely generated by part (1). Therefore, $H^i (\L, \B_\hbar')[\hat{f}\inverse] = 0$ if $U_f\subset \L$. We have proved (2) for $k=1$. 
		
		In general, if $\B_\hbar'$ is annihilated by $\idealI_{ \hbar}^{k+1}$, then apply the exact functor $-[\hat{f}\inverse]$ to the long exact sequence \eqref{part of les of H1 H2}, and we get (2) for all $k$ by induction hypothesis. 
		
		We now prove (3). Suppose $k=1$ so that $\B_\hbar'$ is a coherent $\sheafO_\L$-module; by (1), $H^0(\L,\B_\hbar')$ is finitely generated over $\C[\overline{\L}]$, and may be viewed as a coherent sheaf over $\sheafO_{\overline{\L}}$ that extends $\B_\hbar'$. It follows that $H^0(\L,\B_\hbar')[f\inverse] = H^0(U_f,\B_\hbar')$. This finishes the proof of (3) for $k=1$. 
		
		Suppose $\idealI_{\hbar}^{k+1}\B_\hbar' = 0$. Apply $R\Gamma(U_f,-)$ and $-[\hat{f}\inverse] \circ R\Gamma(\L,-)$ to \eqref{simplified ses when B killed by I to power n+1}, and note that there are natural morphisms $H^0(\L,-)[\hat{f}\inverse] \to H^0(U_f,-)$. We get the following morphism of short exact sequences: 
		\begin{center}
			\begin{tikzcd}[sep=small]
				0 \arrow[r] & H^0(U_f, \idealI_{\hbar}^{k}\B_\hbar') \arrow[r] \arrow[d]& H^0(U_f, \B_\hbar')\arrow[r] \arrow[d]
				& H^0(U_f, \B_\hbar' /  \idealI_{\hbar}^{k}\B_\hbar') \arrow[r]\arrow[d] &  0 \\ 
				0 \arrow[r] & H^0(\L, \idealI_{\hbar}^{k}\B_\hbar')[\hat{f}\inverse] \arrow[r] & H^0(\L, \B_\hbar') [\hat{f}\inverse] \arrow[r] 
				& \arrow[r] H^0(\L,\B_\hbar' /  \idealI_{\hbar}^{k}\B_\hbar') [\hat{f}\inverse] \arrow[r]& 0
			\end{tikzcd}.
		\end{center}
		The first row is exact because $U_f$ is affine; the second row is exact by (2). The first and third vertical arrows are isomorphisms by induction hypothesis, so by the 5-lemma, the middle vertical arrow is also an isomorphism. This proves (3) for $k+1$. 
	\end{proof}
	We will only use the finite-generatedness and the support condition of $H^1$ in the remainder of this paper. 
	
	\subsection{Various functors}\label{subsection various functors}
	In this subsection, we construct functors between the categories defined in \Cref{subsection various categories}, and examine their properties. 
	\subsubsection{Rees module}\label{subsubsection functor taking Rees module}
	Let $\B$ be an object of $\HC(\A)_F$. As we mentioned in \Cref{subsection on coherent poisson bimodule over algebra}, for a good filtration of $\B$ the Rees bimodule 
	$\B_\hbar= \bigoplus_{i=0}^\infty \B_{\le i} \hbar^i$ is an object of $\CP(\A_\hbar)$. Let $\F_1: \HC(\A)_F \to \CP(\A_\hbar)$ be the functor of taking Rees bimodules. 
	Conversely, let $\F_1':  \CP(\A_\hbar) \to \HC(\A)_F $ be the functor $\B_\hbar\mapsto \B_\hbar/(\hbar-1)\B_\hbar$; the resulting bimodule is equipped with a filtration induced from the $\C^*$-action on $\B_\hbar$.
	
	It is clear that $\F_1(\B)$ is $\C[\hbar]$-flat, and $\F_1'\circ\F_1$ is the identity functor on $\HC(\A)_F$. 
	
	\subsubsection{Completion along $\overline{\L}$.}\label{functor completion along L closure}
	Let $\F_2: \CP(\A_\hbar) \to \CP(\A_\hbar^{\wedge \overline{\L}})$ be the functor of taking completion along the closed subset $\overline{\L}$. It restricts to a functor 
	$\CP_{\overline{\L}}(\A_\hbar) \to \CP_{\overline{\L}}(\A_\hbar^{\wedge \overline{\L}})$, 
	which we still denote by $\F_2$. 
	\begin{proposition}
		The functor $\F_2: \CP_{\overline{\L}}(\A_\hbar) \to \CP_{\overline{\L}}(\A_\hbar^{\wedge \overline{\L}})$ is an equivalence. A quasi-inverse $\F_2'$ is given by taking $\C^*$-locally finite submodules. 
	\end{proposition}
	\begin{proof}
		Note that the ideal $\idealI_{\overline{\L},\hbar}$ is positively graded. The proposition is then proved similarly to \cite[Lemma 3.3.1]{losev2011finite}. 
	\end{proof}
	
	\subsubsection{Restriction to $\L$.}\label{subsection functor F3}
	Let $j: \L \to \overline{\L}$ be the open inclusion. 
	Consider the functor
	\begin{equation}
		\F_3:  \CP(\A_\hbar^{\wedge\overline{\L}}) \to \CP(\A_\hbar^{\wedge\L}) 
	\end{equation}
	of restriction to $\L$. Let $\F_3': \CP(\A_\hbar^{\wedge\L})  \to \CP(\A_\hbar^{\wedge\overline{\L}})$ 
	send $\B_\hbar'$ to $\Loc_{\overline{\L}} ( H^0(\L, \B_\hbar'))$. 
	They restrict to functors between $\CP_{\overline{\L}}(\A_\hbar^{\wedge\overline{\L}})$ and $\CP_{\L}(\A_\hbar^{\wedge\L}) $.

	Let $\B_\hbar'$ be an object of $\CP_{\L}(\A_\hbar^{\wedge\L})$. 
	We start with the following important observation. 
	\begin{lemma}\label{lemma hbar and I not too different}
		There is some $k > 0 $ such that 
		\begin{equation}
			\hbar^k \B_\hbar' \subset \idealI_\hbar^k \B_\hbar' \subset \hbar \B_\hbar'\subset \idealI_\hbar \B_\hbar'. 
		\end{equation}
	\end{lemma}
	\begin{proof}
		The first and third inclusions holds because $\hbar \in \idealI_\hbar$. The second inclusion is a consequence of finite generation and the definition of support. 
	\end{proof}
	\begin{theorem}\label{restriction of localization is id over codim 4}
		Let $\B_\hbar'$ be an object of $\CP_{\L}(\A_\hbar^{\wedge\L})$.
		Then 
		\begin{equation}
			\F_3 \circ \F_3' (\B_\hbar') = \Loc_{\overline{\L}} ( H^0(\L, \B_\hbar'))|_{\L} =\B_\hbar'.
		\end{equation}
	\end{theorem}
	
	\begin{proof}
		We prove the theorem in several steps. \\
		\textit{Step 1. } Let $f \in \C[\overline{\L}]$ and $U_f$ be the distinguished affine open subset of $\overline{\L}$. Assume $U_f \subset \L$. Let $\hat{f}$ be a lift of $f$ in $\A_\hbar^{\wedge \overline{\L}}$. By the definition, 
		\begin{equation}
			\Loc_{\overline{\L}} ( H^0(\L, \B_\hbar'))(U_f) = \varprojlim_n H^0(\L, \B_\hbar') / \idealI^n_{\overline{\L},\hbar} H^0(\L, \B_\hbar') [\hat{f}\inverse].
		\end{equation}
		So our goal is to show 
		\begin{equation}\label{goal of microloc equal original Iadic version}
			H^0(U_f, \B_\hbar') = \varprojlim_n H^0(\L, \B_\hbar') / \idealI^n_{\overline{\L},\hbar} H^0(\L, \B_\hbar') [\hat{f}\inverse].
		\end{equation}
		Thanks to \Cref{lemma hbar and I not too different}, we have 
		\begin{equation}
			\idealI^k_{\overline{\L},\hbar} H^0(\L, \B_\hbar') \subset H^0(\L, \idealI^k_{\hbar} \B_\hbar') \subset H^0(\L,  {\hbar} \B_\hbar') = \hbar H^0(\L, \B_\hbar'). 
		\end{equation}
		Since $\hbar \A_\hbar^{\wedge \overline{\L}} \subset \idealI_{\overline{\L},\hbar}$, the $\hbar$-adic topology and the $\idealI_{\overline{\L},\hbar}$-adic topology on $H^0(\L, \B_\hbar')$ are equivalent. Therefore, to prove \eqref{goal of microloc equal original Iadic version}, it suffices to show 
		\begin{equation}\label{goal of microloc equal original hbar version}
			H^0(U_f, \B_\hbar') 
			= \varprojlim_n H^0(\L, \B_\hbar') / \hbar^n H^0(\L, \B_\hbar') [\hat{f}\inverse].
		\end{equation}
		\\
		\textit{Step 2. } 
		The short exact sequence 
		\begin{equation}
			0 \to \hbar^n\B_\hbar' \to \B_\hbar' \to \B_\hbar'/\hbar^n \B_\hbar' \to 0
		\end{equation} 
		induces a long exact sequence 
		\begin{equation}
			0 \to \hbar^n H^0(\L, \B_\hbar') \to H^0(\L,\B_\hbar' ) \to H^0(\L, \B_\hbar'/\hbar^n \B_\hbar') \to \hbar^n H^1 (\L, \B_\hbar') \to \cdots .
		\end{equation}
		Apply the exact functor $-[\hat{f}\inverse]$ to the above exact sequence, we get an inclusion
		\begin{equation}
			H^0(\L, \B_\hbar') / \hbar^n H^0(\L, \B_\hbar') [\hat{f}\inverse] \hookrightarrow H^0(\L, \B_\hbar'/\hbar^n \B_\hbar')[\hat{f}\inverse]
		\end{equation} 
		whose cokernel lies in $\hbar^n H^1 (\L, \B_\hbar')[\hat{f}\inverse]$. 
		
		By \Cref{lemma hbar and I not too different}, $\B_\hbar' / \hbar\B_\hbar'$ is annihilated by $\idealI_\hbar^k$, so by \Cref{finiteness of H1 and H2 for vector bundles}, $H^1(\L,\B_\hbar' / \hbar\B_\hbar')$ is finitely generated.
		
		Then, the argument of \cite[Lemma 5.6.3]{gordon2014category} shows that $H^1( \L,\B_\hbar')$ is also finitely generated. The argument in the proof of part (2) of \Cref{finiteness of H1 and H2 for vector bundles} shows that $ H^1 (\L, \B_\hbar')[\hat{f}\inverse] = 0$. Therefore, we have 
		\begin{equation}
			H^0(\L, \B_\hbar') / \hbar^n H^0(\L, \B_\hbar') [\hat{f}\inverse] \xrightarrow{\sim} H^0(\L, \B_\hbar'/\hbar^n \B_\hbar')[\hat{f}\inverse].
		\end{equation}
		We have reduced proving the theorem to proving the following 
		\begin{equation}\label{goal of microloc equal original B mod hbar version}
			H^0(U_f, \B_\hbar') 
			= \varprojlim_n H^0(\L, \B_\hbar'/\hbar^n \B_\hbar')[\hat{f}\inverse]. 
		\end{equation}
		\\
		\textit{Step 3. } We prove \eqref{goal of microloc equal original B mod hbar version}. By \Cref{lemma hbar and I not too different}, $\B_\hbar'/\hbar^n \B_\hbar'$ is annihilated by some power of $\idealI_\hbar$. Therefore, by (3) of \Cref{finiteness of H1 and H2 for vector bundles}, $H^0(\L, \B_\hbar'/\hbar^n \B_\hbar')[\hat{f}\inverse] = H^0(U_f, \B_\hbar'/\hbar^n \B_\hbar')$. 
		Since $\B_\hbar'$ is $\hbar$-adic complete, $\varprojlim_{n} H^0(U_f, \B_\hbar'/\hbar^n \B_\hbar') =  H^0(U_f, \B_\hbar')$. \eqref{goal of microloc equal original B mod hbar version} follows. This finishes the proof of the theorem. 
	\end{proof}
	\begin{corollary}
		Let $\B_\hbar'$ be a $\C[\hbar]$-flat object of $\CP_{\overline{\L}}(\A_\hbar^{\wedge\overline{\L}})$ and assume that $\B'_\hbar$ is $\C[\hbar]$-flat. 
		Then
		\begin{enumerate}
			\item the kernel $K_\hbar$ and cokernel $C_\hbar$ of the adjunction unit
			\begin{equation}\label{natural map B to H0(L,B)}
				\B'_\hbar \to \F_3'\circ \F_3 (\B'_\hbar)
			\end{equation}
			are supported on $\d\L$;
			\item the following morphism induced from \eqref{natural map B to H0(L,B)}
			\begin{equation}\label{natural map B to H0(L,B) mod h-1}
				\F_1'\circ\F_2' (\B_\hbar') \to \F_1'\circ\F_2'\circ \F_3'\circ \F_3 (\B'_\hbar)
			\end{equation}
			is an isomorphism. 
		\end{enumerate}
		
	\end{corollary}
	\begin{proof}
		Part (1) follows by applying \Cref{restriction of localization is id over codim 4} to $	\B'_\hbar|_{\L}$. 
		
		For part (2), note first  $\F_1'\circ\F_2'(K_\hbar)$ and $\F_1'\circ\F_2'(C_\hbar)$ are Harish-Chandra $\A$-bimodules supported on $\partial\L$. But $\overline{\L}$ is the minimal possible support of a nonzero Harish-Chandra bimodule, so $\F_1'\circ\F_2'(K_\hbar)=0$ and $\F_1'\circ\F_2'(C_\hbar)=0$.

		Note that $\F_1'$ is the tensor functor $- \otimes_{\C[\hbar]}\C[\hbar]/ (\hbar-1)$. 
		Apply $\F_1'$ to the exact sequence 
		\begin{equation}
			\F_2' (\B_\hbar') \to \F_2'\circ \F_3'\circ \F_3 (\B'_\hbar) \to \F_2'(C_\hbar) \to 0,
		\end{equation}
		we see \eqref{natural map B to H0(L,B) mod h-1} is surjective. 
		
		Let $\mathcal{N}_\hbar'$ denote the image of $\F_2' (\B_\hbar') \to  \F_2'\circ \F_3'\circ \F_3 (\B'_\hbar)$, which is $\C[\hbar]$-flat. Apply $\F_1'$ to the short exact sequence 
		\begin{equation}
			0\to \F_2' (K_\hbar) \to \F_2' (\B_\hbar') \to \mathcal{N}_\hbar' \to 0
		\end{equation}
		and use $\C[\hbar]$-flatness of $\mathcal{N}_\hbar'$, we see \eqref{natural map B to H0(L,B) mod h-1} is injective. 
	\end{proof}
	\subsubsection{Restriction to $\L_{I\cup J}$.}
	Let $\iota: \L_{I\cup J} \to \L$ be the open inclusion. 
	\begin{proposition}\label{extension to codim 2 is coherent}
		If $\B_\hbar'$ is an object of $\CP_{\L_{I\cup J}}(\A_\hbar^{\wedge{\L_{I\cup J}}})$, then $\iota_*\B_\hbar'$ is an object of $\CP_{\L}(\A_\hbar^{\wedge{\L}}) $. 
	\end{proposition}
	\begin{proof}
		This is a local statement, so we focus on open inclusions $\iota_{K}: \L_{I\cup J} \cap \L_K \to \L_K$, where $\L_K$ is the open affine defined in \Cref{subsection good affine open cover}. 
		To simplify notation, in this proof we write $X = \L_K$, $X_0 = \L_K \cap (\L_I \cup \L_J)$, and $\iota$ for $\iota_K$. Write $\B_\hbar'$ for $\B_\hbar'|_X$. 
		Let $\idealI_{\hbar,0}$ denote the ideal sheaf of $X_0$ in $ \A_{\lambda,\hbar}^{\wedge X_0}$. 
		We need to show: 
		\begin{enumerate}
			\item $H^0(X_0, \B_\hbar')$ is finitely generated over $ \A_{\lambda,\hbar}^{\wedge\L_K}$;
			\item $\iota_*(\B_\hbar') $ coincides with the microlocalization of $H^0(X_0,  \B_\hbar')$. 
		\end{enumerate} 
		
		To prove both claims, let us first consider the short exact sequence
		\begin{equation}
			0 \to \idealI_{\hbar,0}^n\B_\hbar' \to \B_\hbar' \to \B_\hbar'/\idealI_{\hbar,0}^n\B_\hbar' \to 0
		\end{equation}
		which induces an injection
		\begin{equation}\label{H0B quotient H0IkB inject into H0Bquotient IkB}
			H^0(X_0,\B_\hbar') / H^0(X_0,\idealI_{\hbar,0}^n\B_\hbar') 
			\to 
			H^0(X_0,\B_\hbar'/\idealI_{\hbar,0}^n\B_\hbar').
		\end{equation}
		Thanks to \Cref{Bhbar locally pro O coherent D module}, for any pair $n \ge m$, the natural surjective morphism
		$\B_\hbar'/\idealI_{\hbar,0}^n\B_\hbar' \to \B_\hbar'/\idealI_{\hbar,0}^m\B_\hbar'$ is a morphism of $\sheafO_{X_0}$-coherent $\D_{X_0}$-modules. 
		Recall $\codim_{ X}  X \backslash X_0 = 2$.
		Therefore, all $\iota_*(\B_\hbar'/\idealI_{\hbar,0}^n\B_\hbar') $ are $\sheafO_X$-coherent $\D_X$-modules. 
		We thus get a surjective morphism 
		\begin{multline}
			H^0(X, \iota_*( \B_\hbar'/\idealI_{\hbar,0}^n\B_\hbar') ) \cong	\\
			H^0(X_0, \B_\hbar'/\idealI_{\hbar,0}^n\B_\hbar') \to  H^0(X_0, \B_\hbar'/\idealI_{\hbar,0}^m\B_\hbar')
			\\ \cong H^0(X, \iota_*( \B_\hbar'/\idealI_{\hbar,0}^m\B_\hbar') ) . 
		\end{multline}
		Note also that $H^0(X_0,\B_\hbar') = \varprojlim_n H^0(X_0,\B_\hbar'/\idealI_{\hbar,0}^n\B_\hbar')$. It follows that $H^0(X_0,\B_\hbar') \twoheadrightarrow H^0(X_0,\B_\hbar'/\idealI_{\hbar,0}^n\B_\hbar')$ for any $n$, and \eqref{H0B quotient H0IkB inject into H0Bquotient IkB} is an isomorphism. 
		
		Let us prove (1). Thanks to \Cref{lemma hbar and I not too different}, for any $n$, we have a natural surjection 
		\begin{multline}
			H^0(X_0, \B_\hbar' /\idealI_{\hbar,0} ^{nk}\B_\hbar') \cong
			H^0(X_0,\B_\hbar') /  H^0(X_0,\idealI_{\hbar,0} ^{nk}\B_\hbar')  \\
			\twoheadrightarrow
			H^0(X_0,\B_\hbar') / \hbar^n H^0(X_0,\B_\hbar')
		\end{multline}
		where the first isomorphism is \eqref{H0B quotient H0IkB inject into H0Bquotient IkB}. Since $H^0(X_0, \B_\hbar' \idealI_{\hbar,0} ^{nk}\B_\hbar') $ is finitely generated over $\A_\hbar^{\wedge X}$, so is $H^0(X_0,\B_\hbar') / \hbar^n H^0(X_0,\B_\hbar')$. Since $H^0(X_0,\B_\hbar') $ is $\hbar$-adic complete and separated, we conclude that it is finitely generated over $\A_\hbar^{\wedge X}$. 
		
		Let us prove (2). 
		We need to show for any $f\in \C[X]$, 
		\begin{equation}\label{goal property of coherence}
			H^0(X_f, \iota_*\B_\hbar ) = 	
			H^0(X,  \iota_{*}{\B}_{\hbar}' )[f\inverse]
		\end{equation} 
		where the right hand side is, by definition, 
		\begin{equation}\label{definition Iadic localization}
			\varprojlim_n 
			\left( H^0(X_0,  \iota_{*}{\B}_{\hbar}' )/ \idealI_{\hbar}^n H^0(X,  \iota_{*}{\B}_{\hbar}' ) [f\inverse] \right) .
		\end{equation} 
		By \Cref{lemma hbar and I not too different}, the topology on $H^0(X_0,\B_\hbar')$ given by the filtration $\idealI_\hbar^n H^0(X_0,\B_\hbar') $ and the filtration $H^0(X_0,\idealI_{\hbar,0} ^n\B_\hbar') $ are equivalent -- and both are equivalent to the $\hbar$-adic topology. In particular, the natural morphism 
		\begin{equation}\label{two filtrations are equiv to the hbar adic filtration}
			\varprojlim_n 
			\left( H^0(X_0,  {\B}_{\hbar}' )/ \idealI_{\hbar}^n H^0(X_0,  {\B}_{\hbar}' ) [f\inverse] \right) \to
			\varprojlim_n 
			\left( H^0(X_0,   {\B}_{\hbar}' )/  H^0(X_0,\idealI_{\hbar,0} ^n\B_\hbar') [f\inverse] \right)
		\end{equation}
		is an isomorphism. 
		Therefore, we can identify the inverse limit \eqref{definition Iadic localization} with
		\begin{equation}
			\begin{split}
				\eqref{definition Iadic localization} &=  
				\varprojlim_n 
				\left( H^0(X_0,  {\B}_{\hbar}' )/ \idealI_{\hbar}^n H^0(X_0,  {\B}_{\hbar}' ) [f\inverse] \right)\\
				&\cong 
				\varprojlim_n 
				\left( H^0(X_0,   {\B}_{\hbar}' )/  H^0(X_0,\idealI_{\hbar,0} ^n\B_\hbar') [f\inverse] \right) \\
				&\cong \varprojlim_n 
				\left( H^0(X_0, {\B}_{\hbar}' / \idealI_{\hbar,0} ^n\B_\hbar') [f\inverse] \right)\\
				&\cong \varprojlim_n \left(
				H^0(X, \iota_* ({\B}_{\hbar}' / \idealI_{\hbar,0} ^n\B_\hbar')) [f\inverse] \right)\\
				&\cong \varprojlim_n 
				H^0(X_f, \iota_* ({\B}_{\hbar}' / \idealI_{\hbar,0} ^n\B_\hbar')) \\
				& \cong H^0(X_f, \iota_* {\B}_{\hbar}' )
			\end{split}
		\end{equation}
		which is precisely \eqref{goal property of coherence}. 
		The first equality follows from the definition of $\iota_*$; the second isomorphism follows since \eqref{two filtrations are equiv to the hbar adic filtration} is an isomorphism; 
		the third isomorphism follows since \eqref{H0B quotient H0IkB inject into H0Bquotient IkB} is an isomorphism;
		the forth isomorphism follows from the definition of $\iota_*$; 
		the fifth isomorphism follows since $\iota_* ({\B}_{\hbar}' / \idealI_{\hbar,0} ^n\B_\hbar')$ is $\sheafO_X$-coherent and $X$ is affine; 
		the last isomorphism follows since $\B_\hbar'$ is complete with respect to the $\idealI_{\hbar,0}$-adic topology, and $H^0(X_f,-)$ and $\iota_*$ commutes with inverse limit. 
		
		This finishes the proof of the proposition.  
	\end{proof}
	
	We define functors 
	\begin{equation}
		\F_4' : \CP_{\L_{I\cup J}}(\A_\hbar^{\wedge{\L_{I\cup J}}}) \leftrightarrows \CP_{\L}(\A_\hbar^{\wedge{\L}})  : \F_4. 
	\end{equation}
	given by $\F_4 = \iota^*$ and $\F_4' = \iota_*$. 
	It is clear that $\F_4$ is a left quasiinverse and left adjoint to $\F_4'$. By the next proposition, it is also a right quasiinverse. 
	\begin{proposition}
		Let $\B_\hbar'$ be an object of $ \CP_{\L}(\A_{\hbar} ^ {\wedge \L}) $. Then the adjunction counit map $\iota_*\iota^* \B_\hbar' \to \B_\hbar'$ is an isomorphism. 
	\end{proposition}
	\begin{proof}
		The statement is local, so we prove the corresponding statement for $\iota = \iota_K: X_0 \to X$ and an object $\B_\hbar'$ of $\CP_X(\A_{\hbar} ^ {\wedge X})$, as in the proof of \Cref{extension to codim 2 is coherent}. Since both $\iota_*\iota^* \B_\hbar'$ and $\B_\hbar'$ are coherent Poisson bimodules, it suffices to show the natural map 
		\begin{equation}\label{extension over codim2 is equiv requires H0 X0 = X}
			H^0(X, \B_\hbar') \to H^0(X_0,\B_\hbar')
		\end{equation}
		is an isomorphism. 
		Since each $\B_\hbar'/ \idealI_\hbar^n \B_\hbar'$ is an $\sheafO_X$-coherent $\D_X$-module and thus a vector bundle, $\codim_{ X}  X \backslash X_0 = 2 $, Hartogs's theorem implies $H^0(X_0,\B_\hbar'/ \idealI_\hbar^n \B_\hbar') \cong H^0(X,\B_\hbar'/ \idealI_\hbar^n \B_\hbar')$.
		Since $\B_\hbar'$ is $\idealI_\hbar$-adic complete, $H^0(X_0,\B_\hbar') = \varprojlim_n H^0(X_0,\B_\hbar'/ \idealI_\hbar^n \B_\hbar')$ while $H^0(X,\B_\hbar') = \varprojlim_n H^0(X,\B_\hbar'/ \idealI_\hbar^n \B_\hbar')$. This proves \eqref{extension over codim2 is equiv requires H0 X0 = X}.  
	\end{proof}
	\subsubsection{Regluing.} 
	Let $\B_\hbar'$ be an object of $ \CP(\A_{\hbar} ^ {\wedge \L_{I\cup J}}) $. 
	Restrict it to $\L_I,\L_J$ respectively, we get a finitely generated Poisson $(\A_\lambda) ^{\wedge \L_I} $ -bimodule $\B'_{\hbar,I}$ and a finitely generated Poisson $(\A_\lambda) ^{\wedge \L_J} $ -bimodule $\B'_{\hbar,J}$. 
	Thanks to the isomorphisms $\vartheta_{J,\hbar}^{\red}, \vartheta_{J,\hbar}^{\red}$ from \eqref{isomorphism of completed quantum algebra over affine piece}, 
	we get 
	\begin{equation}
		(  \vartheta_{I,\hbar}^{\red})_*\B'_{\hbar,I}	, \ (  \vartheta_{J,\hbar}^{\red})_*\B'_{\hbar,J}
	\end{equation}
	which are Poisson bimodules over $\D^{\red}_{I,\hbar}$, $\D^{\red}_{J,\hbar}$ respectively. 
	
	Let $F$ be as in \Cref{quantum inner automorphism}. Thanks to \eqref{prounipotency for quantum ad F}, $\exp\left(\frac{1}{\hbar^2} \ad F\right)$ is a well defined operator on $\B_{\hbar,I\cap J}'$, the restrictions of $\B_\hbar'$ to $\L_{I\cap J}$. 
	By (2) of \Cref{quantum inner automorphism}, the operator $\exp\left(\frac{1}{\hbar^2} \ad F\right)$ induces an $\D^{\red}_{I\cap J,\hbar}$-linear isomorphism 
	\begin{equation}
		(  \vartheta_{I,\hbar}^{\red})_*\B'_{\hbar,I} |_{I\cup J}	 \to (  \vartheta_{J,\hbar}^{\red})_*\B'_{\hbar,J}|_{I\cup J}.
	\end{equation}
	We can thus glue $(  \vartheta_{I,\hbar}^{\red})_*\B'_{\hbar,I}$ and $(  \vartheta_{J,\hbar}^{\red})_*\B'_{\hbar,J}$ into a sheaf of coherent Poisson $\D^{\red}_{I\cup J,\hbar}$-bimodules, i.e. an object of $\D^{\red}_{I\cup J,\hbar}$. 
	Let $\F_5:  \CP(\A_{\hbar} ^ {\wedge \L_{I\cup J}})  \to  \CP^{S_n}(\widetilde{\D}^{\red}_{I\cup J,\hbar})$ be the functor sending $\B'_\hbar$ to the above object. 
	Similarly, we get a functor $\F_5':  \CP^{S_n}(\widetilde{\D}^{\red}_{I\cup J,\hbar}) \to \CP(\A_{\hbar} ^ {\wedge \L_{I\cup J}})$ via regluing by $\exp (-\frac{1}{\hbar^2} \ad F)$. The functors $\F_5$ and $\F_5'$ are mutual quasi-inverse, and restricts to an equivalence of categories 
	$ \CP_{\wt{\L}_{I\cup J}}(\A_{\hbar} ^ {\wedge \L_{I\cup J}})  \cong \CP_{\wt{\L}_{I\cup J}}^{S_n}(\widetilde{\D}^{\red}_{I\cup J,\hbar})$. 
	
	\subsubsection{Equivariant descent.} 
	Write $\F_6 $ for $\F$ in \Cref{equivariant descent on D is equivalence} and $\F_6'$ for $\F'$ there. We have seen that $\F_6,\F_6'$ are mutual quasi inverses between the categories $\CP(\D^{\red}_{I\cup J,\hbar}) $ and $ \CP^{S_n}(\widetilde{\D}^{\red}_{I\cup J,\hbar})$, which restrict to quasi inverses between the categories $\CP_{{\L}_{I\cup J}}(\D^{\red}_{I\cup J,\hbar}) $ and $ \CP_{\wt{\L}_{I\cup J}}^{S_n}(\widetilde{\D}^{\red}_{I\cup J,\hbar})$.
	
	\subsubsection{Taking flat sections.}
	
	Let $\mathcal{B}_\hbar' \in \CP^{S_n}( \widetilde{\D}^{\red}_{I\cup J,\hbar})$. 
	By an analogue of \cite[Lemma 2.4.4]{losev2011finite}, its filtration given by
	\begin{equation}
		F_k\mathcal{B}_\hbar': = (\idealI_{\hbar,I\cup J})^ k \mathcal{B}_\hbar'
	\end{equation} 
	is complete and separated. Moreover, each $F_k\mathcal{B}_\hbar'/F_{k+1}\mathcal{B}_\hbar'$ is a coherent $\sheafO_{\wt{\L}_{\hbar,I\cup J}}$-module. 
	
	Define $\F_7(\mathcal{N}_\hbar)$ to be the Poisson centralizer of $\WeylalgebraA_\hbar(\wt{\L}_{I\cup J })$ in $\mathcal{N}_\hbar$. 
	Let $\F_7' : \CP (\underline{\A}_{ \hbar} ^ {\wedge 0}) \to \CP^{S_n}(\widetilde{\D}^{\red}_{I\cup J,\hbar}) $ by $\mathcal{M}_\hbar \mapsto \WeylalgebraA_\hbar(\wt{\L}_{I\cup J }) \completeotimes \M_\hbar$. 
	The functors $\F_7$ and $\F_7'$ are mutual quasi-inverses, thanks to (2) of \Cref{structure of D modules from weyl algebra}. 
	
	\subsubsection{Modifying $\C^*$-action.} 
	Recall we have two $\C^*$-actions on $\sliceA_\hbar$, the model action and the conical action. 
	Let $\alpha_\hbar$ be as in \eqref{definition of alpha relating two torus actions}.
	Let $\B_\hbar'$ be an object of $\CP(\underline{\A}_{\hbar} ^ {\wedge 0})$.
	We are going to equip $\B_\hbar'$ with a $\C^*$-action compatible with the conical $\C^*$-action on $\underline{\A}_{\hbar}$. 
	
	Denote by $\Eu_m$ the Euler derivation on $\B_\hbar'$ induced by the $\C^*$-action compatible with the model $\C^*$-action on $\sliceA_\hbar$. 
	Then, thanks to (2) of \Cref{relation between model and conical Cstar action}, $\Eu_c := \Eu_m-\frac{1}{\hbar^2}[\alpha_\hbar,-]$ is a derivation on $\B_\hbar'$ compatible with the conical $\C^*$-action on  $\underline{\A}_{\hbar}^{\wedge 0}$. It induces a derivation $\Eu_{c,k}$ on each $\B_\hbar' / \idealm^k \B_\hbar'$. We will upgrade this derivation to a $\C^*$-action. 
	\begin{proposition}\label{proposition modify c* action}
		There is a $\C^*$-action on $\B_\hbar'$ that is compatible with respect to the conical $\C^*$-action on $\sliceA_\hbar^{\wedge0}$. 
	\end{proposition}
	\begin{proof}
		We construct such an action in several steps. 	
		\begin{enumerate}[labelsep=0cm,align=left,leftmargin=0cm, label=\textit{Step} \arabic*. \  ,itemindent=0.5cm]
			\item Fix $k\ge 0, n\in\ZZ$. Let $(\B_\hbar'/\idealm^k\B_\hbar')_n$ denote the $\Eu_m$-degree $n$ component of $\B_\hbar'/\idealm^k\B_\hbar'$. We show that $(\B_\hbar'/\idealm^k\B_\hbar')_n$ is finite dimensional. 
			
			Pick a finite dimensional subspace $V_k \subset \B_\hbar'/\idealm^k\B_\hbar'$ that is $\Eu_m$-stable and generates $\B_\hbar'/\idealm^k\B_\hbar'$ as an $\sliceA_\hbar^{\wedge0}$-bimodule. Write $V_k = \bigoplus_{d\in \ZZ} (V_k)_d$, where $(V_k)_d$ is the $\Eu_m$-degree $d$ part of $V_k$; all but finitely $(V_k)_d$ are nonzero. Note that for any $d'\in \ZZ$, the $\eu_m$-degree $d'$ part $(\sliceA_\hbar^{\wedge0}/\idealm^k)_{d'}$ is finite dimensional. Therefore, 
			\begin{equation}
				(\B_\hbar'/\idealm^k\B_\hbar')_n = \bigoplus_{d+d'+d'' = n} 
				(\sliceA_\hbar^{\wedge0}/\idealm^k)_{d'} 
				(V_k)_d (\sliceA_\hbar^{\wedge0}/\idealm^k)_{d''} 
			\end{equation}
			is finite dimensional. This proves our claim. 
			
			\item We claim that $\Eu_{c,k}$ is a locally finite derivation on $\B_\hbar'/\idealm^k\B_\hbar'$. Recall $\Eu_{c,k} = \Eu_m-\frac{1}{\hbar^2}[\alpha_\hbar,-]$. Since $\alpha_\hbar$ has degree 2 with respect to the model $\C^*$-action, it follows that the operator $\frac{1}{\hbar^2}[\alpha_\hbar,-]$, and hence $\Eu_{c,k} $, stabilizes all $(\B_\hbar'/\idealm^k\B_\hbar')_n$. 
			Therefore, $\Eu_{c,k}$ is a sum of linear endomorphisms on finite dimensional subspaces $(\B_\hbar'/\idealm^k\B_\hbar')_n$. Since $\B_\hbar'/\idealm^k\B_\hbar' = \bigoplus_n (\B_\hbar'/\idealm^k\B_\hbar')_n$, we see $\Eu_{c,k}$ is locally finite. 
			
			\item Let $\frac{1}{\hbar^2}[\alpha_\hbar,-]$ act on a vector space $V$. For $z\in\C$, let $V^z$ denote the $z$-generalized eigenspace of  $\frac{1}{\hbar^2}[\alpha_\hbar,-]$. Since $\frac{1}{\hbar^2}[\alpha_\hbar,-]$ preserves all the subspaces $(\B_\hbar'/\idealm^k\B_\hbar')_n$ we have
			\begin{equation}
				\begin{split}
					\B_\hbar'/\idealm^k\B_\hbar' &=  \bigoplus_{n\in \ZZ} (\B_\hbar'/\idealm^k\B_\hbar')_n\\
					&=\bigoplus_{z\in \C, n\in \ZZ} (\B_\hbar'/\idealm^k\B_\hbar')^z_n 
				\end{split}			
			\end{equation}
			It is clear that $\Eu_{c,k}$ acts on $ (\B_\hbar'/\idealm^k\B_\hbar')^z_n$ by generalized eigenvalue $n-z$. 
			Let $f: \C\to \ZZ$ be any function such that $f(z+1)=f(z)+1$; for example, we can take $f = \lfloor \Re(z)\rfloor$. 
			Let $\Eu_{c,k}'$ be the the operator that acts by the scalar $n-f(z)$ on  $(\B_\hbar'/\idealm^k\B_\hbar')^z_n$. 		
			It follows that $\Eu_{c,k}'$ integrates to a $\C^*$-action on $\B_\hbar'/\idealm^k\B_\hbar'$. We write this $\C^*$-action as $(t,b_k) \mapsto t*_k b_k$ for $t\in\C^*, b_k\in \B_\hbar'/\idealm^k\B_\hbar'$. 
			
			Since all the conical $\C^*$-eigenvalues of $\sliceA_\hbar^{\wedge0} / \idealm^k$ are integers and $\Eu_{c,k}$ is compatible with the conical $\C^*$-action, so is $\Eu_{c,k}'$. In terms of $\C^*$-actions, for any $a\in \sliceA_\hbar^{\wedge0}$ and $b_k\in \B_\hbar'/\idealm^k\B_\hbar'$, 
			\begin{equation}\label{modified C* action compatible with conical}
				t*_k(ab_k) = (t\cdot_c a)(t*_kb_k)
			\end{equation}
			where $\cdot_c$ denote the conical action. 
			
			\item We show that the $\C^*$-actions $*_k$ on $\B_\hbar'/\idealm^k\B_\hbar'$ induce a well defined $\C^*$-action on $\B_\hbar'$ that is compatible with the conical $\C^*$-action on $\sliceA_\hbar^{\wedge0}$, thus finishing the proof of the proposition. 
			
			Note that the natural maps  $p_k: \B_\hbar' / \idealm^k \B_\hbar' \to \B_\hbar' / \idealm^{k-1} \B_\hbar'$ intertwines $\Eu_m$, $\Eu_{c,k}$ and $\Eu_{c,k-1}$. Therefore, $p_k$ sends $(\B_\hbar'/\idealm^k\B_\hbar')^z_n $ 
			to $(\B_\hbar'/\idealm^{k-1}\B_\hbar')^z_n $.
			It follows that $p_k$ also sends the $n-f(z)$-eigenspace of $\Eu_{c,k}'$ to the $n-f(z)$-eigenspace of $\Eu_{c,k-1}'$ in $\B_\hbar' / \idealm^{k-1} \B_\hbar'$. In particular, for any $b_k\in \B_\hbar' / \idealm^k \B_\hbar'$, we have
			\begin{equation}\label{modified C* action compatible with pk}
				t*_{k-1} p_k(b_k) = p_k(t*_k b_k).
			\end{equation}
			
			Therefore, write any element of $\B_\hbar'$ as $b = (b_1,b_2,\cdots, b_k,\cdots)$, where $b_k\in \B_\hbar' / \idealm^{k} \B_\hbar'$ and $p_kb_k = b_{k-1}$. Define a $\C^*$-action on $\B_\hbar'$ by 
			\begin{equation}
				(t,b) \mapsto (t*_1b_1,\cdots,t*_kb_k,\cdots), \ t\in \C^*.
			\end{equation}
			By \eqref{modified C* action compatible with pk}, the image is a well defined element of $\B_\hbar'$. 
			By \eqref{modified C* action compatible with conical}, this $\C^*$-action is compatible with the conical $\C^*$-action on $\sliceA_\hbar^{\wedge0}$, as desired. 		
		\end{enumerate}
	\end{proof}
	Therefore, any object of $\CP (\underline{\A}_{ \hbar} ^ {\wedge 0})$ can be equipped with a $\C^*$-action that makes it an object of $^c \CP (\underline{\A}_{ \hbar} ^ {\wedge 0})$. 
	We get a functor $\F_8: \CP(\underline{\A}_{\hbar} ^ {\wedge 0}) \to {^c}\CP(\underline{\A}_{\hbar} ^ {\wedge 0})$, which is the identity on the level of objects if we forget the $\C^*$-actions. It is clear that $\F_8$ descends to functors $\CP^{S_n}(\underline{\A}_{\hbar} ^ {\wedge 0}) \to {^c}\CP^{S_n}(\underline{\A}_{\hbar} ^ {\wedge 0})$ and $\CP^{S_n}_0(\underline{\A}_{\hbar} ^ {\wedge 0}) \to {^c}\CP^{S_n}_0(\underline{\A}_{\hbar} ^ {\wedge 0})$, which we still denote by $\F_8$. 
	By a completely analogous argument, we can use the same function $f:\C\to \ZZ$ in the proof of \Cref{proposition modify c* action}
	to equip an object of $^c\CP (\underline{\A}_{ \hbar} ^ {\wedge 0})$ with a $\C^*$-action that makes it an object of $ \CP (\underline{\A}_{ \hbar} ^ {\wedge 0})$, and get functors $\F_8': ^c\CP (\underline{\A}_{\hbar} ^ {\wedge 0}) \to \CP (\underline{\A}_{\hbar} ^ {\wedge 0})$, $^c\CP^{S_n}(\underline{\A}_{\hbar} ^ {\wedge 0}) \to \CP^{S_n}(\underline{\A}_{\hbar} ^ {\wedge 0})$ and $^c\CP^{S_n}_0(\underline{\A}_{\hbar} ^ {\wedge 0}) \to \CP^{S_n}_0(\underline{\A}_{\hbar} ^ {\wedge 0})$. It is easy to see that $\F_8$ and $\F_8'$ are mutual quasiinverses. 
	
	\subsubsection{Taking finite part. }\label{functor take finite part}
	The algebra $\sliceA_\hbar$ is non-negatively graded with respect to the conical $\C^*$-action, and the ideal $\mathfrak{m}$ corresponding to $0$ is positively graded. Therefore, similar to the functor $\F_2,\F_2'$, 
	we get mutual quasi-inverse functors 
	\begin{equation}
		\F_9: ^c\CP^{S_n} (\underline{\A}_{\lambda,\hbar} ^ {\wedge 0}) \leftrightarrows ^c\CP^{S_n} (\underline{\A}_{\lambda,\hbar} ) :\F_9'
	\end{equation}
	where $\F_9$ is the functor of taking $\C^*$-locally finite part, and $\F_9'$ is the functor of taking $\mf{m}$-adic completion. They restrict to quasi-inverses between the subcategories $^c\CP^{S_n}_0 (\underline{\A}_{\lambda,\hbar} ^ {\wedge 0}) $ and $ ^c\CP^{S_n}_0 (\underline{\A}_{\lambda,\hbar} )$. 
	
	\subsubsection{Setting $\hbar = 1$.}\label{last functor invert hbar}
	For an object $\B'_\hbar$ of $\CP^{S_n}(\underline{\A}_{\lambda,\hbar})$, $ \B'_\hbar/(\hbar-1)\B'_\hbar$ is an object of $\HC^{S_n}(\underline{\A}_{\lambda})_F$. 
	Conversely, for any object $\B$ of $\HC^{S_n}(\underline{\A}_{\lambda})_F$, the Rees bimodule $\B_\hbar$ is an object of $\CP^{S_n}(\underline{\A}_{\lambda})$. 
	We get functors $\F_{10}: \CP^{S_n}(\underline{\A}_{\lambda}) \leftrightarrows \HC^{S_n}(\underline{\A}_{\lambda})_F : \F_{10}'$, similar to $\F_1$ and $\F_1'$. 
	
	They restrict to functors between the subcategories $\CP^{S_n}_0(\underline{\A}_{\lambda,\hbar}) \to \HC^{S_n}_{0}(\underline{\A}_{\lambda,\hbar})$ 
	\subsubsection{Compositions}
	We first compose the functors defined in \Cref{functor completion along L closure} through \Cref{functor take finite part}:
	\begin{equation}
		\begin{split}
			&\F_{\CP}: \CP_{\overline{\L}}(\A_\hbar) \to \CP^{S_n}_0(\underline{\A}_{\lambda,\hbar}) \\
			&\F_{\CP} = \F_9 \circ \F_8\circ\F_7\circ\F_6\circ\F_5\circ\F_4\circ\F_3\circ\F_2. \\
			&\F'_{\CP}: \CP^{S_n}_0(\underline{\A}_{\lambda,\hbar})\to  \CP_{\overline{\L}}(\A_\hbar) \\
			&\F'_{\CP} = \F'_2\circ\F'_3\circ\F'_4\circ\F'_5\circ\F'_6\circ\F'_7\circ\F'_8 \circ \F_9'.
		\end{split}
	\end{equation}
	All the objects involved are supported on suitable closed subsets (e.g. $\overline{\L}, \L, \L_{I\cup J}, 0$). Moreover, all the functors involved preserve $\C[\hbar]$-flatness. Therefore, by the results of \Cref{functor completion along L closure} through \Cref{functor take finite part}, the compositions $\F_{\CP}$ and $\F'_{\CP}$ are well defined quasi-inverses to each other.
	
	Let $\B$ be an object of $\HC_{\L}(\A)$; equip it with a good filtration and view it as an object of $\HC_{\L}(\A)_F$. 
	On the other hand, let $\underline{\mathcal{N}}$ be an object of $\HC_0^{S_n}(\sliceA)$; equip it with a good filtration and view it as an object of $\HC_0^{S_n}(\sliceA)_F$. 
	The following result is proved by an argument analogous to \cite[Section 3.4]{losev2011finite}. 
	\begin{proposition}\label{resulting hc bimodule independent of filtration}
		The $\sliceA$-bimodule $\F_{10}\circ \F_{\CP}\circ \F_1(\B)$ is independent of the choice of a good filtration on $\B$, and the $\A$-bimodule $\F_{1}'\circ\F'_{\CP}\circ\F_{10}'(\underline{\mathcal{N}})$ is independent of the choice of a good filtration on $\underline{\mathcal{N}}$. 
	\end{proposition}
	Thanks to \Cref{resulting hc bimodule independent of filtration} and the results in \Cref{subsubsection functor taking Rees module} through \Cref{last functor invert hbar}, we get the following identification.
	\begin{theorem}\label{bijection with fd equivariant slice modules}
		The functors $\F_{\CP}$ and $\F'_{\CP}$, together with taking Rees bimodule and setting $\hbar-1$, give mutually inverse bijections between objects of $\HC_{\L}(\A)$ and objects of $\HC^{S_n}_0(\underline{\A})$. 
	\end{theorem}

	\section{Modules over quantum slice algebras}\label{Hypertoric quantum slice algebras}
	In this subsection we realize the quantum slice quiver variety $\sliceA$ as a hypertoric enveloping algebra, defined in \cite{braden2012hypertoric}, and use the results there to classify the irreducible objects of $\HC^{S_n}_0(\underline{\A})$. 
	By \Cref{bijection with fd equivariant slice modules}, we get a classification of objects of $\HC_{\L}(\A)$. To emphasize the importance of the quantization parameter $\lambda$, we will write $\A_\lambda$ and $\sliceA_\lambda$ instead of $\A$ and $\sliceA$. 
	
	\subsection{Linear algebraic data}
	We first define necessary linear algebraic data associated with the algebra $\sliceA_\lambda$, following \cite[Section 2]{braden2012hypertoric}. 
	
	Recall we defined the vector space $\sliceR$ in \Cref{subsection the slice quiver} and labeled the vectors. 
	Let $\mathfrak{A} = \{v_{ij}^k,p_i^m|1 \le k \le 2\ell -2, 1\le m \le w \}$ be the set of all arrows of the slice quiver excluding the edge loops. 
	For each arrow $a\in \mathfrak{A}$, let $h_a$ be a formal variable. 
	Let $P = \C[h_a|a\in \mathfrak{A}]$ be the polynomial ring in $|\mf{A}|$ variables. Let 
	\begin{equation}\label{definition of quantum weight ring}
		\mathbf{P} = \C[h_a^+,h_a^-|a\in \mathfrak{A}]/(h_a^+ - h_a^- +1).
	\end{equation}
	There is a filtration on $\mathbf{P}$ given by $\deg h^+_a = \deg h^-_a = 1$ such that $\gr \mathbf{P} = P$. 
	Define
	\begin{equation}
		\begin{split}
			&W = \Spec P, \\ &\mathbf{W} = \Spec\mathbf{P}; \\
			&W_\ZZ = \{w\in W| h_a(w)\in \ZZ, \forall a\in\mathfrak{A}\};\\
			&\mathbf{W}_\ZZ = \{w\in \mathbf{W}| h^+_a(w)\in \ZZ, \forall a\in\mathfrak{A}\}.
		\end{split}
	\end{equation}
	So $W_\ZZ$ is a lattice in $W$. Let $\C[W_\ZZ]$ be its group ring, and let $T = \Spec\C[W_\ZZ]\cong (\C^*)^{\dim\sliceR}$. 
	Then $W_\ZZ$ is the character lattice of $T$, and $T$ act on $\sliceR$ with weights in $W_\ZZ$. 
	
	Let $D = D(\underline{R})$ denote the algebra of differential operators on $\underline{R}$. More explicitly, $D(\underline{R})$ is generated by $\{x_a, \partial_{a}|a\in\mathfrak{A}\}$ modulo the relations
	\begin{align}
		\begin{split}
			[x_a,x_b] &= 0 \\
			[\partial_a,\partial_b] &= 0	\\
			[x_a,\partial_b] &= -\delta_{ab}
		\end{split}
	\end{align}
	for all $a,b \in \mathfrak{A}$. Here, we understand $x_a$ as the linear function dual to $a$. For example, for $a = v_{ij}^k$, $x_a = x_{ij}^k$ defined in \Cref{subsection the slice algebra}. 
	
	The $T$-action on $\underline{R}$ induces a $T$-action on $D$. We decompose $D$ into $T$-weight spaces,
	\begin{equation}
		D = \bigoplus_{z\in W_\ZZ} D_z
	\end{equation}
	i.e. for $t\in T, d\in D_z$, $t.d = z(t) d$. In particular, we have a subalgebra
	\begin{equation}
		D_0 = \C[x_a\partial_a|a\in \mathfrak{A}] = \C[\partial_a x_a|a\in\mathfrak{A}]. 
	\end{equation}
	We identify $D_0 \cong \mathbf{P}$ via $x_a\partial_a \mapsto h_a^+, \partial_a x_a \mapsto h_a^-$. 
	
	The group $H\cong (\C^*)^n$ is a subgroup of $T$. Define 
	\begin{equation}
		\Lambda_0 = \{w\in W_\ZZ|w(t)=0 \text{ for all } t\in H\}
	\end{equation}
	which is a direct summand of the lattice $W_\ZZ$. 
	Define 
	\begin{equation}
		U = D^H = \bigoplus_{z\in \Lambda_0}D_z.
	\end{equation}
	Since $D_0\subset U$, we see $\mathbf{P}$ is also a subalgebra of $U$. 
	
	Recall that we computed $r(\lambda) = (\lambda,\lambda, \cdots , \lambda)\in \cartanh^*$ in \Cref{subsection morphism of Quantization parameters} and by definition, 
	\begin{equation}
		\underline{\A}_{\lambda} = U/U(\underline{\Phi} (\eta) - \<r(\lambda),\eta\>|\eta\in\cartanh).
	\end{equation}
	
	Let us write down $\underline{\Phi}(\eta)$ for $\eta = (\eta_1,\cdots,\eta_n) \in \cartanh$ as an element of $\mathbf{P}$. By the definition of the $H$-action on $\sliceR$, we have, for all $1 \le i<j \le n$, 
	\begin{align}
		\begin{split}
			& \eta.v_{ij}^k = (\eta_j-\eta_i)v_{ij}^k, 1\le k\le {\ell}-1;    \\
			& \eta.q_{i}^m = -\eta_i q_{i}^m, 1\le m \le w.
		\end{split}
	\end{align}
	Therefore, 
	\begin{equation}\label{formula for slicePhi eta}
		\underline{\Phi}(\eta) = 
		\sum_{1\le i\neq j\le n} 
		\sum_{k=1}^{{\ell}-1} (\eta_j-\eta_i)h^+_{v_{ij}^k} 
		-\sum_{i=1}^n\sum_{k=1}^{w} \eta_i h^+_{q_i^k}.
	\end{equation}
	
	\subsection{Finite dimensional simple modules}
	In this subsection we recall the classification of finite dimensional simple $\underline{\A}_{\lambda}$-modules in terms of bounded chambers of certain hyperplane arrangements, which is a special case of the more general classification of the hypertoric category $\sheafO$ in \cite{braden2012hypertoric}. 
	Then, we give some conditions on $\lambda$ for $\underline{\A}_{\lambda}$ to have a finite dimensional simple module. 
	
	Let $M$ be a simple finite dimensional $U$-module. Let $v\in \mathbf{W} = \Spec\mathbf{P}$ be a closed point, and let $\mathfrak{m}_v\subset \mathbf{P}$ be the corresponding maximal ideal. Define the \textit{$v$-weight space} of $M$ to be 
	\begin{equation}
		M_v = \{m\in M| \mathfrak{m}_v^km=0 \text{ for large enough }k \}.
	\end{equation}
	Since $M$ is finite dimensional, we have a direct sum decomposition 
	\begin{equation}
		M = \bigoplus_{v\in \mathbf{W}} M_v.
	\end{equation}
	Note that $\Lambda_0$ acts freely by translation on the affine space $\mathbf{W}$. 
	It is easy to verity that $D_z M_v \subset M_{z+v}$ for any $ z\in \Lambda_0 $ (so that $D_{z}\subset U)$. 
	Since $M$ is simple, the set of weights $\{v\in \mathbf{W} | M_v\neq 0\}$ must be a single $\Lambda_0$-orbit in $\mathbf{W}$. Denote this orbit by $\mathbf{\Lambda}$, so $	M = \bigoplus_{v\in \mathbf{\Lambda}} M_v$, with all but finitely many $M_v$ being 0. 
	
	The $U$-module $M$ descends to an $\underline{\A}_{\lambda}$-module if and only if $\underline{\Phi}(\eta) - \<r(\lambda),\eta\>$ acts on $M$ by 0 for all $\eta\in \cartanh$. The equations $\underline{\Phi}(\eta) - \<r(\lambda),\eta\>=0$ for all $\eta\in\cartanh$ cut out an affine subspace $\mathbf{V}$ of $\mathbf{W}$, and $M$ carries an induced (simple) $\underline{\A}_{\lambda}$-module structure if and only if $\mathbf{\Lambda}\subset \mathbf{V}$. 
	
	Fix a $\Lambda_0$-orbit $\mathbf{\Lambda}\subset \mathbf{V}$. Define the following subset of $\mathfrak{A}$: 
	\begin{equation}
		I_\mathbf{\Lambda} = \{a\in \mathfrak{A} | h^+_a(\mathbf{\Lambda}) \subset \ZZ \}. 
	\end{equation}
	Pick any $v\in \mathbf{\Lambda}$.
	For any sign vector $\alpha\in \{+.-\}^{I_\mathbf{\Lambda}}$, define the simple $D$-module
	\begin{equation}
		L_\alpha : = 
		D / D(\partial_a|\alpha({a}) = +) + D(a|\alpha({a}) = -) + D(h_a^+ - h_a^+(v)| a\not\in I_\mathbf{\Lambda}).
	\end{equation} 
	The module $L_\alpha$ is independent of the choice of $v$. 
	Define 
	\begin{equation}
		L_\alpha^{\mathbf{\Lambda}} = \left\{ m\in \bigoplus_{v\in \mathbf{\Lambda}} (L_\alpha)_v  | (\underline{\Phi}(\eta) - \<r(\lambda),\eta\>) m =0   \right\} 
	\end{equation}
	so that $L_\alpha^{\mathbf{\Lambda}}$ is an $\underline{\A}_{\lambda}$-module. 
	
	Fix $\alpha\in \{+.-\}^{I_\mathbf{\Lambda}}$. For $a\in\mathfrak{A}$, define the half spaces of $\mathbf{V}$: \begin{equation}
		\begin{split}
			R(a,+) = \{v\in \mathbf{V}|h^+_a(v)\ge 0\},\\
			R(a,-) = \{v\in \mathbf{V}|h_a^-(v)\le 0\}.
		\end{split}
	\end{equation}
	In view of the relation $ h_a^- = h_a^++1$, we can also write
	\begin{equation}
		R(a,-) = \{v\in \mathbf{V}|h_a^+(v)\le -1\}.
	\end{equation}
	These half spaces cut out the chamber 
	\begin{equation}
		R(\alpha): = \bigcap_{a\in I_\mathbf{\Lambda }} R(a,\alpha(a))\subset \mathbf{V}.
	\end{equation}
	which can be bounded, unbounded or empty. 
	\begin{proposition}[{\cite[Proposition 3.5]{braden2012hypertoric}}]\label{BLPW classification of fd simple modules}
		The set of isomorphism classes of finite dimensional simple $\underline{\A}_{\lambda}$-modules are precisely $L_\alpha^{\mathbf{\Lambda}}$ where $\alpha\in\{+,-\}^{I_\mathbf{\Lambda }}$ is such that $R(\alpha)$ is bounded and $R(\alpha) \cap \mathbf{\Lambda}\neq \emptyset$. 
	\end{proposition}
	
	\begin{remarklabeled}\label{chambers contain eigenvalues of ad hij}
		From the construction, it is easy to see that the elements of $R(\alpha) \cap \mathbf{\Lambda}$ are the eigenvalues of $ h_a^+= x_a\partial_a$ in $L_\alpha^{\mathbf{\Lambda}}$.
	\end{remarklabeled}
	
	\subsection{Finding bounded chambers}
	Thanks to \Cref{BLPW classification of fd simple modules}, our goal now is to find bounded chambers in $\mathbf{V}$ of the form $R(\alpha)$.
	We first reduce this task to a similar one in a smaller affine space, then use some combinatorial arguments to finish the classification. 
	
	For each pair $1\le i<j\le n$, define
	\begin{equation}h_{ij} =\sum_{k=1}^{{\ell}-1} \left( h^+_{v_{ij}^k} - h^+_{v_{ji}^k} \right).
	\end{equation}
	In particular, $h_{ij} = - h_{ji}$. 
	For each $1\le i\le n$, define 
	\begin{equation}q_i =\sum_{k=1}^{w} h^+_{q_k^i} .\end{equation}
	The defining equation  $\underline{\Phi}(\eta) - \<r(\lambda),\eta\>=0$ of $\mathbf{V}$ is equivalent to the condition that the coefficient of each $\eta_i$ in \eqref{formula for slicePhi eta} equals $\lambda$. 
	That is,  
	\begin{equation}\label{more simplified equation for VV} 
		-\sum_{j=1}^{n} h_{ij}-q_i = \lambda , \ \text{ for any } 1\le i\le n.
	\end{equation}
	
	View $h_{ij}, q_i$ as independent variables (we assume $i < j$ for $h_{ij}$) and define $\WW' = \Spec\C[h_{ij},q_i]$. This is an affine space and $\dim \WW' = \binom{n}{2}+n$. Let $\VV'$ be the affine subspace of $\WW'$ cut out by the equations \eqref{more simplified equation for VV}. By construction, there is a natural affine linear epimorphism $p: \VV \twoheadrightarrow \VV'$. Let $\mathbf{\Lambda}' = p(\mathbf{\Lambda})$. 
	Let $I_{\mathbf{\Lambda}'}$ be the subset of $\{h_{ij},q_i\}$ that takes integral value on $\mathbf{\Lambda}'$. 
	Fix a sign vector $\alpha'\in \{+,-\}^{I_{\mathbf{\Lambda}' }}$. For $h_{ij} \in I_{\mathbf{\Lambda}' }$, define half spaces 
	\begin{equation}
		\begin{split}
			R'(h_{ij},+) &= \{v\in \VV'| h_{ij}(v) \ge {\ell}-1 \},\\
			R'(h_{ij},-) &= \{v\in \mathbf{V}'| h_{ij}(v) \le -\ell+1 \}.
		\end{split}
	\end{equation} 
	For $q_i\in I_{\mathbf{\Lambda}'} $, define the half spaces 
	\begin{equation}
		\begin{split}
			R' (q_i,+) &= \{v\in \mathbf{V}'| q_i(v) \ge 0 \},\\
			R' (q_i,-)& = \{v\in \mathbf{V}'| q_i(v) \le -w_i \}.
		\end{split}
	\end{equation}
	Define the corresponding chamber of intersections 
	\begin{equation}
		R'(\alpha') = \bigcap_{h_{ij} \in I_{\mathbf{\Lambda}' } } R'(h_{ij},\alpha'(h_{ij})) \cap \bigcap_{q_i\in I_{\mathbf{\Lambda}'}} R'(q_i, \alpha'(q_i)) \subset \VV'.
	\end{equation}
	The next proposition relates the bounded chambers $R(\alpha)$ and $R'(\alpha')$. 
	\begin{proposition}\label{proposition bijection of bounded chambers}
		The following are true. 
		\begin{enumerate}
			\item If $R(\alpha)$ is nonempty and bounded for some $\alpha\in  \{\pm\}^{I_\mathbf{\Lambda }}$, then $h_{ij}\in I_{\mathbf{\Lambda }'}$ for all pairs $1\le i\neq j \le n$. 
			If we further assume $w\ge 2$, then $q_i\in I_{\mathbf{\Lambda }'}$ for all $i$. 
			\item There is a bijection between
			\begin{equation}\label{set of bounded chambers before grouping variables}
				\{\alpha\in  \{\pm\}^{I_\mathbf{\Lambda }} | R(\alpha) \text{ is bounded and nonempty}\},
			\end{equation}
			and 
			\begin{equation}\label{set of bounded chambers after grouping variables}
				\{\alpha' \in  \{\pm\}^{I_{\mathbf{\Lambda }'}} | R'(\alpha')\text{ is bounded and nonempty}\}.
			\end{equation}
		\end{enumerate}
	\end{proposition}
	\begin{proof}
		We prove (1). Assume $\alpha\in  \{\pm\}^{I_\mathbf{\Lambda }}$ is such that $ R(\alpha)\neq \emptyset$ is bounded. Let $z \in R(\alpha)$. If $h_{ij}\not\in I_{\mathbf{\Lambda }'}$, then $h_{ij}(z)\not\in\ZZ$, and therefore some $h^+_{v^k_{ij}}(z)\not\in\ZZ$, i.e. $h^+_{v^k_{ij}} \not\in I_{\mathbf{\Lambda }}$. Without loss of generality, let $h^+_{v^1_{ij}} \not\in I_{\mathbf{\Lambda }}$ and $h^+_{v^2_{ij}}(v) \ge 0 $. Then we can define a ray $\gamma(t)\in \mathbf{V}$ by coordinates 
		\begin{equation}
			h^+_{v^1_{ij}}(\gamma(t)) = h^+_{v^1_{ij}}(z)-t; \ \  h^+_{v^{2}_{ij}}(\gamma(t)) = h^+_{v_{ij}^2}(z)+t, t\in \RR^+
		\end{equation}
		and all other coordinates of $\gamma(t)$ are the same as those of $z$. This defines an unbounded subset of $R(\alpha)$, contradiction. A similar construction shows $q_i\in I_{\mathbf{\Lambda }'}$ for all $i$. This proves (1).  
		
		We prove (2) in several steps. We first construct a map from \eqref{set of bounded chambers before grouping variables} to \eqref{set of bounded chambers after grouping variables}, so assume the former set is nonempty. 
		
		\textit{Step 1.} We show that 
		\begin{enumerate}[label=\alph*)]
			\item if $h_{ij}\in I_{\mathbf{\Lambda }'}$, then for all $1\le k\le {\ell}-1$, $h^+_{x_k^{ij}}\in I_{\mathbf{\Lambda }}$;
			\item if $q_i \in I_{\mathbf{\Lambda }'}$, then for all $1\le k\le {\ell}-1$, $h^+_{x_k^{ij}}\in I_{\mathbf{\Lambda }}$
		\end{enumerate} 
		The idea is similar to the proof of (1).
		Assume the contrary. Let $z\in R(\alpha)$. Since $h_{ij} (z) \subset \ZZ$, we may assume without loss of generality that there exists a pair, say 
		$h^+_{v^1_{ij}},  h^+_{v^2_{ji}}$, 
		such that $h^+_{v^1_{ij}}(z),  h^+_{v^2_{ji}}(z)\not\in \ZZ$. 
		Then we define the ray $\gamma(t)\in \mathbf{V}$ through $v$ by 
		\begin{equation}h^+_{v^1_{ij}}(\gamma(t)) = h^+_{v^1_{ij}}(z)+t, h^+_{v^2_{ji}}(\gamma(t)) = h^+_{v^2_{ji}}(z)+t, t\in \RR^+,\end{equation} 
		while all other coordinates are the same as those of $v$. This is an unbounded line in $R(\alpha)$ whose image under $p$ is a single point, contradiction. This proves a). Part b) is proved similarly. 
		Note that we do not need to assume $w \ge 2$; in the case $w=1$ we directly have $q_i = h^+_{q_i^1}$. 
		
		\textit{Step 2. }
		Let $z\in R(\alpha)$. We show that, for any $1\le i\neq j\le n$, 
		\begin{enumerate}[label = \alph*)]
			\item if $h_{ij}(z) \ge 0$, then $h^+_{v^k_{ij}}(z) \ge 0$ and $h^+_{v^k_{ij}}(z) \le -1$ for all $1\le k \le {\ell}-1$.
			In particular, $h_{ij}(z) \ge {\ell}-1$; 
			\item if $h_{ij}(z) \le -1$ then $h^+_{v^k_{ij}}(z) \le -1$ and $h^+_{v^k_{ji}}(z)\ge 0$ 
			for all $1\le k \le {\ell}-1$. In particular, $h_{ij}(z) \le - l+1$.
		\end{enumerate}
		We prove a); b) is proved completely analogously. If $h_{ij}(z) \ge 0$ then by definition, $h^+_{v^k_{ij}}(z) \ge 0$ for {some} $1\le k \le {\ell}-1$ or $h^+_{v^k_{ji}}(z) \le -1$ for {some} $1\le k \le {\ell}-1$. 
		Without loss of generality, let $h^+_{v^1_{ij}}(z) \ge 0$. 
		Then, if $h^+_{v^m_{ij}} \le -1$ for some $1\le m\le {\ell}-1$,  
		we can define a ray $\gamma(t)\subset \mathbf{V}$ by 
		\begin{equation}h^+_{v^1_{ij}}(\gamma(t)) = h^+_{v^1_{ij}}(z )+t,  h^+_{x_{k_1}^{ij}}(\gamma(t)) = h^+_{x_{k_1}^{ij}}(v )-t, t\in \RR^+\end{equation} 
		while all other coordinates are the same as those of $v$. This is an unbounded ray in $R(\alpha)$, contradiction. 
		A similar contradiction arises if $h^+_{v_{ij}^{k}}(z ) \ge 0$ for some $1\le k \le {\ell}-1$. This proves statement a) of Step 2. 
		Statement b) is proved by by a completely analogous argument. Moreover, we see if $q_i\in I_{\mathbf{\Lambda }'} $ then
		\begin{enumerate}[label = \alph*)]\setcounter{enumi}{2}
			\item If $q_i(z )\ge 0$ then $h^+_{q^k_i}(z)\ge 0$ for all $1\le k\le w$. 
			\item If $q_i(z )\le -1$ then $h^+_{q^k_i}(z)\le -1$ for all $1\le k\le w$. In particular $q_i(v)\le -w$. 
		\end{enumerate}
		
		\textit{Step 3. } Let $\alpha $ be such that $R(\alpha)$ is bounded nonempty. We construct a sign vector $\alpha'\in  \{\pm\}^{I_\mathbf{\Lambda '}}$ from $\alpha$. 
		Fix $z\in R(\alpha)$. By the previous steps, we have 
		\begin{enumerate}
			\item $h^+_{v^1_{ij}}(z) \ge 0 \iff h_{ij}(z)\ge {\ell}-1$;
			\item $h^+_{v^1_{ij}}(v) \le -1 \iff h_{ij}(z)\le  -\ell+1$;
			\item $h^+_{q^k_i}(z)\ge 0 \iff q_i(z)\ge 0$;
			\item $h^+_{q^k_i}(z)\le -1 \iff q_i(z)\le -w$;
		\end{enumerate}
		therefore, for $h_{ij}\in I_{\mathbf{\Lambda }'}$, set $\alpha(h_{ij}) = +$ if $h^+_{v^1_{ij}}(z)$ and $-$ if $h^+_{v^1_{ij}}(z) \le -1$; for $q_i \in I_{\mathbf{\Lambda }'} $, set $\alpha(q_i) = +$ if $h^+_{q^k_i}(z)\ge 0$ and $-$ if $h^+_{q^k_i}(z)\le -1$. Then  $R'(\alpha')$ is the image of $R(\alpha)$ under the natural projection $p$, so $R'(\alpha')$ is compact and hence bounded. 
		
		\textit{Step 4.}
		We construct the inverse map. Suppose $\alpha'$ is such that $R'(\alpha')$ is nonempty and bounded. 
		Let $z'\in R'(\alpha')$.
		Pick $z \in \mathbf{V}$ whose coordinates $h^+_{v^k_{ij}}(z)$ and $h^+_{q^k_i}(v)$ satisfy 
		\begin{enumerate}[label = \alph*)]
			\item if $\alpha(h_{ij})= +$, i.e. $h_{ij}(z') \ge {\ell}-1$, then $h^+_{v^k_{ji}}(z) \ge 0$ and $h^+_{v^k_{ij}}(z) \le -1 $ for all $1\le k \le {\ell}-1$;
			\item if $\alpha(h_{ij})= -$, i.e. $h_{ij}(z') \le -{\ell}-1$, then $h^+_{v^k_{ij}}(z) \le -1$ and $h^+_{v^k_{ji}}(z) \ge 0 $ for all $1\le k \le {\ell}-1$;
			\item if $\alpha(q_i) = +$, i.e. $q_i(z')\ge 0$ then $h^+_{q_i^k}(z)\ge 0$ for all $1\le k\le w$;
			\item if $\alpha(q_i) = -$, i.e. $q_i(z')\le -w$ then $h^+_{q_i^k}(z)\le -1$ for all $1\le k\le w$. 
		\end{enumerate}
		Such element $z$ exists by the definition of $R'(\alpha')$, and they define a unique chamber 
		$R(\alpha)$, which defines the sign vector $\alpha$. 
		Note that $p(R(\alpha)) = R'(\alpha')$. 
		It is straightforward to show $R(\alpha)$ is bounded and the maps defined in Step 3 and 4 are mutual inverse. The proof is now complete. 
	\end{proof}
	
	Therefore, it suffices to find bounded nonempty chambers in $\VV'$ for different  ${\mathbf{\Lambda}'}$. We first consider the existence problem. 
	
	\begin{proposition}\label{propositon existence of fd for min leaf}
		If $\underline{\A}_{\lambda}$ has a finite dimensional simple module, then $\lambda\in\ZZ$. If $\lambda\in\ZZ$ and
		\begin{equation}
			\lambda \ge ({\ell}-1)(n-1)+wn,
		\end{equation} 
		then $\underline{\A}_{\lambda}$ has a finite dimensional simple module. 
	\end{proposition}
	\begin{proof}
		Assume $\underline{\A}_{\lambda}$ has a finite dimensional simple module. Then by \Cref{proposition bijection of bounded chambers}, there is some $\mathbf{\Lambda}$ such that  $I_{\mathbf{\Lambda }'}$ contains all $h_{ij}$, 
		and a sign vector $\alpha'$ such that $R'(\alpha')\subset \mathbf{V}'$ is nonempty and bounded. 
		Let $z' \in R'(\alpha')$. 
		If $\lambda\not\in \ZZ$, then $w=1$ and $q_i(z')$ $\not\in\ZZ$ for all $i$. 
		Without loss of generality, let $h_{12}(z)\ge {\ell}-1$. Define a ray $\gamma(t),t\in \RR^+$ in the chamber $R'(\alpha')$ by 
		\begin{equation}h_{12}(\gamma(t)) = h_{12}(z') +t; q_1(\gamma(t)) = q_1(z')+t; q_2(\gamma(t)) = q_2(z')-t
		\end{equation}
		and all other coordinates being the same as those of $v$. This is an unbounded subset of $R'(\alpha')$, contradiction. Therefore, $\lambda$ must be integral. 
		
		Conversely, suppose $\lambda\in\ZZ$ and $\lambda \ge ({\ell}-1)(n-1)+wn$. Then it is easy to check that $R'(\alpha')$, where $\alpha'(q_i) = -, \alpha'(h_{ij}) = -$ for $1\le i<j\le n$, is a nonempty bounded chamber. 
	\end{proof}
	
	We are going to describe all the sign vectors that give bounded chambers, for different $\lambda$'s. First, the following lemma imposes restriction on the sign of $q_i$. 
	\begin{lemma}
		Let $\lambda\in\ZZ$. 
		Suppose $h_{ij},q_i \in I_{\mathbf{\Lambda }'}$ for all $1\le i<j\le n$, and $R'(\alpha')\neq \emptyset$ is bounded. Then 
		\begin{enumerate}
			\item $\alpha'(q_i) = +$ for all $i$ if $\lambda \le 0$;
			\item $\alpha'(q_i) = -$ for all $i$ if $\lambda\ge 1$.
		\end{enumerate}
	\end{lemma}
	\begin{proof}
		Assume $\lambda\ge 1$ and let $I\subset \{1,2....,n\}$ be the set of $i$ such that $\alpha'(q_i) = +$. Assume, for a contradiction, that $I$ is nonempty. Let $J$ be the complement of $I$. Note that the sum of \eqref{more simplified equation for VV} through $1\le i\le n$ implies $\sum q_i = -n\lambda <0 $. Therefore $J$ is nonempty. 
		
		Take the sum of \eqref{more simplified equation for VV} through $i\in I$, we see
		\begin{equation}\label{equation in lemma for all q same sign}
			-\sum_{k\in I,k'\in J}h_{kk'} -\sum_{i\in I} q_i= |I|\lambda.
		\end{equation}
		In particular, 
		\begin{equation}\label{contradictory inequality proving all alpha qi same sign}
			\sum_{k\in I,k'\in J}h_{kk'}<0.
		\end{equation}
		If $\alpha'(q_k)=+,\alpha'(q_{k'})=-$ and $R'(\alpha')$ is nonempty and bounded, then $\alpha'(h_{kk'})=+$; otherwise we can increase $q_k$ and decrease $q_{k'},h_{kk'}$ by same amount while staying in $R'(\alpha')$, and $R'(\alpha')$ cannot be bounded. 
		But this contradicts \eqref{contradictory inequality proving all alpha qi same sign}. Therefore, $I$ is empty and we conclude $\alpha(q_i) = -$ for all $i$. 
		The case for $\lambda\le 0$ is proved completely analogously. 
	\end{proof}
	
	The following notations will be convenient. 
	\begin{definition}
		Suppose all $h_{ij} \in I_{\mathbf{\Lambda}'}$. 
		We define 
		\begin{equation}B(\mathbf{\Lambda}') := \{\alpha \in  \{\pm\}^{I_{\mathbf{\Lambda }'}} | R'(\alpha)\neq \emptyset\text{ is bounded}\}.\end{equation} 
		For each $1\le i\le n$, define 
		\begin{equation}B_i({\mathbf{\Lambda }'}) = \{\alpha' \in  \{\pm\}^{I_{\mathbf{\Lambda }'}} | R'(\alpha')\neq \emptyset\text{ is bounded},\alpha(h_{ij}) = +\text{ for all } j\neq i\},
		\end{equation}
		where, as before, we set $h_{ij} = -h_{ji}$, and $\alpha'(h_{ij}) = + $ if and only if $\alpha'(h_{ji}) = -$. Clearly, each $B_i({\mathbf{\Lambda }'})$ is a subset of $B(\mathbf{\Lambda}')$. 
	\end{definition}
	\begin{lemma}
		If $\lambda\in \ZZ$ and $B(\mathbf{\Lambda}') \neq \emptyset$, then $I_{\mathbf{\Lambda}'} = \{h_{ij}, q_i|1\le i<j\le n\}$, i.e. the full set of variables. 
	\end{lemma}
	\begin{proof}
		By \Cref{proposition bijection of bounded chambers}, all $h_{ij}\subset I_{\mathbf{\Lambda}'}$. Suppose $q_i \not\in I_{\mathbf{\Lambda}'}$. Then, take the sum of \eqref{more simplified equation for VV} over all $i$ we see that $q_j \not\in I_{\mathbf{\Lambda}'}$ for some $j\neq i$ as well. Then we can produce a ray involving $h_{ij}, q_i, q_j$ in the chamber $R'(\alpha')$, contradiction. 
	\end{proof}
	\begin{theorem}\label{proposition counting integral parameter simples}
		The following are true. 
		\begin{enumerate}
			\item If $\lambda\not\in\ZZ$, then $|B(\mathbf{\Lambda}')| = 0$. 
			\item If $1\le \lambda < (n-1)({\ell}-1) + w$ and $\lambda\in\ZZ$, then $|B(\mathbf{\Lambda}')| = 0$. 
			\item If $-(n-1)({\ell}-1) <\lambda \le 0$ and $\lambda\in\ZZ$, then $|B(\mathbf{\Lambda}')| = 0$. 
			\item If $\lambda\ge (n-1)({\ell}-1) + w$ and $\lambda\in \ZZ$, then  $|B(\mathbf{\Lambda}')| = n!$. 
			\item If $\lambda\le - (n-1)({\ell}-1)$ and $\lambda\in\ZZ$, then  $|B(\mathbf{\Lambda}')| = n!$. 
		\end{enumerate}
	\end{theorem}
	\begin{proof}
		(1) has been proved in \Cref{propositon existence of fd for min leaf}. 
		
		We prove (2). The proof of (3) is completely analogous. 
		First, we claim that $B(\mathbf{\Lambda}') = \bigsqcup_{1\le i\le n} B_i(\mathbf{\Lambda}')$. The sets $B_i(\mathbf{\Lambda}')$ are disjoint for different $i$ by the definition. Suppose $\alpha'\in B(\mathbf{\Lambda}')$ but $\alpha'\not\in B_i(\mathbf{\Lambda}')$ for any $i$. 
		Then for each $i$ there exists some $i'$ such that $\alpha'(h_{ii'}) = -$, i.e. for any $v\in R'(\alpha')$, $h_{ii'} (v) \le -({\ell}-1)$. We can thus find a sequence $i_1,i_2,...,i_k$ without repetition such that for any $v \in R'(\alpha')$, the following cycle of inequalities hold:
		\begin{equation}
			h_{i_1,i_2}(v)\le -\ell+1, 
			h_{i_2,i_3}(v)\le -\ell+1,...,h_{i_k,i_1}(v) \le -\ell+1.\end{equation}
		Therefore, we can define a ray starting at some $v$ in $R'(\alpha')$ by increasing all $h_{i_j,i_{j+1}}$ coordinates by the same amount $t>0$. This contradicts the boundedness of $R'(\alpha')$. Therefore, $B(\mathbf{\Lambda}') = \bigsqcup_{1\le i\le n} B_i(\mathbf{\Lambda}')$. 
		
		Suppose $B(\mathbf{\Lambda}')$ is nonempty; without loss of generality, let $B_1(\mathbf{\Lambda}')$ be nonempty. Then the equation \eqref{more simplified equation for VV} for $i=1$ shows $\lambda\ge (n-1)({\ell}-1)+w$, and (2)  is proved. 
		
		We prove (4) by induction on $n$. The proof of (5) is completely analogous. 	
		When $n=2$, we can check directly that $B(\mathbf{\Lambda}') = \{\alpha'_1,\alpha'_2\}$, where 
		\begin{equation}
			\begin{split}
				\alpha'_1(h_{12}) = +, \alpha'_1(q_1) = \alpha'_1(q_2) = - ,\\
				\alpha'_2(h_{12}) = -, \alpha'_2(q_1) = \alpha'_2(q_2) = - .
			\end{split}
		\end{equation}
		
		Assume (4) holds for $n-1$. 
		We claim that $|B_i(\mathbf{\Lambda}')| = (n-1)!$ for all $1\le i\le n$. Together with the decomposition 
		$B(\mathbf{\Lambda}') = \bigsqcup_{1\le i\le n} B_i(\mathbf{\Lambda}')$, we get (4).	
		
		Without loss of generality, let us show $|B_n(\mathbf{\Lambda}')| = (n-1)!$. 
		The first $n-1$ equations of \ref{more simplified equation for VV} can be written as 
		\begin{equation}\label{modified more simple equation for VV}
			-\sum_{j=1}^{n-1} h_{ij}-q_i = \lambda+h_{in}, 1\le i \le n-1.
		\end{equation}
		A sign vector $\alpha' \in \{\pm\}^{\mathbf{\Lambda}'}$ lies in $|B(\mathbf{\Lambda}')_n|$ if and only if $\alpha'(h_{in}) = -$ for all $i<n$, $\alpha'(q_i) = - $ for all $i$, and the chamber $R'(\alpha') \subset \VV'$ is nonempty and bounded. 
		If $\alpha'$ is such a vector, then there is an element $z' \in R'(\alpha')$ such that $h_{in}(z') = -\ell+1$ for all $i<n$. 
		Therefore, consider the following system of equations obtained from \eqref{modified more simple equation for VV}: 
		\begin{equation}\label{equation for VV''}
			-\sum_{j=1}^{n-1} h_{ij}-q_i = \lambda-\ell+1,\ \  1\le i \le n-1.
		\end{equation}
		When $\lambda \ge (n-1)({\ell}-1)+w$, the right hand side of \eqref{equation for VV''} is $\ge (n-2)({\ell}-1) + w$. By induction hypothesis, there are exactly $(n-1)!$ sign vectors $\alpha'' \in \{\pm\}^{\{h_{ij}, q_1,\cdots, q_{n-1} | 1\le 1<j\le n-1\} }$ such that the chamber of vectors $z\in \VV'$ satisfying 
		\begin{itemize}
			\item \eqref{equation for VV''};
			\item $h_{in}(v) = -\ell+1$; 
			\item $h_{ij}(v) \ge {\ell}-1$ if $\alpha''(h_{ij}) = +$; 
			\item $h_{ij}(v) \le -\ell+1$ if $\alpha''(h_{ij}) = -$; 
			\item $q_i(v) \le 0$ if $\alpha''(q_i) = -$;
			\item $q_i(v) \ge w$ if $\alpha''(q_i) = +$
		\end{itemize}
		is nonempty and bounded. Moreover, all such $\alpha''$ satisfies $\alpha''(q_i) = -$. 
		Denote the set of such dimension vectors by $B''$.
		Let us show there is a bijection between $B''$ and $B_n(\mathbf{\Lambda}')$. 
		
		Suppose $\alpha'\in B_n(\mathbf{\Lambda}')$. We can remove its $h_{ij}$-components to get an $\alpha''\in B''$. 
		Conversely, suppose $\alpha''\in B''$. Let $\alpha' \in \{\pm \}^{I}$ be defined by extending $\alpha''$ such that $\alpha'(h_{in}) = -$ for $1\le i \le n-1$ and $\alpha'(q_n) = -$. We show that $\alpha' \in B_n(\mathbf{\Lambda}')$. 
		
		In fact, it is clear that $R'(\alpha')$ is nonempty. To see it is bounded, take the sum of \eqref{more simplified equation for VV} for $1\le i\le n$, we see $\sum_{i=1}^n q_i = n\lambda$. Since all $\alpha'(q_i)=-$, we see $q_i(R'(\alpha'))$ is bounded for all $i$. Thanks to the boundedness of $q_n$, \eqref{more simplified equation for VV} for $i=n$ shows $h_{in}(R'(\alpha'))$ is bounded as well. 
		Finally, consider the linear map $R'(\alpha') \to \C^{n}, v\mapsto (h_{1n}(v),\cdots,h_{n-1,n}(v),q_n(v))$. Its image is closed and bounded, and by \eqref{modified more simple equation for VV}, the fiber over each point is bounded (or empty). 
		Therefore, $R'(\alpha')$ is bounded.
		This establishes the bijection between $B''$ and $B_n(\mathbf{\Lambda}')$. 
		We have now proved $|B_n(\mathbf{\Lambda}')| = (n-1)!$, and hence $|B(\mathbf{\Lambda}')| = n!$. 
	\end{proof}
	
	\begin{remarklabeled}\label{bijection bounded chambers and s-ordering}
		There is a natural bijection between $B(\mathbf{\Lambda}')$ and the set of orderings of $\{1,2,...,n\}$. Namely, In the proof of \Cref{proposition counting integral parameter simples}, we see each $\alpha'\in B(\mathbf{\Lambda}')$ lies in a unique $B_{i_1}(\mathbf{\Lambda}')$. 
		The set $B_{i_1}(\mathbf{\Lambda}')$ is in bijection with $B(\mathbf{\Lambda}_{\neq i_1}')$ where $\mathbf{\Lambda}_{\neq i_1}' = \{h_{ij}, q_k|1\le i<j\le n, 1\le k\le n, \text{ and } i,j,k\neq i_1\}$. 
		Truncate $\alpha'$ and view it as an element of $ B(\mathbf{\Lambda}_{\neq i_1}')$, there is a unique $i_2\neq i_1$ such that $\alpha'\in B_{i_2}(\mathbf{\Lambda}_{\neq i_1 }')$. Repeat this process, we get a unique ordering $(i_1,i_2,...,i_n)$ on the set $\{1,2,...,n\}$. 
	\end{remarklabeled}
	
	\begin{example}
		Suppose $n=4$ and $\lambda\in \ZZ, \lambda\ge (n-1)({\ell}-1) + w$. Then the element $\alpha' \in B(\mathbf{\Lambda}')$ corresponding to $(3,1,2,4)$ is given by
		\begin{equation}
			\begin{array}{ccc}
				\alpha'(h_{13}) = - & \alpha'(h_{23}) = - & \alpha'(h_{34}) = +\\
				\alpha'(h_{12}) = + & \alpha'(h_{14}) = + & \                  \\
				\alpha'(h_{24}) = + & \                  & \
		\end{array}\end{equation}
		while all $\alpha'(q_i) = -$. 
	\end{example}
	Recall $S_n$ acts on $\underline{R}$ by permuting the vertices of the slice quiver, i.e. by permuting $1\le i\le n$. 
	We can now combine \Cref{chambers contain eigenvalues of ad hij} and \Cref{bijection bounded chambers and s-ordering} 
	to get the following corollary. 
	\begin{corollary}\label{Sn transitive free action on f.d. simple left}
		When $\lambda\in\ZZ$ and $\lambda \ge (n-1)({\ell}-1)+w$, 
		or $\lambda\le -(n-1)({\ell}-1)$,  the group $S_n$ acts freely transitively on the set of finite dimensional simple left $\sliceA_\lambda$-modules. 
	\end{corollary}
	
	By completely analogous arguments, we see that the set of finite dimensional simple right $\sliceA_\lambda$-modules is also in bijection with $S_n$, and is equipped with a transitive $S_n$-action. 
	Therefore, the set of finite dimensional simple $\sliceA_\lambda$-bimodules is in bijection with $S_n \times S_n$. 
	We have the following bimodule version of \Cref{Sn transitive free action on f.d. simple left}. 
	\begin{corollary}
		When $\lambda\in\ZZ$ and $\lambda \ge (n-1)({\ell}-1)+w$ 
		or $\lambda\le -(n-1)({\ell}-1)$,  
		there are precisely $n!$ non-isomorphic finite dimensional $S_n$-equivariant simple bimodules over $\A_{\lambda}$. When $-(n-1)({\ell}-1) <\lambda < (n-1)({\ell}-1)+w$, $\A_{\lambda}$ has no simple finite dimensional bimodules. 
	\end{corollary}
	By \Cref{bijection with fd equivariant slice modules}, we conclude the following classification of simple objects of $\HC_{\L}(\A_\lambda)$, where $\L$ denotes the minimal relevant leaf of the flower quiver variety. 
	\begin{corollary}
		When $\lambda\in\ZZ$ and $\lambda \ge (n-1)({\ell}-1)+w$, 
		or $\lambda\le -(n-1)({\ell}-1)$, there are precisely $n!$ non-isomorphic simple objects in $\HC_{\L}(\A_\lambda)$. Otherwise, $\HC_{\L}(\A_\lambda)$ is empty. 
	\end{corollary}
	
	\section{Future work}
	In an ongoing work, we aim to generalize the construction of this paper to other relevant leaves. 
	Namely, let $\L$ be a relevant leaf and $\HC_\L(\A)$ be as defined in \Cref{definition of cats of HC bimodules with support condition}.
	Let $\sliceA$ be the quantized slice quiver algebra associated to $\L$. 
	We expect to construct a restriction functor 
	\begin{equation}
		\HC_{{\L}}(\A)_F \to \HC^{\pi_1(\L)}_{0}(\sliceA)_F 
	\end{equation} 
	similar to the constructions in \Cref{section Minimally-supported HC bimodules}. Here $F$ indicates we consider the category of Harish-Chandra bimodules equipped with a good filtration. 
	
	The first step will be to compute the fundamental groups $\pi_1(\L)$. Then, to construct functors similar to $\F_2$ through $\F_9$, we will need to 
	\begin{enumerate}
		\item construct model varieties $X_K$ for a nice affine open cover $\L_K$ of $\L$, and obtain an analogue of \Cref{the final affine symplectic isomorphism from XI to TstarRI};
		\item pick affine open subsets $\L_I,\L_J$ of $\L$ such that the complement of their union has codimension 2 in $\L$, and obtain an analogue of \Cref{structure of nbhd after reduction};
		\item to carry out the reglueing construction, we need to show $H^2_{\dR}(\L) = 0$. This condition does not hold in general. For example, suppose $\M_0^0$ has no codimension 2 leaf, and let $\L$ be the open leaf of $\M_0^0$. Since the morphism $\rho: \M_0^0 \to \M_0^0$ is semismall, the complement of $\rho\inverse(\L)$ in $\M^\theta_0$ has complex codimension at least 2. Therefore, $H^2_{\dR}(\L )\cong H^2_{\dR}(\M^\theta_0)$. The latter space classifies the quantizations of $\M_0^0$ via the period map, hence must be nonzero. 
	\end{enumerate}
	To avoid the issue in (3), we need to restrict ourselves to a smaller family of leaves. For example, we may require the representation type to have the form
	\begin{equation}
		\tau = (\alpha_\infty ,1; n_1\alpha,1; n_2\alpha,1;\cdots n_k\alpha,1),
	\end{equation}
	i.e. the framing and its dual are 0. On the other hand, Harish-Chandra bimodules whose support is the whole $\M_0^0$ (more generally, Harish-Chandra bimodules with maximal support over any quantized symplectic singularities) are classified in \Cite{losev2021harish}. 
	
	The slice quiver variety to general relevant leaves are still hypertoric varieties, so the combinatorial tools of classifying finite dimensional simple modules over the quantized slice algebra are still available.
	
	As remarked in \Cref{remark cod 4 condition can fail}, the codimension condition
	\Cref{boundary of minimal leaf has codim 4} for the minimal leaf no longer holds for other relevant leaves. Therefore, the functor $\F_3$ defined in \Cref{subsection functor F3} may not be an equivalence. We expect $\F_3$, as well as the restriction functor, to be a fully faithful embedding. It is not yet clear to us what the image should be. 
	
	\printbibliography

@article{nakajima1998quiver,
  title={Quiver varieties and Kac-Moody algebras},
  author={Nakajima, Hiraku},
  journal={Duke Math J.},
  volume={91},
  number={3},
  pages = {515-560},
  year={1998}}

@article{beauville2000symplectic,
title={Symplectic singularities},
author={Beauville, Arnaud},
journal={Inventiones mathematicae},
volume={139},
number={3},
pages={541--549},
year={2000},
publisher={Springer}
}

@article{kaledin2006symplectic,
title={Symplectic singularities from the Poisson point of view},
author={Kaledin, Dmitry},
year={2006},
publisher={Walter de Gruyter}
}

@article{bellamy2021symplectic,
title={Symplectic resolutions of quiver varieties},
author={Bellamy, Gwyn and Schedler, Travis},
journal={Selecta Mathematica},
volume={27},
number={3},
pages={1--50},
year={2021},
publisher={Springer}
}

@article{nakajima1994instantons,
title={Instantons on ALE spaces, quiver varieties, and Kac-Moody algebras},
author={Nakajima, Hiraku},
journal={Duke Mathematical Journal},
volume={76},
number={2},
pages={365--416},
year={1994},
publisher={Duke University Press}
}

@article{bozec2016Quiverswithloops,
title={Quivers with loops and generalized crystals},
volume={152},
journal={Compositio Mathematica},
publisher={London Mathematical Society},
author={Bozec, Tristan},
year={2016},
pages={1999–2040}}

@article{bezrukavnikov2021etingof,
title={Etingof’s conjecture for quantized quiver varieties},
author={Bezrukavnikov, Roman and Losev, Ivan},
journal={Inventiones mathematicae},
volume={223},
number={3},
pages={1097--1226},
year={2021},
publisher={Springer}
}

@article{crawley2001geometry,
title={Geometry of the moment map for representations of quivers},
author={Crawley-Boevey, William},
journal={Compositio Mathematica},
volume={126},
number={3},
pages={257--293},
year={2001},
shorthand={CB01},
publisher={London Mathematical Society}
}

@article{losev2006symplectic,
title={Symplectic slices for actions of reductive groups},
author={Losev, Ivan},
journal={Sbornik: Mathematics},
volume={197},
number={2},
pages={213--224},
year={2006},
publisher={Turpion Ltd}
}

@article{losev2021harish,
title={Harish-Chandra bimodules over quantized symplectic singularities},
author={Losev, Ivan},
journal={Transformation Groups},
volume={26},
number={2},
pages={565--600},
year={2021},
publisher={Springer}
}

@article{losev2012completions,
title={Completions of symplectic reflection algebras},
author={Losev, Ivan},
journal={Selecta Mathematica},
volume={18},
number={1},
pages={179--251},
year={2012},
publisher={Springer}
}

@article{kaledin2000symplectic,
title={Symplectic resolutions: deformations and birational maps},
author={Kaledin, Dmitry},
journal={arXiv preprint math/0012008},
year={2000}
}

@article{losev2011finite,
title={Finite-dimensional representations of W-algebras},
author={Losev, Ivan},
journal={Duke Mathematical Journal},
volume={159},
number={1},
pages={99--143},
year={2011},
publisher={Duke University Press}
}

@book{chriss1997representation,
title={Representation theory and complex geometry},
author={Chriss, Neil and Ginzburg, Victor},
volume={42},
year={1997},
publisher={Springer}
}

@article{braden2012hypertoric,
title={Hypertoric category O},
author={Braden, Tom and Licata, Anthony and Proudfoot, Nicholas and Webster, Ben},
journal={Advances in Mathematics},
volume={231},
number={3-4},
pages={1487--1545},
year={2012},
publisher={Elsevier},
shorthand={BLPW12}
}

@book{hartshorne1977algebraic,
title={Algebraic geometry},
author={Hartshorne, Robin},
volume={52},
year={1977},
publisher={Springer Science \& Business Media}
}

@article{grothendieck1968sga,
title={SGA 2},
author={Grothendieck, Alexander},
journal={S{\'e}minaire de G{\'e}om{\'e}trie Alg{\'e}brique du Bois Marie-1962-Cohomologie locale des faisceaux coh{\'e}rents et th{\'e}oremes de Lefschetz locaux et globaux (North-Holland, Amsterdam)},
year={1968}
}

@book{bott1982differential,
title={Differential forms in algebraic topology},
author={Bott, Raoul and Tu, Loring W},
volume={82},
year={1982},
publisher={Springer}
}

@article{losev2012isomorphisms,
title={Isomorphisms of quantizations via quantization of resolutions},
author={Losev, Ivan},
journal={Advances in Mathematics},
volume={231},
number={3-4},
pages={1216--1270},
year={2012},
publisher={Elsevier}
}

@article{gordon2014category,
title={On category $\mathcal{O}$ for cyclotomic rational Cherednik algebras},
author={Gordon, Iain G and Losev, Ivan},
journal={Journal of the European Mathematical Society},
volume={16},
number={5},
pages={1017--1079},
year={2014}
}
	
\end{document}